\documentclass[11pt]{article}  
\usepackage[latin1]{inputenc}
\usepackage{amsmath,indentfirst,qsymbols,array}
\usepackage{graphicx,psfrag,amsfonts,bbm}
\usepackage[usenames]{color}

\DeclareMathOperator{\Numb}{Numb}

\DeclareMathOperator{\Lip}{Lip}
\DeclareMathOperator{\ODE}{ODE}
\DeclareMathOperator{\SOL}{SOL}
\DeclareMathOperator{\Time}{Time}

\DeclareMathOperator{\arcsinh}{arcsinh}

\DeclareMathOperator\NAVMAX{Navmax}
\DeclareMathOperator\MAXBALL{Maxball}
\DeclareMathOperator\HL{HL}
\DeclareMathOperator\Eq{Eq}
\DeclareMathOperator\sol{sol}

\def \Yao{{\sf Yao}}
\def \Squ{{\sf Squ}}

\def \BB{{\sf BB}}
\def \Nb{{\sf Nb}}

\def \1{{\bf 1}}%\textrm{\dsrom{1}}}
\def \Cost{{\sf Cost}}

\def \X{{\bf X}}
\DeclareMathOperator\Nav{\X}

\def \Path{{\sf Path}}
\def \Poi{{\sf Poisson}}
\def \Pos{{\sf Pos}}
\def \Grid{{\sf Grid}}

\def \app#1#2#3#4#5{\begin{array}{rccl} #1:&#2&\longrightarrow&#3\\ &#4&\longmapsto&#5\end{array}}
\def \as{\xrightarrow[n]{a.s.}}

\def \bD{{\bf D}}

\def \bY{{\bf Y}}
\def \bZ{{\bf Z}}
\def \bC{{\bf C}}
\def \bar{\overline}

\def \ben{\begin{eqnarray}}
\def \be{\begin{eqnarray*}}

\def \bq{\begin{equation}}

\def \build#1#2#3{\mathrel{\mathop{\kern 0pt#1}\limits_{#2}^{#3}}}

\def \cD{{\cal D}}
\def \captionn#1{\begin{center}\begin{minipage}{15cm}\sf\caption{\small #1}\end{minipage}\end{center}}

\def \cro#1{\llbracket#1\rrbracket}

\def \cross{{\sf Cross}}

\def \dis{\displaystyle}
\def \een{\end{eqnarray}}
\def \ee{\end{eqnarray*}}

\def \eq{\end{equation}}
\def \eref#1{(\ref{#1})}
\def \floor#1{\lfloor#1\rfloor}

\def \lip{\Lip^\star}

\def \l{\left}

\def \ov#1{\overrightarrow{#1}}

\def \r{\right}

\def \sous#1#2{\mathrel{\mathop{\kern 0pt#1}\limits_{#2}}}
\def \sur#1#2{\mathrel{\mathop{\kern 0pt#1}\limits^{#2}}}

\def \sv{{\sf v}}

\def\cq{$\hfill \square$}

\def\prooff#1{\noindent{\bf Proof of #1. \rm}}
\def\proof{\noindent{\bf Proof. }}

\def\B{B}
\def\C{C}
\def\Nu{w}

\newcommand{\CC}{{\cal C}}

\newcommand{\bK}{{\bf K}}
\newcommand{\bP}{{\bf P}}
\newcommand{\bS}{{\bf S}}

\newcommand{\ind}{\mathbbm{1}}
\newqsymbol{`E}{\mathbb{E}}
\newqsymbol{`I}{\textrm{\dsrom{1}}}
\newqsymbol{`N}{\mathbb{N}}
\newqsymbol{`O}{\Omega}
\newqsymbol{`P}{\mathbb{P}}
\newqsymbol{`Q}{\mathbb{Q}}
\newqsymbol{`R}{\mathbb{R}}
\newqsymbol{`Z}{\mathbb{Z}}
\newqsymbol{`a}{\alpha}
\newqsymbol{`e}{\varepsilon}
\newqsymbol{`o}{\omega}
\newqsymbol{`t}{\tau}

\RequirePackage{ifpdf}
\ifpdf%
   \let\igcommande=\includegraphics%
   \renewcommand\includegraphics[2][]{%
   \IfFileExists{#2.psfrag}{
       \immediate\write18{pdffrag -l "\the\linewidth" -t "\the\textheight" -w "\the\textwidth" -p "#1" -f "#2"}%
       \igcommande{#2}% les paramtres doivent tre limins
   }{\igcommande[#1]{#2}}}% pas de modification
\fi

\graphicspath{{../DESSINS/}}

\def \sect{{\sf Sect}}
\def \tri{{\sf Tri}}
\def \cam{{\sf Cam}}
\def \ovl#1#2{{\tri}[#1,#2]}

\def \dovl#1{\overrightarrow{\tri}#1}
\def \ovc#1#2{{\cam}[#1,#2]}
\def \dovc#1{\overrightarrow{\cam}#1}
\def \ov#1#2{{\sf Sect}[#1,#2]}
\def \dov#1{\overrightarrow{\sect}#1}

\def \Cbis#1{{\bf C_{bis}^{#1}}}
\def \Cbor#1{{\bf C_{bor}^{#1}}}
\def \Qbor#1{{\bf Q_{bor}^{#1}}}
\def \Qbis#1{{\bf Q_{bis}^{#1}}}
\def \bQ{{\bf Q}}

\def \ST{{\bf ST}}
\def \SY{{\bf SY}}
\def \RNY{{\bf RNY}}
\def \RNT{{\bf RNT}}

\def \Y{{\bf CY}}
\def \T{{\bf CT}}
\def \DT{{\bf DT}}
\def \DY{{\bf DY}}

\def \YN{$p_\theta$ Yao navigation}

\def \TN{$p_\theta$ \T~navigation}

\def \STN{$\theta$ Straight T~navigation}
\def \SYN{$\theta$ Straight Yao navigation}
\begin{document}
\renewcommand{\baselinestretch}{1.15}

\newtheorem{lem}{Lemma}
\newtheorem{defi}[lem]{Definition}
\newtheorem{pro}[lem]{Proposition}
\newtheorem{theo}[lem]{Theorem}
\newtheorem{cor}[lem]{Corollary}
\newtheorem{rem}[lem]{Remark\rm}{\rm}
\newtheorem{com}[lem]{Comments\rm}{\rm}
\newtheorem{exa}[lem]{Examples \rm}{\rm}
\sloppy
\begin{center}

\LARGE{\bf Asymptotic of geometrical navigation \\
on a random set of points of the plane}\medskip\normalsize
\[\begin{array}{ll}
\textrm{\Large Nicolas Bonichon}& \textrm{\Large ~~~~~~~~~~Jean-Fran\c{c}ois Marckert}\end{array}\]
\textrm{CNRS, LaBRI, Universit\'e de Bordeaux}\\
\textrm{351 cours de la Lib\'eration}\\
\textrm{33405 Talence cedex, France}\\
\textrm{email: name@labri.fr}
 \end{center}
\begin{abstract}

A navigation on a set of points $S$ is a rule for choosing which point to move to 
from the present point in order to progress toward a specified target. We study 
some navigations in the plane where $S$ is a non uniform Poisson point process  (in a finite domain) with intensity going to $+\infty$. We show the convergence of the traveller path lengths, the number of stages done, and the geometry of the traveller trajectories, uniformly for all starting points and targets, for several navigations of geometric nature. Other costs are also considered. This leads 
to asymptotic results on the stretch factors of random Yao-graphs and random $\theta$-graphs.
\end{abstract}
\small
The first author is partially supported by the ANR project ALADDIN and
the second author is partially supported by ANR-08-BLAN-0190-04 A3.
\normalsize
\section{Introduction}
\subsection{Navigations}

Often, a traveller who can be a human being, a migratory animal, a letter, a radio message, a message in a wireless ad hoc network, ... wanting to reach a point $t$ starting from a point $s$ has to stop along the route where, according to the case, he can sleep, eat, be sorted, be amplified, or routed. Generally the traveller can not stop everywhere: only some special places offer what is needed (a hostel, a river, a post-office, a radio relay station, a router, ...). Often the traveller can not compute the optimal route from its initial position: it has to choose the next point to move to using only some local information. This paper deals with this problem: a random set $S$ of possible stops given, what happens if a traveller stops ``at the first point in $S$'' which is in the direction of $t$ up to an angle $\theta$? How many steps are done? What is the total distance done? 
In this paper we answer these questions in the asymptotic case, when the number of points in $S$ goes to $+\infty$.\par
 
Formally, consider a traveller on $`R^2$ beginning its travel at the \emph{starting position} $s$ and wanting to reach the \emph{target} $t$ using as \emph{set of possible stopping places} $S$, a finite subset of $`R^2$. We call navigation\footnote{also called in the literature  \emph{memoryless routing algorithm}} with set of stopping places $S$, a function $\Nav$ from $`R^2\times `R^2$ onto $`R^2$ such that for any $(s,t)$, $\Nav(s,t)$ belongs to $S\cup\{t\}$, and satisfies moreover $\Nav(s,s)=s$ for any $s\in `R^2$. The position $\X(s,t)$ corresponds to the first stop of the traveller in its travel from $s$ to $t$. Hence, 
\[\Nav(s,t,j):=\Nav(\Nav(s,t,j-1),t),\qquad j\geq 1\] 
are the successive stops of the traveller, where by convention $\Nav(s,t,j)=s$ for $j= 0$. When no confusion is possible on $s$, $t$ and $\Nav$, we will write $s_j$ instead of $\Nav(s,t,j)$. The quantity 
\[\Delta^{\Nav}(s,t,j):=\Nav(s,t,j)-\Nav(s,t,j-1)\]
is called the \emph{$j$th stage}. For a general navigation algorithm $\Nav$, if $\#S<+\infty$, either $\Nav(s,t,j)=t$ for $j$ large enough, or $\Nav(s,t,j)\neq t$ for all $j$. In the first case, the global navigation from $s$ to $t$ succeeds, whereas in the second case, it fails.  In case of success, the  (global) path from $s$ to $t$ is
\begin{equation}\label{eq:road}
\Path^{\Nav}(s,t):=\l(s_0,\dots,s_{\Nb^{\Nav}(s,t)}\r)
\end{equation}
where $\Nb^{\Nav}(s,t):=\min \l\{j,~ s_j=t\r\}$ is the number of stages needed to go from $s$ to $t$.\par
We are interested in navigations in $`R^2$ where the point to move to is chosen according to some rules of geometrical nature: we consider two classes of so-called \emph{compass navigations}; these navigations select the next stopping place to move to as the ``nearest'' point of $s$ in the set $S\cup\{t\}$, in the ``direction'' of $t$ (see Section~\ref{ssec:TTN}).\par
All along the paper, $\cD$ refers to a bounded and simply connected open domain in $`R^2$. 
The sets of considered stopping places $S$ are finite random subsets of $\cD$ taken according to two models: $S$ will be either the set $\{p_1,\dots, p_n\}$ where the points $p_i$'s are picked independently according to a distribution having a regular density $f$ (with respect to the Lebesgue measure) in $\cD$ (see Section \ref{sssec:niidpositions}),  or $S$ will be a Poisson point process with intensity $nf$, for some $n>0$ (see Section~\ref{ssec:mrp}). \par   

The main goal of this paper is to study the global asymptotic behaviour of the paths of the traveller. Global means all the possible trajectories corresponding to all starting points $s$ and targets $t$ of $\cD$ are considered all together. Several quantities are then studied, describing the ``deviation'' of the paths of the traveller (or functionals of the path, as the length) to a deterministic limit (see Section~\ref{ssec:QI}). The asymptotic is made on the number of points of $S$, which will go to $+\infty$ (that is $n\to+\infty$ in one of the models). 

\paragraph{Convention}
Throughout the paper, the two dimensional real plane $`R^2$ is identified with the set of complex numbers $\mathbb{C}$ and according to what appears simpler, the complex notation or the real one is used without warning. The real part, the imaginary part, the modulus of $z$ are respectively written $\Re(z)$, $\Im(z)$, and $|z|$; the argument  $\arg(z)$ of any real number $z\neq 0$ is the unique real number $\theta$ in $[0,2\pi)$ such that $z=\rho e^{i\theta}$, for some $\rho> 0$ (we set $\arg(0)=0$). The characteristic function of the set $A$ is denoted by $\1_A$. Notation $\cro{x,y}$ refers to the set of integers included in $[x,y]$ and $B(x,r)=\{ y \in \mathbb{C},~|x-y|<r\}$ to the open ball with centre $x$ and radius $r$ in $\mathbb{C}$, the closed ball is $\bar{B}(x,r)$. For $x\in\mathbb{C}$, $A\subset \mathbb{C}$, $d(x,A)=\inf\{|x-y|,y\in A\}$.

\subsection{Two types of navigations}

\label{ssec:TTN}
The two types of navigation introduced below, namely \emph{cross navigation} and \emph{straight navigation}, may appear very similar, but their asymptotic behaviours as well as their analysis are quite different.

\begin{figure}[ht]
\psfrag{h}{$h$}
\psfrag{b}{$\beta/2$}
\centerline{\includegraphics[height=2.4cm]{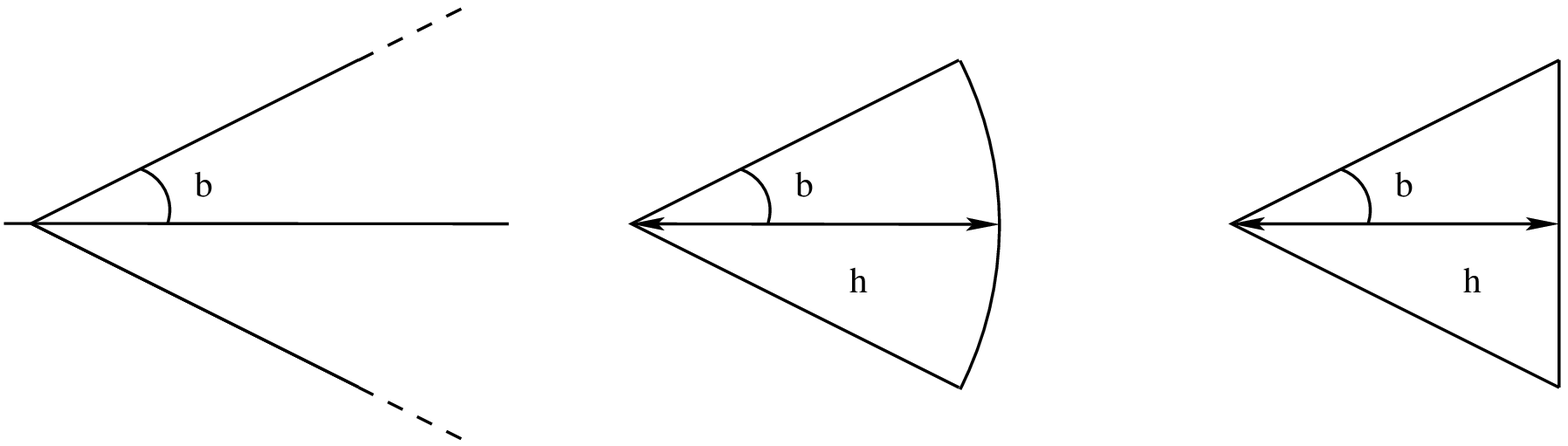}}
\captionn{\label{fig:def-ang-2}Representation of $\sect(\beta), \tri(\beta)(h)$ and $\cam(\beta)(h)$. Notice that each peak point of the figures is excluded from the corresponding set.}
\end{figure}

For any $\beta\in(0,2\pi)$, as illustrated on Fig.~\ref{fig:def-ang-2}, let 
\ben 
\nonumber\sect(\beta)&:=&\l\{\rho e^{i\nu} ,~ \rho> 0, |\nu|\leq \beta/2\r\},\\
\label{cam}\cam(\beta)(h)&:=&\sect(\beta)\cap \bar{B}(0,h),\\
\nonumber\tri(\beta)(h)&:=& \l\{x+iy ,~x\in (0,h], y\in \mathbb{R}, |y| \leq x\tan(\beta/2)\r\}.
\een 
Notation $\cam$ and $\tri$ are short version for ``Camembert portion'' and ``triangle''. 
\vspace{-6 mm}

\subsubsection{First type of navigations: Cross navigations}
\label{sssec:FTN}

\begin{figure}[ht]
\psfrag{x}{$x$}
\psfrag{h}{$h$}
\psfrag{h1}{$\HL_1$}
\psfrag{h0}{$\HL_0$}
\psfrag{h2}{$\HL_2$}
\psfrag{S1}{$\ov0x$}
\psfrag{S1h}{$\ovl0x(h)$}
\psfrag{S2h}{$\ovc{0}{x}(h)$}
\psfrag{S2}{$\ov{1}{x}$}
\psfrag{S3}{$\ov{2}{x}$}
\psfrag{S4}{$\ov3x$}
\psfrag{S5}{$\ov4x$}
\psfrag{S6}{$\ov5x$}
\centerline{\includegraphics[height=3.4cm]{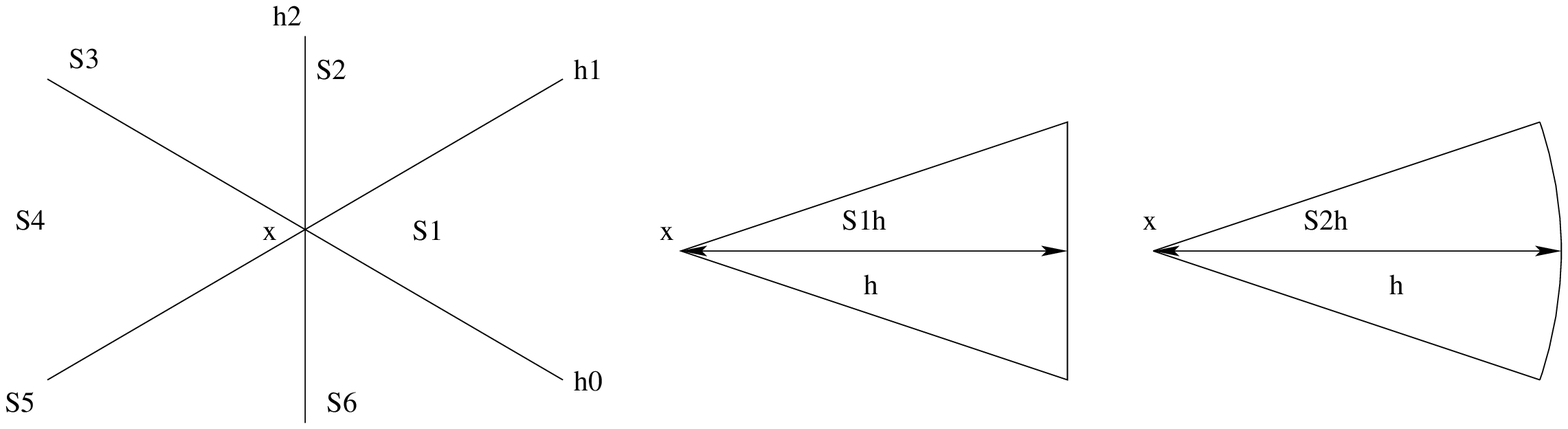}}
\captionn{\label{fig:def-ang}Representation of the 6 angular sectors around a point $x$ with $p_{\theta}=6$; on the second picture, representation of $\ovl0x(h)$, and on the third $\ovc0x(h)$.}
\end{figure}

Cross navigations are parametrised by a parameter $\theta$ satisfying $\theta=2\pi/p_{\theta}$ for some $p_{\theta} \in\{3,4,\dots\}$. For any $x\in \mathbb{C}$, for $\kappa$ in$\cro{0,p_{\theta}-1}$ the \emph{$\kappa$th angular sector} around $x$ is
\[\ov{\kappa}{x}:= x+e^{i\kappa\theta}\sect(\theta).\]
The two half-lines $\HL_\kappa(x)$ and $\HL_{\kappa+1}(x)$ defined by
\begin{equation}
\label{eq:HL}
\HL_j(x):=x+\l\{\rho e^{i\theta(j-1/2)},~\rho>0\r\},~~~ j\in \cro{0,p_\theta},
\end{equation}
are called the \emph{first} and \emph{last border} of  $\ov{\kappa}{x}$. As illustrated in Fig.~\ref{fig:def-ang}, for $h>0$, the $\kappa$th  triangle and $\kappa$th Camembert section around $z$ with height $h$ are respectively: 
\be
\ovl{\kappa}{x}(h)&:=&x+e^{i\kappa\theta}\tri(\theta)(h),~~~~
\ovc{\kappa}{x}(h):=x+e^{i\kappa\theta}\cam(\theta)(h),
\ee
where for any $z_1,z_2 \in\mathbb{C}$, and any  $A \subset \mathbb{C}$, $z_1+z_2A$ is the set $\{z_1+z_2y,~y\in A\}$. \par
As one can see on Fig.~\ref{fig:def-ang}, the half-lines $(\HL_\kappa(t),\kappa\in\cro{1,p_\theta})$, forms a cross around $t$ that we denote by $\cross(t):=\bigcup_{\kappa=1}^{p_\theta} \HL_\kappa(t).$ This justifies the terminology ``cross navigation''.\par

Let $s$ and $t$ in $\cD$, and let $\kappa$ be such that $t \in \ov{\kappa}{s}$. Consider ${\bf B}_\kappa(s)=\{s+re^{i\kappa\theta},~r\geq 0\}$ the bisecting half-line of $\ov{\kappa}{s}$. Consider the lines parallel to $\HL_\kappa(s)$ and $\HL_{\kappa+1}(s)$ passing via $t$. These two lines intersect the half-line ${\bf B}_\kappa(s)$ in two points. 
\begin{figure}[ht]
\psfrag{delta}{${\bf B}_\kappa(s)$}
\psfrag{hk}{$\HL_\kappa(s)$}
\psfrag{hk+1}{$\HL_{\kappa+1}(s)$}
\psfrag{s}{$s$}
\psfrag{t}{$t$}
\psfrag{alpha}{$\alpha$}
\psfrag{theta/2}{$\theta/2$}
\psfrag{ist}{$I(s,t)$}
\centerline{\includegraphics[height=3.6cm]{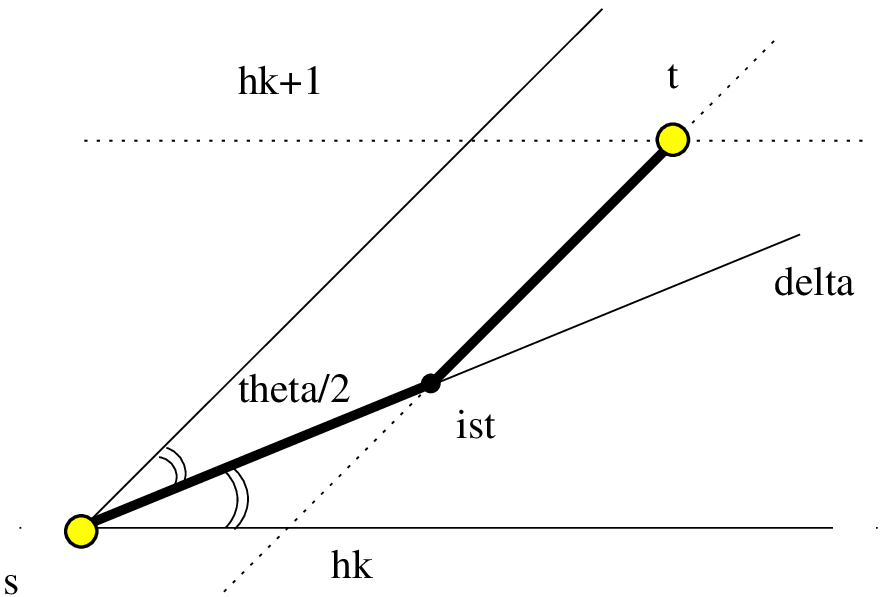}}
\captionn{\label{fig:ist}Representation $I(s,t)$, ${\bf B}_\kappa(s)$. $\Gamma(s,t)$ is represented by the heavy segments.}
\end{figure}
The point which is the closest from $s$ is called $I(s,t)$, as represented on Fig.~\ref{fig:ist}. The compact  $\Gamma(s,t)$ is:
\begin{equation}
\Gamma(s,t):=[s,I(s,t)]\cup[I(s,t),t].
\end{equation}

\YN\ and \TN, noted $\Y$ and $\T$ are defined given a set of stopping places $S$, a finite subset of $`R^2$.

\subsubsection*{Definition of \YN}

\begin{figure}[ht]
\centerline{\includegraphics[width=7cm,angle=22.5,origin=c]{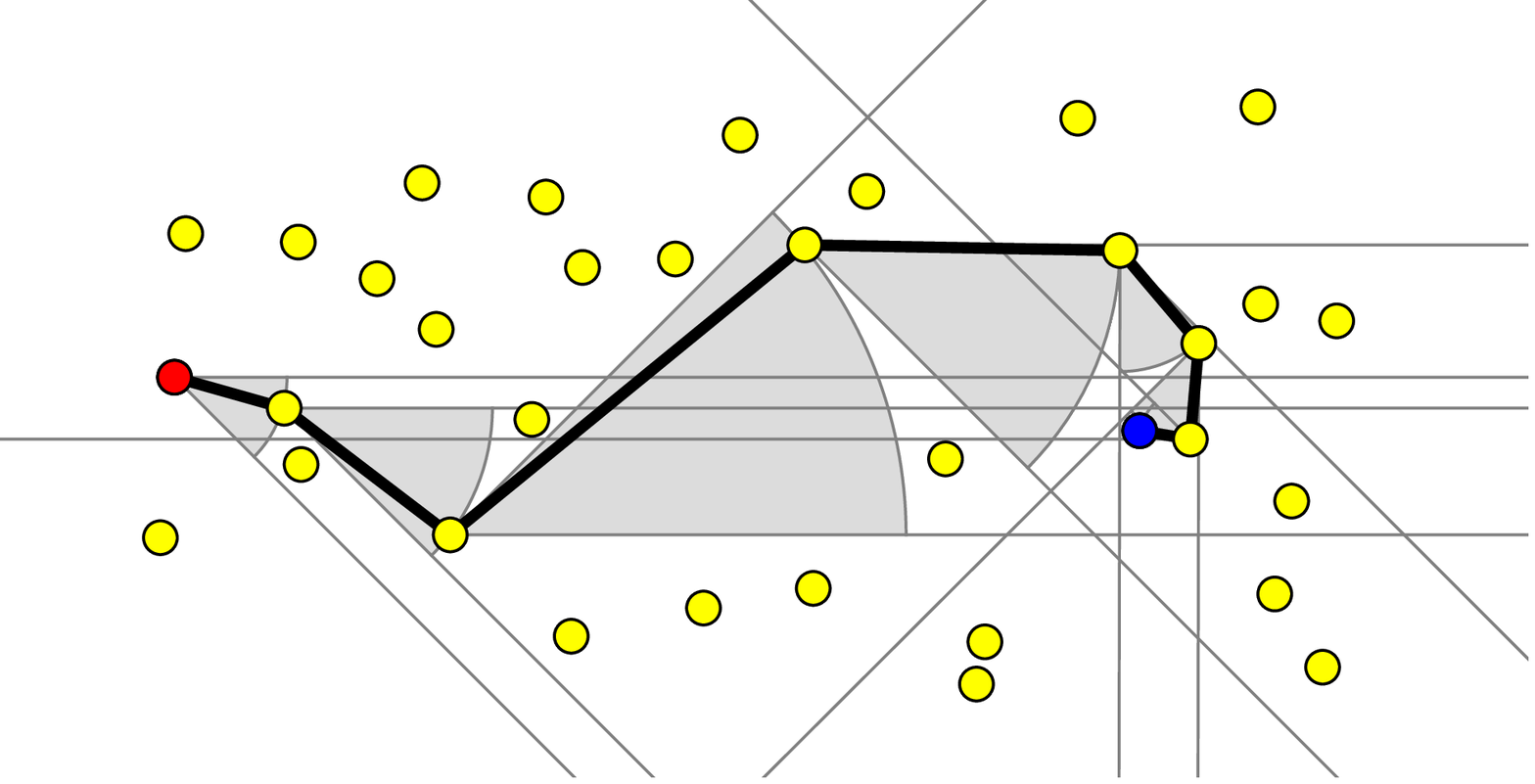} \includegraphics[width=7cm]{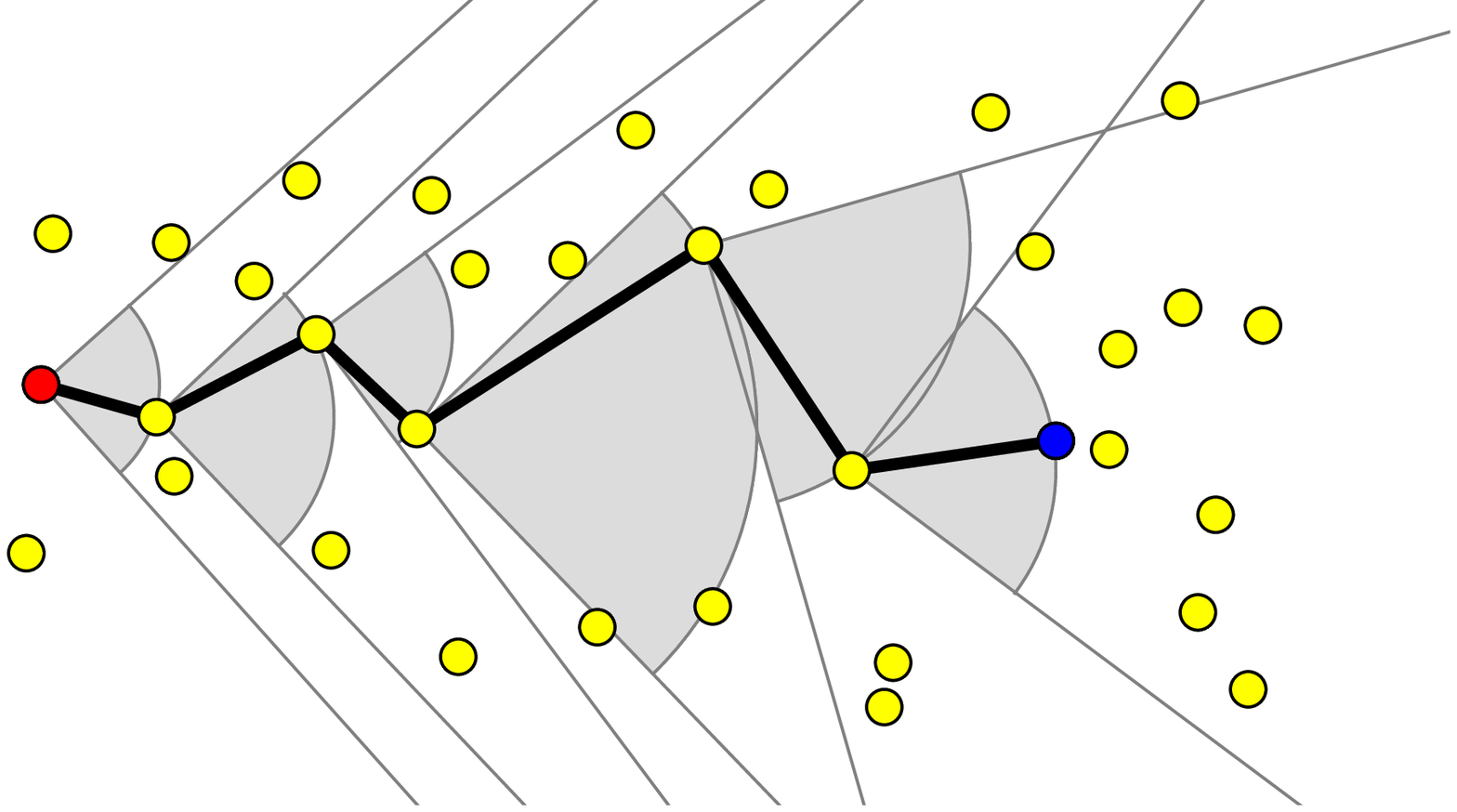}}
\vskip -1cm 
\captionn{\label{fig:example4}(a) An example of Yao navigation with $p_\theta = 8$ the starting point is on the left. Are drawn the decision domains for each stage, and the sectors that contain the destination. (b)  An example of Yao straight navigation with $\theta = \pi/2$ the starting point is on the left. Are drawn the decision domains for each stage, and the sectors that contain the destination.}
\end{figure}

\noindent For $(s,t)\in \mathbb{`R}^2\times  \mathbb{`R}^2$, consider the smallest integer $\kappa$ in $\cro{0,p_\theta-1}$ such that $t$ lies in $\ov{\kappa}{s}$.  Consider the smallest $r$ such that $\ovc{\kappa}{s}(r)\cap \l(S\cup\{t\}\r)$ is not empty. 
We set $\Y(s,t)=z$, the element of $\ovc{\kappa}{s}(r)\cap \l(S\cup\{t\}\r)$ that is the closest of the first border of  $\ov{\kappa}{s}$  (see Fig.~\ref{fig:example4} (a)).\par

This navigation appears to be the canonical navigation in Yao's graphs~\cite{Yao82}.

\subsubsection*{Definition of \TN}

\noindent \TN~is a simple variant of \YN. It is defined as \YN, except that $ \ovl{\kappa}{s}(r)$ replaces  $\ovc{\kappa}{s}(r)$. 

\begin{com}If $t$ belongs to $\HL_\kappa(s)$ then $t$ lies in both $\ov{\kappa}{s}$ and $\ov{\kappa-1 \mod p_\theta}{s}$; in the other cases, $t$ belongs to a unique sector $\ov{\kappa}{s}$. In the sequel, the domains $\ovc{\kappa}{s}(r)$ and  $\ovl{\kappa}{s}(r)$ which are, depending on the algorithm considered, the minimum domains allowing to decide to which point to move to, are called the \emph{decision domains}.\par
\end{com}

\subsubsection{Second type of navigations: Straight navigation}

Straight navigations are closer to real traveller strategies: the stopping places are chosen not too far from the segment leading to the target. We set two straight navigations, parametrised by a real number $\theta\in(0,2\pi)$. The difference between the following objects and those before defined is that the sectors are oriented such that their bisecting lines are the straight line $(s,t)$: 
\be
\dov{(s,t)}&:=&s+\exp(i\arg(t-s))\,\sect(\theta),\\
\dovc{(s,t)}(h)&:=& s+\exp(i\arg(t-s))\,\cam(\theta)(h),\\
\dovl{(s,t)}(h)&:=&s+\exp(i\arg(t-s))\,\tri(\theta)(h).
\ee
\SYN\ and \STN\ will be denoted $\SY$ and $\ST$. 

\subsubsection*{\SYN}
\noindent  For $(s,t)\in S^2$, consider the smallest $r$  such that the set $\dovc{(s,t)}(r)\cap (S\cup \{t\})$ is not empty. Set $\SY(s,t)=z$, where $z$  is the closest of the first border of  $\dovc{(s,t)}$ (see Fig.~\ref{fig:example4} (b)).

\subsubsection*{\STN} 
\noindent The \STN~$\ST$ is similarly defined except that $ \dovl{(s,t)}(r)$ replaces  $\dovc{(s,t)}(r)$.

\subsection{The model of random stopping places set: a Poisson point process}
\label{ssec:mrp}

Denote by $\lip(\cD)$ the set of Lipschitz functions $f:\cD \to `R^{+}$ with a positive infimum on $\cD$: a function $f$ is in $\lip(\cD)$, if there exists $\alpha_f>0$, such that $|f(x)-f(y)|\leq \alpha_f |x-y|$ for all $x,y\in \cD$, and  $m_f:=\inf \{f(x), x\in \cD\}>0$. Since $\cD$ is bounded, $M_f:=\max\{f(x), x\in \cD\}<+\infty$. \par
In this paper, we consider a Poisson point processes $\bS(f)$ whose intensity measure $\mu_f$ on $\cD$ is $\mu_f(A)=\int_A f(z)dz$, for some $f\in \lip(\cD)$ (in words, $\mu_f$ has density $f$ with respect to the Lebesgue measure on $`R^2$): the necessity of the Lipschitz condition will clearly appear later on, with the appearance of a differential equation. The distribution of $\bS(f)$ is denoted by $`P_f$. Since $\bS(f)$ has no multiple points a.s., $\#\bS(f)\sim \Poi(\int_\cD f(u) du)$. For any disjoint Borelian subsets $B_1,\dots,B_k$ of $`R^2$, $(\bS(f)\cap B_1),\dots, (\bS(f)\cap B_k)$ are independent and $\#\l(\bS(f)\cap B_j\r)\sim \Poi\l(\int_{B_j}f(u)\ du\r)$ (for more explanations on PPP, see Chap. 10 in \cite{KAL}). As usual, a representation of the set $\bS(f)$ is as follows: $\bS(f)=\{x_1,\dots,x_{\bf n}\}$, where ${\bf n}\sim\Poi(\int_{\cD}f)$, and where the  $x_1,x_2, \dots, $ are i.i.d. and independent of ${\bf n}$,  and each $x_i$ has density $f/\inf_{\cD}f$. \par

In the paper, we assume the measure $\sum_{x \in \bS(f)} \delta_x$ defined on a probability space $(\Omega,{\cal A},`P)$, and is considered as a random variable taking its values in ${\sf N}$ the set of counting measures in $\mathbb{C}$, equipped with the $\sigma$-field $\mathcal{N}$ generated by the sets $E_{B,k}=\{\mu, \mu(B)=k\}$ for compacts sets $B\subset \mathbb{C}$, $k$ integer.

The model consisting of $n$ i.i.d. points chosen according to a density $f$ is considered in Section~\ref{sssec:niidpositions}.

\begin{com} The navigations  $\Nav$ we study, as well as the function $\Path^{\Nav}$, $\Pos^{\Nav}$, and various cost functions are all defined given a stopping places set $S$, and then, there are some functions of $S$ (one should write $\Nav(S)$, etc. but we choose to delete this $S$ to lighten the notation). On $(\Omega,{\cal A},`P)$, these functions are random variables with values in some functional spaces.
\end{com}

\subsection{Quantities of interest}
\label{ssec:QI}

In~\eref{eq:road} is defined  $\Path^\X(s,t)$. In many applications, a quantity of interest is the comparison between the Euclidean distance $|s-t|$ and the \emph{path length}, the total distance done by the traveller in case of success:
\[|\Path^\X(s,t)|:=\sum_{j=1}^{\Nb^\X(s,t)}|s_j-s_{j-1}|.\]
The associated \emph{trajectory} is the compact subset of $`R^2$ formed by the union of the segments $[s_j,s_{j+1}]$:
\[[\Path^\X(s,t)]:=\bigcup_{j=1}^{\Nb^\X(s,t)}[s_{j-1},s_j].\]
 One of the results of the paper is the comparison between $[\Path^\X(s,t)]$ and a limiting object.
 \par
If the point process is not homogeneous, the evaluation of the traveller's trajectory calls for a precise study of the traveller's speed. We introduce the function $\Pos_{s,t}^{\Nav}:[0,\Nb^\X(s,t)]\to`R^2$ whose values coincides with those of  $\Nav(s,t,.)$ at integer points, and interpolated in between. Hence, $\Pos_{s,t}^\X$ gives the position of the traveller according to the number of stages; for any $s,t\in\cD$, $[\Path^\X(s,t)]$ equals the range of $\Pos_{s,t}^\X$. From our point of view, the asymptotic results obtained for $\Pos_{s,t}$ are among the main contributions of this paper.\par

In some applications,  the sum of a function of the stage lengths appears to be the relevant quantity instead of the length. Formally a ``unitary cost function'' $H:\mathbb{C}\to `R$ is considered. The total cost associated with $H$ corresponding to $\Path^\X(s,j)=\l(s_0,\dots,s_{\Nb^\X(s,t)}\r)$ is 
\begin{equation}
\Cost_H^\X(s,t):=\sum_{j=1}^{\Nb^\X(s,t)} H\l(s_j-s_{j-1}\r).
\end{equation}
Important elementary cost functions are $H:x\to |x|$ which gives $\Cost_H(s,t)=|\Path(s,t)|$, the function $H:x\to 1$ which gives $\Cost_H^\X(s,t)=\Nb^\X(s,t)$, and functions of the type $H(x)=|x|^\beta$ which corresponds, at least for some $\beta$, to some real applications. 

\section{Statements of the main results}
\label{sec:main-results}
In the rest of the paper,  $f$ is a fixed element of $\lip(\cD)$ and $a$ a fixed positive number. The set $\cD[a]:= \l\{x \in\cD~: d(x,\complement \cD)\geq a\r\}$ is the set of points in $\cD$ at distance at least $a$ to $\complement\cD$, the complement of $\cD$. 
The asymptotic behaviour of $\Path^{\Nav}(s,t)$ (when the set of stopping places $S$ increases) is difficult to predict if the limiting trajectories between $s$ and $t$ are too close of $\complement\cD$. We then set a notation to designate the pairs $(s,t)$ that can be treated. Set 
$$\cD'[a]:=\{(s,t)\in \cD^2 ~:  Z(s,t)\subset\cD[a]\},$$
where according to the context $Z(s,t)=\Gamma(s,t)$ (in case cross navigations are concerned) or $Z(s,t)=[s,t]$ (in case straight ones are concerned). 
Some other restrictions concerning $\theta$ will be added in order to avoid situations where the navigation can fail.\par 
 Before giving the results, we define some constants, computed in Sections~\ref{sssec:DNUUI} and~\ref{ssec:CYN}, and related to the speed of the traveller along its trajectory, and to some ratio ``mean length of a stage'' divided by ``length of the projection of a stage'' with respect to some ad hoc directions.
 \renewcommand{\arraystretch}{1.9}

\begin{equation}
\begin{array}{|l|l|}
\hline
\Cbis{\ST}=\Cbis{\T}=\frac{1}2\sqrt{\frac{\pi}{\tan(\theta/2)}}&\Cbis{\SY}=\Cbis{\Y}=\frac{\sqrt{2\pi}\sin(\theta/2)}{\theta^{3/2}}\\
\hline
\Cbor{\T}=\sqrt{\frac{\pi\cos^3(\theta/2)}{4\,\sin(\theta/2)}}&\Cbor{\Y}=\frac{\sqrt{\pi/2}\sin(\theta)}{\theta^{3/2}}\\\hline
\Qbis{\ST}=\Qbis{\T}=\frac{1}{2}\left( \frac{1}{\cos(\theta/2)}+\frac{\arcsinh(\tan(\theta/2))}{\tan(\theta/2)} \right)&
\Qbis{\SY}=\Qbis{\Y}=\frac{\theta/2}{\sin(\theta/2)}\\\hline
\Qbor{\T}= \frac{1}{2}\left( \frac{1}{\cos^2(\theta/2)}+\frac{\arcsinh(\tan(\theta/2) ) }{\sin(\theta/2)} \right)&
\Qbor{\Y}=\frac{\theta}{\sin(\theta)}
\\
\hline
\end{array}
\end{equation}
\renewcommand{\arraystretch}{1.2}

\subsection{Limit theorems for the straight navigations}
\label{ssec:SN}

The first theorem  uniformly compares the path length  $|\Path^\X(s,t)|$ with a multiple of $|s-t|$.
\begin{theo} Let $\X=\SY$ and $\theta < \pi/2$, or  $\X=\ST$ and $\theta \leq \pi/2$.
\label{theo:main1}
For any $\alpha\in(0,1/8)$, any $\beta>0$, for $n$ large enough
\[`P_{nf}\l(\sup_{(s,t)\in \cD'[a]}  
\big||\Path^{\X}(s,t)|-\Qbis{\X}\times|s-t|\big|\geq n^{-\alpha}\r)\leq n^{-\beta}. \] 
\end{theo}
The terms $\Qbis{\ST}$ and $\Qbis{\SY}$ ``measures'' the efficiency of the traveller with respect to the direction of the bisecting lines of the decision sectors used. \par

We now describe the asymptotic behavior of $\Pos_{s,t}^\X$. In the case of straight navigation, the limiting position function $\Pos^{\infty,c}$ is the deterministic solution of a differential equation, and depends on a real parameter $c$ function of $\theta$ and of $\X$. \par
For any $(\lambda,\nu)\in (0,+\infty)\times[0,2\pi]$, let $F_{\lambda,\nu}$ be the function from  $\cD$ into $\mathbb{C}$ defined by
\begin{equation}
\label{eq:F-l-m}
F_{\lambda,\nu}(z):={\lambda e^{i\nu}}/{\sqrt{f(z)}}.
\end{equation}
For $s_0\in \cD$ given, consider the following ordinary differential equation 
\begin{equation}
\label{eq:ODE}
\ODE(\lambda,\nu,s_0):=\left\{\begin{array}{l}
                               \rho(0)=s_0,\\
\dis\frac{\partial \rho(x)}{\partial x}=F_{\lambda,\nu}(\rho(x)).
\end{array}\right.
\end{equation}
By Cauchy-Lipschitz Theorem, $\ODE(\lambda,\nu,s_0)$ admits a unique solution $x\mapsto \SOL_{s_0}^{\lambda,\nu}(x)$ valid for $x\geq 0$ and while $\SOL_{s_0}^{\lambda,\nu}(x)$ stays inside $\cD$; indeed since $f$ is in $\lip(\cD)$, so do $F_{\lambda,\nu}$. Since $f$ takes its values in $`R^+$, $\arg(F_{\lambda,\nu})=\nu$. Hence, $\SOL_{s_0}^{\lambda,\nu}$ appears to be a monotone parametrisation of the segment $\{s_0+xe^{i\nu},~ x>0\}\cap \cD$ (more precisely,  of the component containing $s_0$ in this intersection), the speed along this segment at position $z$ being $F_{\lambda,\nu}(z)$. Now, choose two points $(s,t)$ in $\cD$ such that $[s,t]\subset \cD$ and consider $\SOL_s^{\lambda, \arg(t-s)}$ for some fixed $\lambda>0$. In this case, the range of  $\SOL_{s}^{\lambda,\arg(t-s)}$ contains $[s,t]$ and by continuity there exists a time $\Time_{s,t}^\lambda$ such that
\[\SOL_{s}^{\lambda,\arg(t-s)}\l(\Time_{s,t}^\lambda\r)=t.\]
We then define the function $\Pos^{\infty,\lambda}_{s,t}$ as the function $\SOL_{s}^{\lambda,\arg(t-s)}$ stopped at $t$.
\begin{theo}\label{theo:main2}
Let  $\X=\SY$ and $\theta < \pi/2$, or $\X=\ST$ and $\theta \leq \pi/2$.
For any $\alpha\in(0,1/8)$, any $\beta>0$ for $n$ large enough
\[`P_{nf}\l(\sup_{(s,t)\in \cD'[a]} \sup_{x\geq 0} \l|\Pos^{\X}_{s,t}(x\sqrt{n})-\Pos^{\infty,\Cbis{\X}}_{s,t}(x)\r| \geq n^{-\alpha}\r)\leq n^{-\beta}.\]
\end{theo}
A corollary of this, is that the asymptotic traveller trajectories are segments: for $d_H$ the Hausdorff distance between compact subset of $\mathbb{`R}^2$, $ `P_{nf}\l(\sup_{(s,t)\in \cD'[a]} d_H([\Path(s,t),[s,t]) \geq n^{-\alpha}\r)\leq n^{-\beta}$.\par
The asymptotic behaviour of different cost functions are studied in Section~\ref{ssec:aco}.

\subsection{Limit theorems for cross navigation}
\label{ssec:NN}

For any $c_1>0,c_2>0$, and any $s,t \in \mathbb{C}$, let 
\begin{equation}
\bD_{c_1,c_2}(s,t):=c_1 \times |I(s,t)-s|+c_2\times|t-I(s,t)|
\end{equation}
be a kind of weighted length of $\Gamma(s,t)$.\par
\begin{figure}
\centerline{\includegraphics[height=4.2cm,angle=-30]{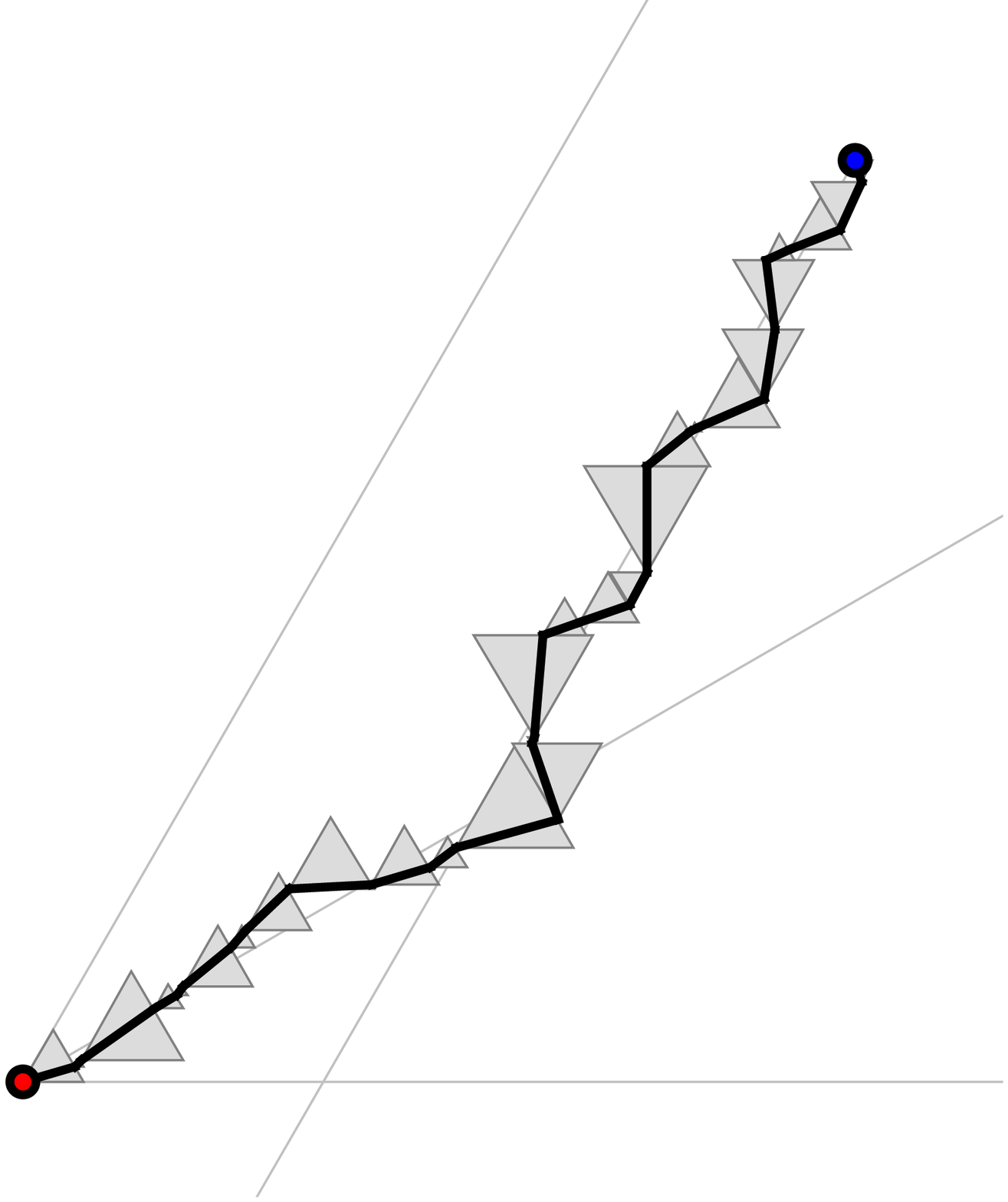} \includegraphics[height=4.2cm,angle=-30]{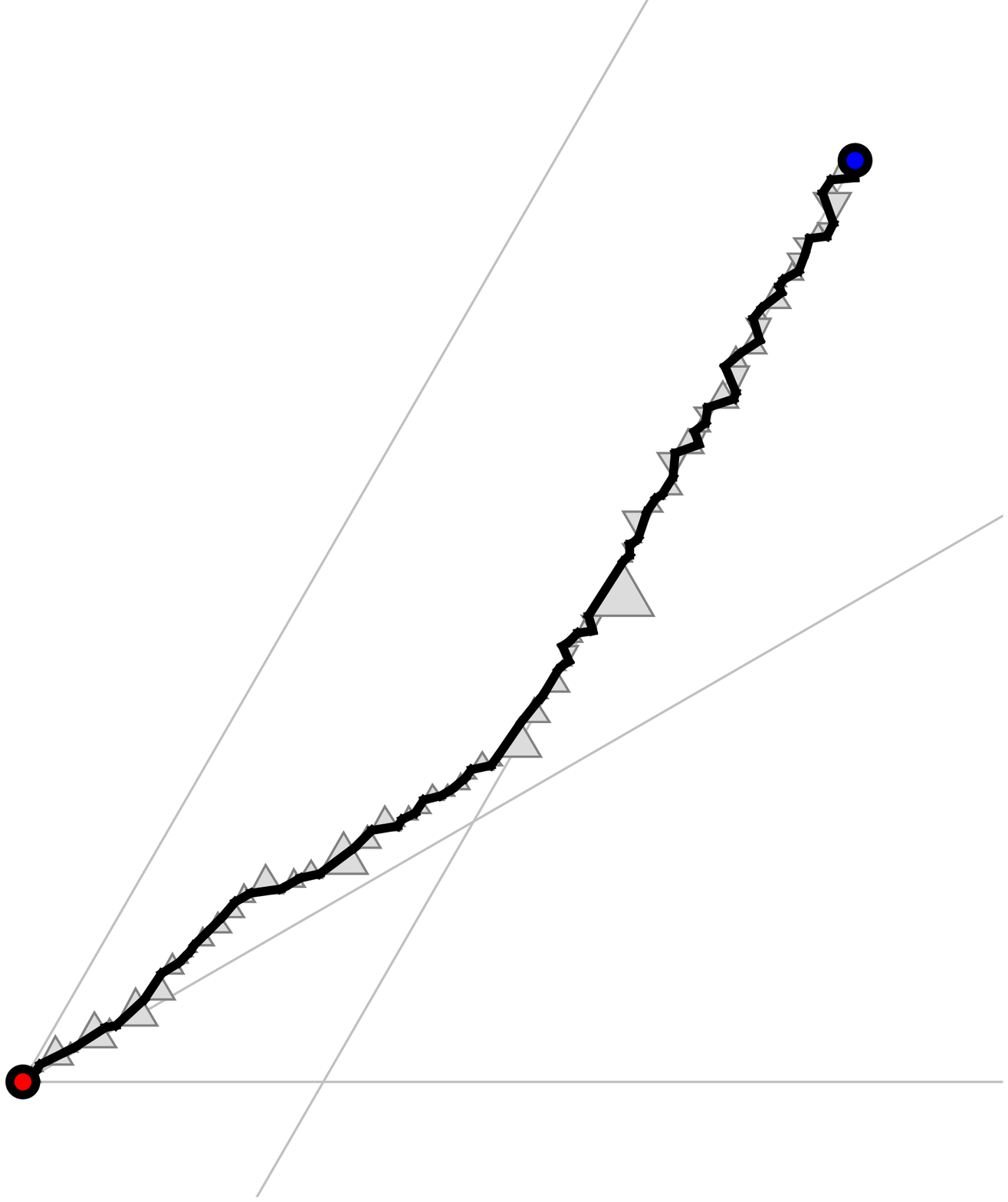} \includegraphics[height=4.2cm,angle=-30]{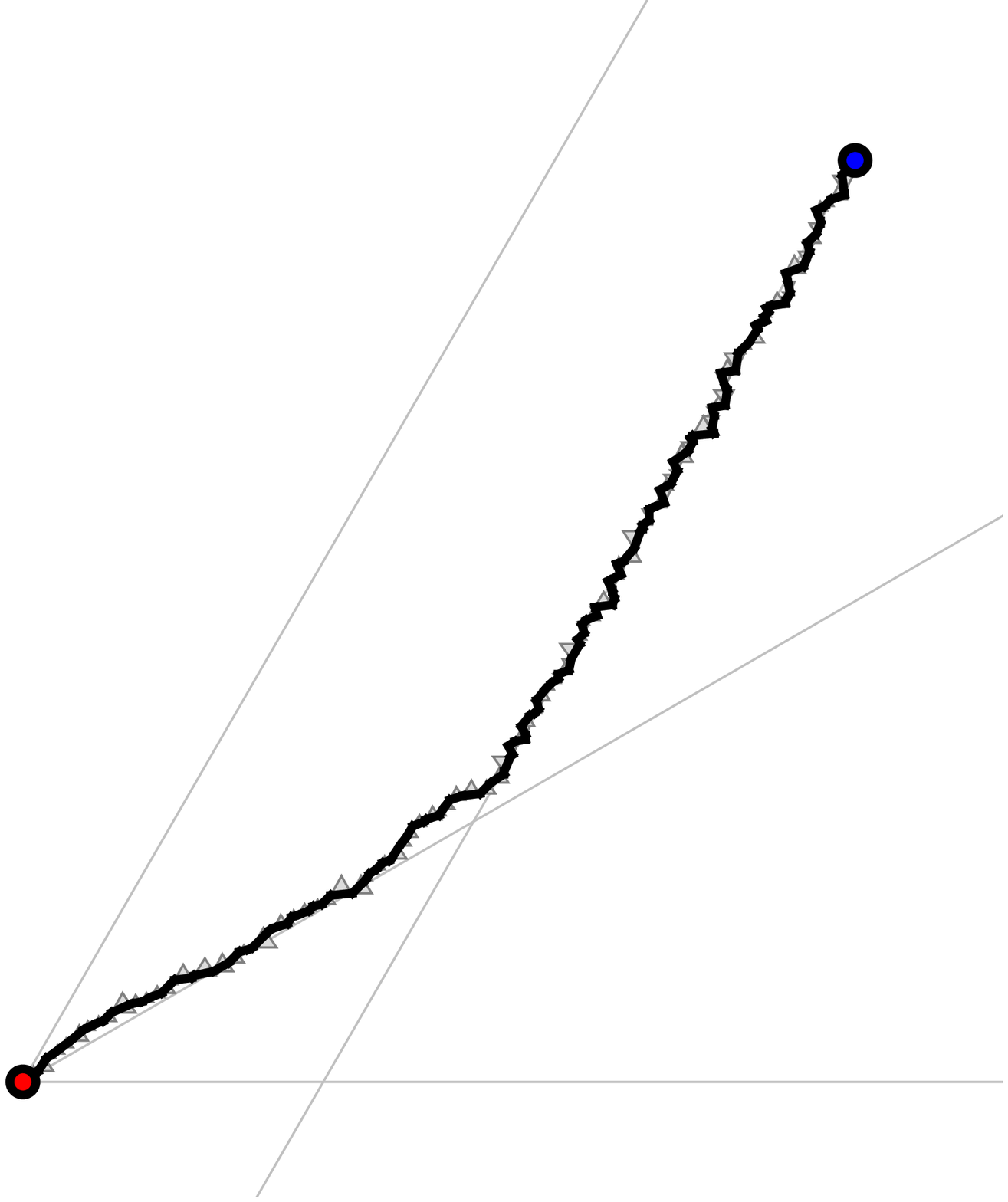}}
\captionn{\label{fig:examplebig}Examples of \TN\ (for $p_\theta=6$) with 500 nodes, 5000 nodes and 20000 nodes uniformly distributed in the unit square.}
\end{figure}

The ``cross navigations'' analogous of Theorem~\ref{theo:main1} is Theorem~\ref{theo:main3} below (see Fig.~\ref{fig:examplebig}). 
\begin{theo}Let $\X\in\{\T,\Y\}$ and $\theta\leq \pi/3$.
\label{theo:main3}
 For any $\alpha\in(0,1/8)$, any $\beta>0$, for $n$ large enough
\[`P_{nf}\l(\sup_{(s,t)\in \cD'[a]}  
\l||\Path^{\X}(s,t)|-\bD_{\Qbis{\X},\Qbor{\X}}(s,t)\r|\geq n^{-\alpha}\r)\leq n^{-\beta}. \] 
\end{theo}
The terms $\Qbor{\ST}$ and $\Qbor{\SY}$ measure some local efficiency  of the traveller with respect to the direction of the border of the decision sectors. \par 
As said above, $\Y$ is the canonical navigation on Yao's graph, and $\T$ is also the canonical navigation on the so-called $\Theta$ graphs \cite{Clarkson87,keil88, RS91a}.  A worst case study have been made showing that  the stretch factor of theses navigations is at most $\frac{1}{1-2\sin{\pi/p_\theta}}$, for every $p_\theta > 6$.  For $p_\theta \leq 6$ the stretch factor of these navigations can be unbounded. However, it has been proved that  $\Theta_6$ is a 2-spanner~\cite{BGHI10}, $\Yao_6$ is a $c$-spanner (for some $c$)~\cite{OR10} and $\Yao_4$ is a $8(29+23\sqrt{2})$-spanner~\cite{BDDRSS10}. The theorem obtained here says that ``for a typical set of points'', the navigation distance, and then also the graph distance, between any two points $s$ and $t$ is smaller than $\bD_{\Qbis{\X},\Qbor{\X}}(s,t)+n^{-\alpha}$ with huge probability. For far away points $s$ and $t$, this is much smaller than the worst case bounds.

Again, for $\X\in\{\T,\Y\}$, the function $x\mapsto \Pos_{s,t}^{\X}(.\sqrt{n})$ admits a deterministic limit $\Pos^{\infty,a,b}_{s,t}$ depending on two parameters $a$ and $b$ (which depend on $\X$ and $\theta$) related to the speed of the traveller along the two branches of $\Gamma(s,t)$. Set 
\begin{equation}
\label{eq:cbis-cbor}
\Pos^{\infty,\lambda,\lambda'}_{s,t}(x):= 
\left\{
\begin{array}{ll}
\SOL_s^{\lambda,\arg(I(s,t)-s)}(x) &\textrm{ for } x \in[0, \Time_{s,I(s,t)}^{\lambda}],\\
\SOL_{I(s,t)}^{\lambda',\arg(t-I(s,t))}(x-\Time_{s,I(s,t)}^{\lambda(1)}) &\textrm{ if } x \in\Time_{s,I(s,t)}^{\lambda}+[0,\Time_{I(s,t),t}^{\lambda'}],\\
t&\textrm{ if }x \geq \Time_{s,I(s,t)}^{\lambda}+\Time_{I(s,t),t}^{\lambda'}.
\end{array} \right.
\end{equation}
\begin{theo}\label{theo:main4}
Let $\X \in \{\Y,\T\}$ and $\theta\leq \pi/3$, for any $\alpha\in(0,1/8)$, any $\beta >0$, any $n$ large enough
\[`P_{nf}\l(\sup_{(s,t)\in \cD'[a]} \sup_{x\geq 0} \l|\Pos^{\X}_{s,t}(x\sqrt{n})-\Pos^{\infty,\Cbis{\X},\Cbor{\X}}_{s,t}(x)\r| \geq n^{-\alpha}\r)\leq n^{-\beta}.\]
\end{theo} As for the straight case, this entails the convergence for the Hausdorff distance of $\Path(s,t)$ to $\Gamma(s,t)$.\par
Again, other results concerning other cost functions are studied in Section~\ref{ssec:aco}.

\subsection{Extensions}

The following extensions can be treated with the material available in the present paper. We just provide the main lines of their analysis.

\subsubsection{Random north navigation}

This is a version of the cross navigation where each point $s\in S$ has its own (random) north $n(s)$ used to compute $\cross(s)$. The random variables $(n(s),s\in S)$ are assumed to be i.i.d. and takes their values in $[0,2\pi]$. This can be used to model some imprecisions in the Yao's construction, where the north is not exactly known\footnote{Each distribution of $n(s)$ leads to a different behaviour for the traveller. If the value of $n(s)$ is 0 a.s., then this models coincides with the standard cross navigation. If the distribution of $n(s)$ owns several atoms, then the speed of the traveller  may be different along numerous directions.} by the points of $S$. The corresponding navigation algorithm is defined as follows: A traveller at $s$, choose the smallest integer $\kappa$ in $\cro{0,p_{\theta}-1}$ such that $t$ lies in $s+e^{in(s)+ik\theta}\sect(\theta)$ and consider the smallest $r$ such that $s+e^{in(s)+ik\theta}\tri(\theta)(r)\cap \l(S\cup\{t\}\r)$ is not empty. Then set $\RNT(s,t)=z$, the element of this set with smallest argument. \par
Using the Camembert sections instead leads to the random north Yao's navigation, $\RNY$.\medskip
The uniform random north navigation is not much different to the straight one, the limiting path being segments, and the limiting path length being a multiple of the Euclidean distance. Using the same arguments than those given in the present paper, one can prove that the distance done by the traveller satisfies, for any $\theta\leq \pi/3$, for any $\alpha\in(0,1/8)$ and $\beta>0$, if $n$ is large enough
\[`P_{nf}\l(\sup_{(s,t)\in \cD'[a]}\l|\Path^{\X}(s,t)-\bQ^{\X}|s-t|\r|\geq n^{-\alpha}\r)\leq n^{-\beta}\]
where 
\be
\bQ^{\RNT}=\frac{\theta/2}{2\sin({\theta}/2)}\left( \frac{1}{\cos(\theta/2)}+\frac{\arcsinh(\tan(\theta/2))}{\tan(\theta/2)} \right),~~~
\bQ^{\RNY}&=&\frac{\theta^2}{2-2\cos(\theta)}.
\ee
The computation of these constants are done in Section \ref{seq:RNM}; the proof of the globalisation of the bounds do not present any problem in this case.

\subsubsection{Model of $n$ i.i.d. positions}
\label{sssec:niidpositions}
Consider $g\in\lip(\cD)$ such that $\ind_{\cD} g= 1$. Let $p_1,\dots,p_n$ be $n$ i.i.d. random points chosen in $\cD$ under the distribution having $g$ as density. Set $\bS^n_g:=\{p_1,\dots,p_n\}$, and by $\bP^n_g$ the distribution of this set.
The analysis of the navigations under $`P_{ng}$ is simpler than under $\bP^n_g$ since under $`P_{ng}$ Markovian properties can be used; the derivation of the results under $\bP^n_g$ will appear to be simple consequences of those on $`P_{ng}$. Indeed, the two models $\bP_g^n$ and $`P_{ng}$ are related via the classical fact
\begin{equation}
\label{eq:cond}
\bP_g^n(~.~)=`P_{ng}(~.~|\, \#\bS=n);
\end{equation}
in other words $\bK_{ng}$ conditioned by $\#\bS(ng)=n$ has the same distribution as $\bK_g^n$.\par
For any measurable event $A$ (element of $\mathcal{N}$ as defined in Section \ref{ssec:mrp}) and any $c>0$,
\[\bP_g^n(A)=`P_{ng}(A~|~ \#\bS=n)\leq \frac{`P_{ng}(A)}{`P_{ng}(\#\bS=n)}.\]
Since $\#\bS(ng)$ is $\Poi(n)$ distributed, one sees that  $`P_{ng}(\#\bS=n)\geq c_1 n^{-1/2}$, for $n$ large enough, for a constant $c_1>0$ (this is an application of the Stirling formula) and therefore
\[\bP_g^n(A)\leq (1/c_1) n^{1/2} `P_{ng}(A),\]
for $n$ large enough. All the results of the present paper under the form $`P_{ng}(A_n)=o(x_n)$ for some decreasing function $x_n$ can be transferred under $\bP_g^n$ to the form $\bP_g^n(A_{n})=O(n^{1/2} `P_{ng}(A_{ng}))=o(n^{1/2}x_{n})$, that is a factor $n^{1/2}$ on the bound must be added. Therefore, the main theorems of the present paper can be transferred to $\bP_g^n$ without any problems since all results (or intermediate results) have the form $`P_{ng}(A_n)=o(x_n)$ for $x_n=o(n^{-1/2})$.

\subsection{Link with other random objects}
\label{ssec:RWC}

\begin{figure}[ht]
\centerline{\includegraphics[height=5cm]{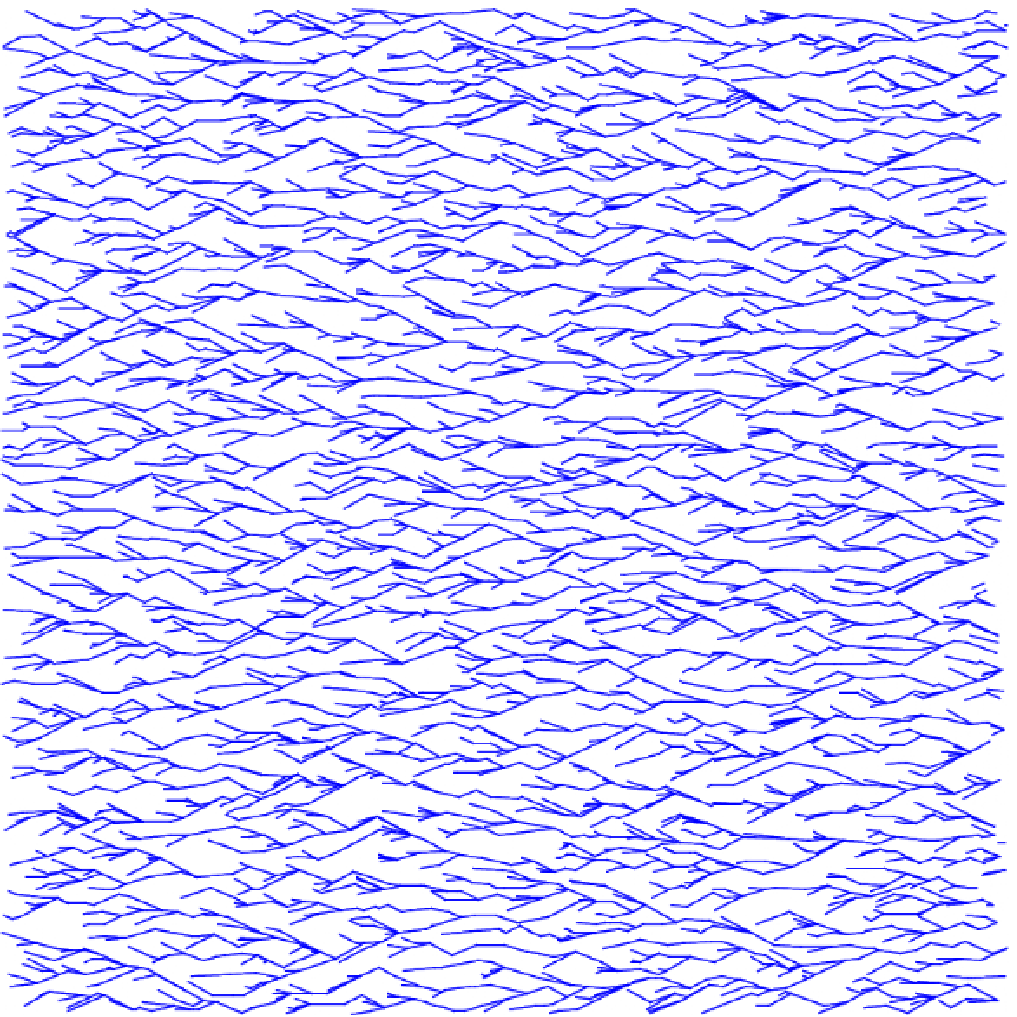} \includegraphics[height=5cm]{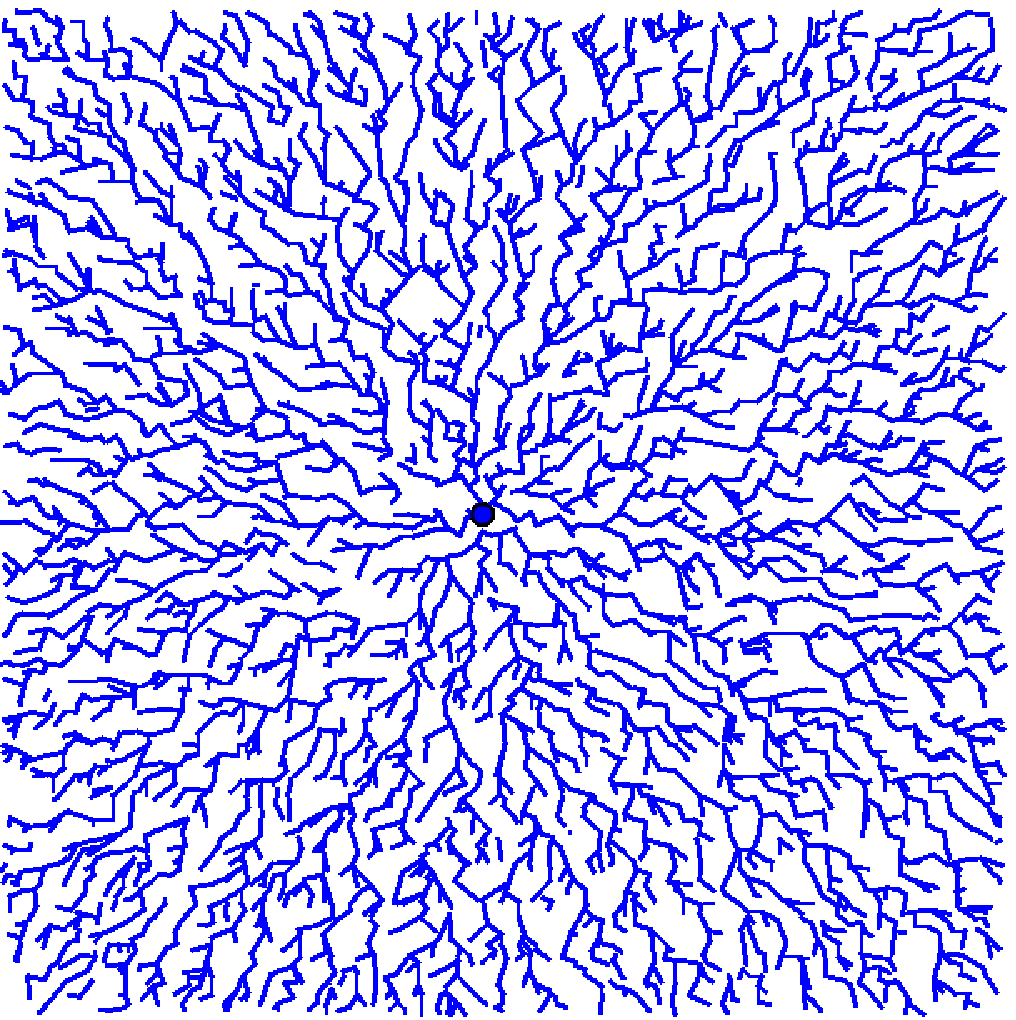} \includegraphics[height=5cm]{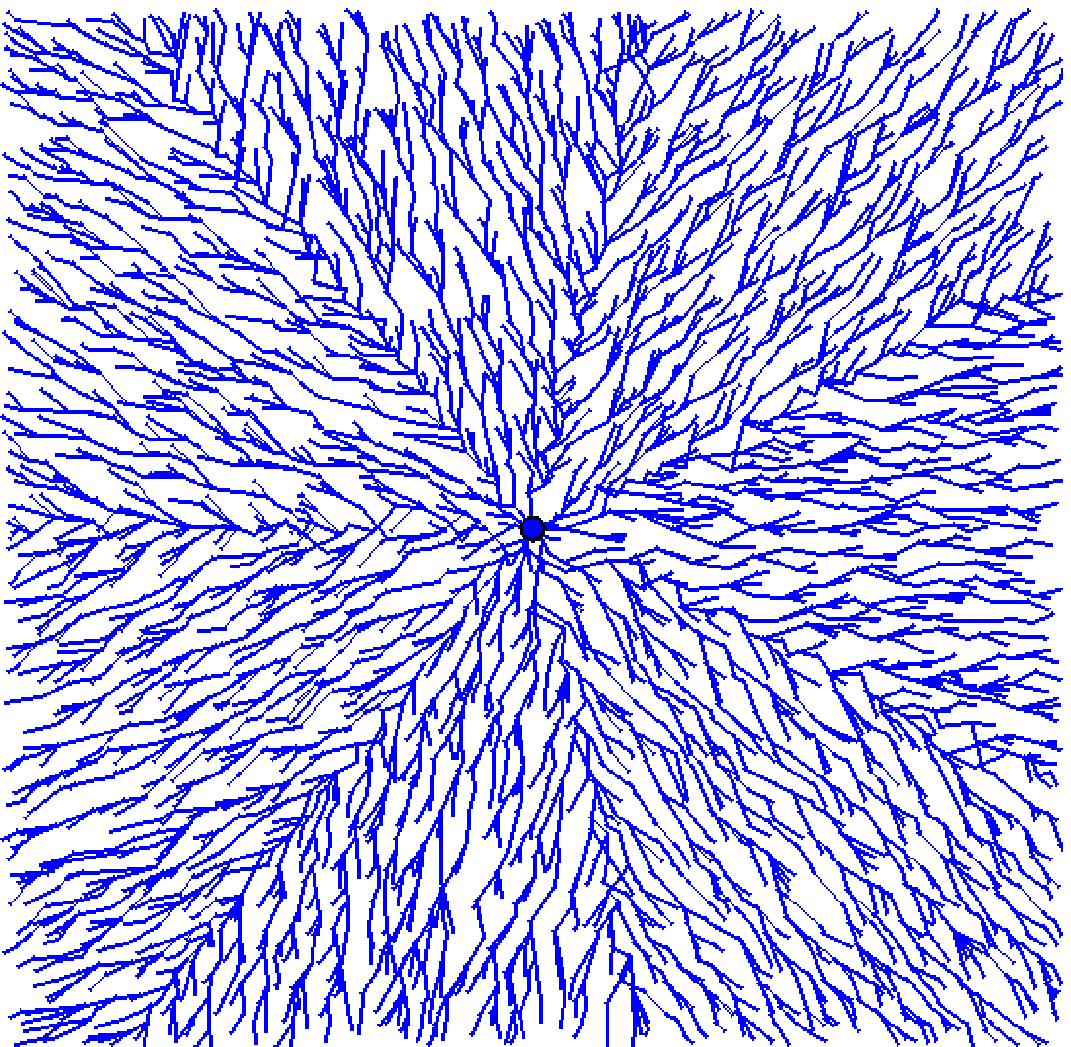}}
\captionn{\label{fig:trees} Union of all paths (toward in a unique target) on 5000 uniform random points for: (a) a \DT\ navigation with $\theta =2\pi/7$. (b) a \SY\ navigation with $\theta =2\pi/3$. (c) a \T\ navigation with $p_\theta =7$.}
\end{figure}

-- In \cite{BB07}, the authors investigate a navigation ${\bf N}$ on a homogeneous PPP $\bS$ on $\mathbb{C}$, where the traveller wants to reach the origin 0 of the plane:  for a given $s$, ${\bf N}(s)$ is the closest element of $s$ in $\bS\cup\{0\}$, which has the additional property to be closer to 0 than $s$. Adding an edge between each $s \in \bS$ and ${\bf N}(s)$ one gets a tree, called the ``radial spanning tree'' of the PPP. \footnote{On the second picture in Fig.\ref{fig:trees} straight \YN\  with $\theta =2\pi/3$; this picture is quite similar to the Fig.1 in   \cite{BB07} even if the navigation is less homogeneous in that paper} The authors provide numerous results concerning this tree; some local properties concerning the degree of the vertices, the length of edges, and some more global properties, as the behaviour of the subset of elements $s$ of $\bS$ such that ${\bf N}(s,k)=0$.  In \cite{BORD}, the author goes on this study. Among other, functional along a path are studied (the tree has infinite number of ends, paths going to infinity).\\
-- In \cite{BR, PW} the authors study the property of a so-called ``minimal directed spanning tree'' (MDST). To each finite subset $S$ of $(0,1]^2$ is associated a tree as follows : $x\in S$ is connected by a directed edge to the nearest $y\in S\cup\{(0,0)\}$, having both coordinates smaller. They then study the asymptotic total length of the MDST when $S$ follows various distributions, with $|S|\to +\infty$.  \rm \\
-- In \cite{ARS}, the authors examine a directed like navigation, where the set of points is a random subset of $\mathbb{Z}^d$: each point of this set is kept with probability $p\in(0,1)$. They then examine the connectivity of the construction according to the dimension. \\ 
-- Another object that may be related to directed navigation is the Brownian web.  In \cite{FFW05}, the authors construct a binary tree using as set of vertices the points of a homogeneous PPP in the plane. The parent-children relation is then induced by a kind of navigation, where, instead of Camembert section, a rectangle section is used; this is deeply similar to what is represented in picture 1 of Fig.\ref{fig:trees}. They then show that in a certain sense, their object converges to the Brownian web (see the references therein for definitions and criteria of convergences to the Brownian web). What is done in the present paper concerning directed navigation let us conjecture that what is done in \cite{FFW05} could be extended to Yao's graphs or $\Theta$ graphs associated to directed navigations. The main structural difference here is that the infinite tree constructed by directed navigation has an unbounded degree, when it was binary in~\cite{FFW05}. \\
-- In a sequence of papers, Aldous and coauthors \cite{AL,AK,AS10}... investigate numerous questions related to navigation (or traveller salesman type problem) in a PPP. In particular, in \cite{AL} sufficient conditions on navigation algorithms are given to have an asymptotic shape (where the shape is roughly, the set of points at distance smaller than $r$ to the origin, properly rescaled).

\section{Toward the proofs, first considerations}
\label{sec:tpfc}

\subsection{Presentation of the analysis}
\label{sec:pa}

We present here the main ideas used in the paper  before entering into the details. 
\begin{itemize}
 \setlength{\itemsep}{2pt}
  \setlength{\parskip}{2pt}
  \setlength{\parsep}{2pt}
\item  Under $`P_{nf}$, a ball of radius $r$ included in $\cD$ contains $O(nr^2)$ point of $\bS(nf)$ in average, and with huge probability less than $n^{1+`e}r^2$, uniformly on all balls (Lemma~\ref{lem:max-boule}). We restrict ourselves to the case where the navigation has the property to force the traveller to come closer to $t$ at each stage. Hence for a right $r=r_n$, Lemma~\ref{lem:max-boule} allows one to show that the contribution of the stages of the traveller in the final ball of radius $r_n$ is negligible. See Section~\ref{ssec:fb}. 
\item To study the behaviour of the traveller far from its target, a local argument is used: under $`P_{nf}$ for a non constant function $f$, the stages $(s_j-s_{j-1},j\geq 1)$ are not identically distributed, neither independent since the value of $s_{j-1}$ affects the distribution of $s_j-s_{j-1}$. But, if one considers only $a_n$ successive stages, with $a_n\to+\infty$ and $a_nn^{`e}/\sqrt{n}\to 0$, these stages stay in a small window around $s_0$ with a huge probability; these stages appear moreover to roughly behave as i.i.d. random variables under the homogeneous PPP $`P_{nf(s_0)}$; the behavior of these $a_n$ stages are seen to depend at the first order, only of $f(s_0)$. Moreover, since $a_n\to+\infty$ some regularizations of the type ``law of large numbers'' occur. See Section~\ref{ssec:Dn}. This local theorem is one of the cornerstones of the study: in some sense, the global trajectory between any two points is a concatenation of these parts of length $a_n$; the successive local values $(f(s_{k\times a_n}),k\geq 0)$ yields directly to a differential equation. See section~\ref{ssec:DE}.
\item The speed of convergence stated in the theorems is roughly given for each part of length $a_n$: it is shown that well chosen deviations are exponential rare, and then the deviation between an entire trajectory and the limiting solution of the ODE is shown to have exponentially rare deviations. Hence, for free, this result maybe extended to all trajectories starting and ending on the two dimensional grid  $n^{-\Nu}\mathbb{Z}^2\cap \cD$, for some $\Nu>0$. The final globalisation consists in the comparison between the paths between any two points $s$ and $t$ of $\cD$, and some well chosen points of the grid  $n^{-\Nu}\mathbb{Z}^2\cap \cD$. This is possible if $\Nu$ is large enough. See Section~\ref{ssec:GB}.
\end{itemize}

\subsection{About the termination of the navigations}
\label{sec:term}

The main of this section is to state and prove the following proposition.
\begin{pro}\label{pro:boule} 
Let $\X=\SY$ and  $\theta \leq 2\pi/3$, or $\X=\ST$ and $\theta \leq \pi/2$, or $\X\in\{\T,\Y\}$ and $\theta\leq \pi/3$. Under $`P_{nf}$, a.s. for all $(s,t)\in\cD'[a]$
\begin{enumerate}\setlength{\itemsep}{2pt}
   \setlength{\parskip}{2pt}
   \setlength{\parsep}{2pt}
\item the navigation from $s$ to $t$ succeeds (i.e. $\exists k$ such that $\X(s,t,k)=t$);
\item the traveller comes closer to the target at each step, $\Path^\X(s,t) \subset B(t,|s-t|)$;
\item $|\Path^\X(s,t)| \leq 2 |s-t| \times \#S\cap B(t,|s-t|)$;
\item if $|s_1-s| \leq |t-s|/2$ then $|t-s| - |t-s_1| \geq (2-\sqrt{3}) |s_1 - s|$.
\end{enumerate}
\end{pro}
\proof
Let $\kappa$ be the integer such that $t \in \ov{\kappa}{s}$ and let $t_{\kappa}$ be the orthogonal projection of $t$ on the bisecting line of the sector $\kappa$.
Let ${\sf Area}$ be $\ovc{\kappa}{s}(|s-t|)$ (resp. $\ovl{\kappa}{s}(|s-t_{\kappa}|)$, $\dovc{(s,t)}(|s-t|)$ and $\dovl{(s,t)}(|s-t|)$) if $\X$ is $\Y$ (resp. $\T$, $\SY$ and $\ST$). For the values of $\theta$ considered, ${\sf Area} \subset B(t,|s-t|)$.  In each case, $s_1 \in {\sf Area}$, and so $|s_1-t| \leq  |s-t|$, which implies the second item of the proposition. The previous inequality is strict except in one particular case for each navigation considered and only for the maximal values of $\theta$ considered: in case ${\sf Area}$ contains only the two points $t$ and $s_1$, and $s_1$ is on the first border of ${\sf Area}$ such that $|s_1-t|=|s-t|$. Observe that in that case, $\arg(s_1-t)-\arg(s-t)$ equals $\pi/2$ if $X=\ST$ and $\pi/3$ for the other navigations (still for the maximal values of $\theta$ considered). Hence the navigation fails only if there exists $(s,t)\in \cD^2$ such that $S\cup \{s\}$ contains 4 points (resp. 6 points) forming an empty square (resp. an empty hexagon) centred at $t$ with no other points of $S$ inside this polygon. Under $`P_{nf}$, almost surely $S$ doesn't contain such a configuration. This implies the first item of the proposition.

The third item comes directly from the two first ones: the length of each stage of $\Path^\X(s,t)$ is at most $2 |s-t|$ and $\Path^\X(s,t)$ is composed of at most $\#S\cap B(t,|s-t|)$ stages.

Now let us prove the last item. Let $\alpha:=\arg(t-s)-\arg(s_1-s)$. For the navigation considered, $\alpha \in [0,\pi/3]$. By the cosine law, we have:
\begin{equation}
\label{sec:alka}
|t-s_1|^2 = |t-s|^2 + |s_1-s|^2 - 2|t-s|.|s_1-s|\cos \alpha.
\end{equation}
Let $x\geq 2$ such that $|t-s|=x|s_1-s| $. By \eref{sec:alka}, $|t-s_1| = |s_1-s| \sqrt{x^2+ 1 - 2x \cos \alpha}$ and then 
$|t-s|-|t-s_1| = |s_1-s| (x - \sqrt{x^2+ 1 - 2x \cos \alpha})$.
For $x\geq 2$, using that $\alpha\mapsto \cos \alpha$ is decreasing on $[0,\pi/3]$,
$x-(x^2+1-x)^{1/2}\leq x - \sqrt{x^2+ 1 - 2x \cos \alpha}\leq x-(x-1)$.Therefore $(x - \sqrt{x^2+ 1 - 2x \cos \alpha}) \in [2-\sqrt{3} , 1]$, we get the last item.~$\Box$

\subsection{A notion of directed navigation}
\label{ssec:Dn}

In the straight navigation under $`P_{nf}$, when the traveller is far from its target, the fluctuations of $\arg(t-s_i)$  stay small for the first values of $i$ (for large $n$). It is then intuitive that the trajectory of the traveller would not be much changed if he'd use the constant direction given by $\arg(t-s_0)$ instead of $\arg(t-s_i)$; similarly, in the cross navigation, when the traveller is far from $\cross(t)$, it is easily seen that along its first stages, the bisecting lines of its decision sectors are parallel to the first one, in other words, he follows the direction of the bisecting line of the sector around $s_0$ containing $t$. In order to make clear these phenomena, the directed navigation $\DT$ is defined below.\par
Again $\theta$ is a fixed parameter chosen in $(0,2\pi)$. 
\begin{defi}$\bullet$ $\DT$ is a map from $\mathbb{C}$ onto $S$ defined as follows. Let $s\in\mathbb{C}$. If $\dov{(s,+\infty)}\cap S$ is empty, then $\DT(s)$ is defined to be $s$. Otherwise consider the smallest $r$ such that $\l(s+\tri(\theta)(r)\r)\cap S$  is not empty.  Then $\DT(s)$ is defined to be the element of this intersection, that is the closest of the first border of $\dov{(s,+\infty)}$.\\
$\bullet$ \DY is defined similarly, except that  $\l(s+\cam(\theta)(r)\r)\cap S$ replaces  $\l(s+\tri(\theta)(r)\r)\cap S$.

\end{defi}
For any $s\in \mathbb{C}$, the successive stopping places $(\DT(s,j),j\geq 0)$ of the traveller satisfy $\DT(s,0)=s$, $\DT(s,j)=\DT(\DT(s,j-1))$. Informally, \DT\ coincides with \ST\ if the target is $t=+\infty$. If the directed navigation is done on a homogeneous PPP on $\mathbb{C}$, the stages
\[\Delta^{\DT}(s,j):=\DT(s,j)-\DT(s,j-1),~~~j\geq 1\]
are i.i.d..
Under $`P_{f}$, for a non constant $f$, the process $(\DT(s,j),j\geq 0)$ is Markovian but the stages are not i.i.d.. (see the asymptotic behaviour of the directed navigation on Fig.~\ref{fig:trees}(a)).

\subsubsection{Directed navigation on a  PPP with constant intensity} 
\label{sssec:DNUUI}

\begin{figure}[ht]
\psfrag{x}{$x_1$}\psfrag{y}{$y_1$}\psfrag{f}{$\xi_1$}\psfrag{g}{$\xi'_1$}\psfrag{y}{$y_1$}\psfrag{p}{$p$}\psfrag{d}{$\Delta_1$}
\psfrag{s}{$s$}\psfrag{n}{$\T(s)$}\psfrag{H0}{$\HL_0(s)$}\psfrag{H1}{$\HL_1(s)$}
\centerline{\includegraphics[width=6cm]{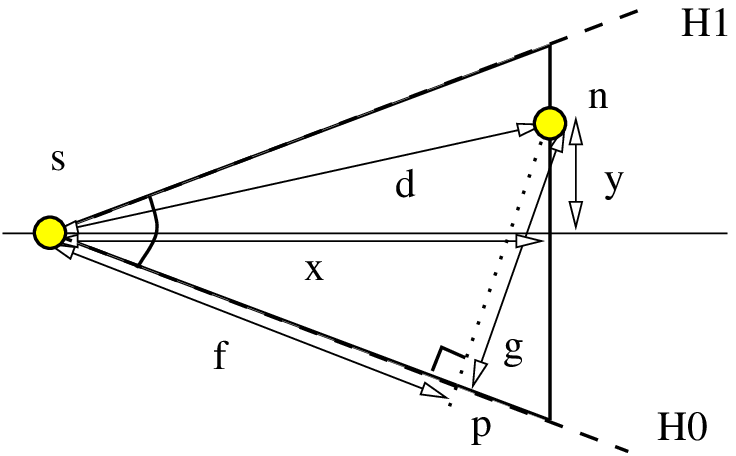}}
\captionn{\label{fig:F} Representation of $\xi_1$ and $\xi'_1$.}
\end{figure}

We consider now \DT~ under $`P_c$ on all $\mathbb{C}$. We write the intensity in index position to let appear the rescaling arguments. Notice further that $\Delta^{\DT}(s,j)$ does not depend on $s$; we then omit it.  We remove now the superscript \DT\ from everywhere.  For any $j\geq 1$, write 
\ben\label{eq:single-stage}
\Delta_c(j)&=&x_c(j)+i\,y_c(j)~~~\textrm{for }x_c(j),y_c(j)\in `R,
\een
and set $l_c(j)=\l|\Delta_c(j)\r|$, the length of this stage. The stages $(\Delta_c(j),j\geq 1)$ are i.i.d., as well as the lengths $(l_c(j),j\geq 1)$ and the pairs $((x_c(j),y_c(j)),j\geq 1)$. By a clear space rescaling argument, we have the following equality in distribution for any $j\geq 1$,
\begin{equation}
\label{eq:Poisson1}
\Delta_c(j)\sur{=}{(d)}\frac{1}{\sqrt{c}}\Delta_1(j),~~~~
(x_c(j),y_c(j))\sur{=}{(d)}\frac{1}{\sqrt{c}}(x_1(j),y_1(j)),~~~~l_c(j)\sur{=}{(d)}\frac{1}{\sqrt{c}}l_1(j).
\end{equation}
Note $x_1$ instead of $x_1(1)$, $y_1$ instead of $y_1(1)$, etc. For any $r>0$, 
\ben\label{eq:zeza}
`P(x_1>r)&=&`P_1(\#\ovl{\kappa}{0}(r)\cap \bS=0)=\exp\l(-r^2\tan(\theta/2)\r),
\een
and then, by integration, one finds $`E(x_1^{\DT})=(1/2)\sqrt{{\pi}/{\tan(\theta/2)}}.$ Notice also that $`E(x_1^{\DT})=\Cbis{\ST}=\Cbis{\T}.$
Now, since $y_1=x_1U$ where $U$ is independent of $x_1$ and uniform on $\l[-\tan(\frac\theta2),\tan(\frac\theta2)\r]$,
\[l_1=|\Delta_1|=\sqrt{x_1{}^2+y_1{}^2}=x_1\l(1+\tan^2(\theta/2)V^2\r)^{1/2}\]
for $V$ uniform in $[0,1]$, independent of $x_1$. A trite computation leads to 
$`E\l(\sqrt{1+\tan^2(\theta/2)V^2}\r)=  \Qbis{\ST},$
as given in the beginning of Section~\ref{sec:main-results}, and then
$`E(l_1^{\DT})=`E(x_1^{\DT})\times \Qbis{\ST}= \Qbis{\ST}\ \Cbis{\ST}.$
We consider another quantity that will play a special role during the analysis of cross navigation.  
Consider the coordinates  $(\xi_1,\xi_1')$ of a stage $\Delta_1$ in the coordinate system $(0,e^{-i\theta/2},e^{-i\theta/2+i\pi/2})$:
\[\Delta_1=\xi_1e^{-i\theta/2}+\xi_1'e^{-i\theta/2+i\pi/2}.\] 
As seen on  Fig.~\ref{fig:F}, $\xi_1$ is the length of the projection of $\Delta_1$ on $\HL_0(s)$.
The middle of the projections of $(x_1,y_1)$ and $(x_1,-y_1)$ coincides with the projection of $(x_1,0)$, and by symmetry of the law of $y_1$,
\begin{equation}
`E(\xi_1^{\DT})=`E(x_1^{\DT}) \cos(\theta/2)=\Cbor{\T}.
\end{equation}
Consider now the distance and the position of the traveller after $k$ stops
\[X_1^{\DT}(k)=\sum_{j=1}^k x_1^{\DT}(j),~~~Y_1^{\DT}(k)=\sum_{j=1}^k y_1^{\DT}(j),~~~L_1^{\DT}(k)=\sum_{j=1}^k l_1^{\DT}(j),~~~\Xi_1^{\DT}(k)=\sum_{j=1}^k \xi_1^{\DT}(j).\]
By the law of large numbers, the following 4-tuple converges a.s.
\be
\l(\frac{X_1^{\DT}(k)}{k}, \frac{Y_1^{\DT}(k)}{k}, \frac{L_1^{\DT}(k)}{k}, \frac{\Xi_1^{\DT}(k)}{k}\r)
&\xrightarrow[k]{a.s.}&\l(\Cbis{\T},\ 0,\ `E(l_1^{\DT}),\ \Cbor{\T}\r)
\ee 
and also then, $\dis\frac{L_1^{\DT}(k)}{X_1^{\DT}(k)}\xrightarrow[k]{a.s.}\frac{`E(l_1^{\DT})}{`E(x_1^{\DT})}=\Qbis{\ST},$ and $\dis\frac{\Xi_1^{\DT}(k)}{X_1^{\DT}(k)}\xrightarrow[k]{a.s.}\frac{`E(l_1^{\DT})}{`E(\xi_1^{\DT})}=\Qbor{\T}.$
These numbers are then the crucial coefficients appearing in Theorems~\ref{theo:main1} and~\ref{theo:main3}. 

\subsubsection{Control of the fluctuations of  $X_1^{\DT}(k), Y_1^{\DT}(k)$ and $L_1^{\DT}(k)$}
\label{sssec:DLN}
First, it is easily checked that the random variables $l_1$, $x_1$, $y_1$, and $\xi_1$ have exponential moments. Indeed, by~\eref{eq:zeza}, $`E(e^{t|x_1|})<+\infty $ for any $t>0$, and the other variables are controlled as follows:  there exist some constants $\beta, \beta',\beta''$ (that depends only on $\theta$) such that
\[0\leq |y_1|\leq \beta x_1, 0\leq  l_1 \leq \beta' x_1; 0\leq \xi_1\leq  \beta'' x_1.\]
Therefore the variables $|x_1-`E(x_1)|$, $|y_1-`E(y_1)|$ and $|l_1-`E(l_1)|$, $|\xi_1-`E(\xi_1)|$ own also exponential moments, which permits to control the fluctuations of $X_1(k), Y_1(k), L_1(k), \Xi_1(k)$ around their means by Petrov's~\cite{PET} Theorem~15~p.52, and Lemma~5~p.54 that we adapt a bit; if $Z_k=z_1+\dots+z_k$, and the $z_i$ are centred, have exponential moments ($`E(e^{t|z_1|})<+\infty$, for some $t>0$), then there exist constants $\alpha>0$ and $c>0$ such that
\bq\label{eq:kez}
\begin{array}{lllll}
`P\l(\l|Z_m \r|\geq x \sqrt{m}\r)&\leq& 2\exp(-c\, x^2)  &\textrm{~~~if }x\leq \alpha\sqrt{m}\\
`P\l(\l|Z_m \r|\geq x \sqrt{m}\r)&\leq& 2\exp(-c\, x\sqrt{m})&\textrm{~~~if }x\geq \alpha\sqrt{m};
\end{array}
\eq
constants $c,\alpha$ depend only on the distribution of $z_1$, but not on $m$ (the same constant $c$ can be chosen for both bounds; we do this). 
A consequence is the following proposition.
\begin{pro}\label{pro:uc}
Let $(x_n)$ be a sequence of real numbers such that $x_n\to +\infty$ and $x_n=O(\sqrt{n})$. There exists a constant $\gamma>0$, such that for any $m$ large enough
\[`P\l(\sup_{k\in\cro{1,m}}\l|\frac{X_1^\DT(k)-k`E(x_1^\DT)}{\sqrt{m}} \r|\geq x_m\r)\leq 2m\exp(-\gamma x_m);\]
The same results hold for $Y_1^\DT(k)-k`E(y_1^\DT)$, $L_1^\DT(k)-k`E(l_1^\DT)$ and $\Xi_1^\DT(k)-k`E(\xi_1^\DT)$ and the same $\gamma$ can be chosen for all the cases.
\end{pro}
\proof~ The proofs for $X_1(k),Y_1(k), L_1(k)$ and $\Xi_1(k)$ are identical. For $X_1(k)$, write
\ben\label{eq:ub}
`P\l(\sup_{k\in\cro{1,m}}\l|\frac{X_1(k)-k`E(x_1)}{\sqrt{m}} \r|\geq x_m\r)&\leq &\sum_{k=1}^{m}`P\l(\l|\frac{X_1(k)-k`E(x_1)}{\sqrt{k}} \r|\geq \frac{m^{1/2}}{\sqrt{k}}x_m\r)
\een
then use the first or second bounds of~\eref{eq:kez}  according to whether $\frac{m^{1/2}x_m}{\sqrt{k}}\leq \alpha \sqrt{k}$ or not.~$\Box$\medskip

Now, let $(x_n)$ be a sequence of real numbers, such that $x_n\to +\infty$ and $x_n=O(\sqrt{n})$.  There exists $\gamma>0$, such that for any $m$ large enough 
\begin{equation}
\label{eq:uc3}
`P\l(\sup_{k\in\cro{1,m}}\l|\frac{X_1^\DT(k)+iY_1^\DT(k)-k`E(x_1^\DT)}{\sqrt{m}} \r|\geq x_m\r)\leq 4m\exp(-\gamma x_m).
\end{equation}
To see this, use Proposition~\ref{pro:uc} and that for $u$ and $v$ in $\mathbb{R}$, if $|u+iv|\geq a$ then $|u|\geq a/\sqrt{2}$ or $|v|\geq a/\sqrt{2}$. \par

In case the intensity of the PPP is constant and equal to $c>0$, for any $k\geq 1$
\begin{equation}
\label{eq:CE}
X_c(k)\sur{=}{d}\frac{1}{\sqrt{c}}X_1(k),~~~
Y_c(k)
\sur{=}{d}\frac{1}{\sqrt{c}}Y_1(k),
~~~ L_c(k)
\sur{=}{d}\frac{1}{\sqrt{c}}L_1(k)~~~ 
\Xi_c(k)
\sur{=}{d}\frac{1}{\sqrt{c}}\Xi_1(k).
\end{equation}
This allows one to transfer results obtained under $`P_{1}$ to $`P_c$, since for any Borelian $A$ in $`R^k$,
\ben\label{eq:uc4}
`P((X_1(1),\dots,X_1(k))\in A)
&=&`P((X_c(1),\dots,X_c(k))\in \sqrt{c}A).
\een

\subsubsection{Computations for directed \YN}
\label{ssec:CYN}
We use the same notation as in the previous part; this time, for $r>0$,
\ben\label{eq:zeza2}
`P\l(l_1^{\DY}>r\r)&=&`P_1(\#\ovc{k}{0}(r)\cap \bS=0)=
\exp\l(-r^2\theta/2\r),
\een
and $x_1^{\DY} = \cos(\alpha) l_1^{\DY}$ and $y_1^{\DY}=\sin(\alpha)l_1^{\DY}$ for a variable $\alpha$ uniform on $[-\theta/2, \theta/2]$, independent of $l_1^{\DY}$. Again, $l_1^{\DY}$, $x_1^{\DY}$ and $y_1^{\DY}$ own exponential moments, and on gets $`E\l(y_1^{\DY}\r)=0$ and 
\[
`E\l(l_1^{\DY}\r)=\sqrt{\frac{\pi}{2\theta}},~~~
`E\l(x_1^{\DY}\r)=`E(l_1^{\DY})\frac{\sin(\theta/2)}{\theta/2}=\Cbis{\SY}=\Cbis{\Y}.\]
Moreover since $\xi_1^{\DY}=l_1^{\DY} \cos(\theta/2+V)$
where $V$ is uniform on $[-\theta/2,\theta/2]$ and independent of $l_1^{\DY}$, 
\be
`E(\xi_1^{\DY})&=&`E(l_1^{\DY}) \frac{\int_{-\theta/2}^{\theta/2} \cos(\alpha+\theta/2)d\alpha}{\theta}
=`E(l_1^{\DY}) \frac{\sin(\theta)}{\theta}=\Cbor{\Y}
\ee
and again, notice that $\Qbor{\Y}=`E(l_1^{\DY})/`E(\xi_1^{\DY})$, $\Qbis{\Y}=`E(l_1^{\DY})/`E(x_1^{\DY})$. This leads to
\begin{pro} In the \DY\ case, we have
\ben
\l(\frac{L_1^{\DY}(k)}{k},\frac{X_1^{\DY}(k)}{k},\frac{Y_1^{\DY}(k)}{k},\frac{\Xi_1^{\DY}(k)}{k}\r)&\xrightarrow[k]{a.s.}\l(`E(l_1^{\DY}),\Cbis{\Y},0,\Cbor{\Y} \r),
\een
and then $ \dis\frac{L_1^{\DY}(k)}{X_1^{\DY}(k)}\xrightarrow[k]{a.s.} \Qbis{\Y}$, and $\dis\frac{L_1^{\DY}(k)}{\Xi_1^{\DY}(k)}\xrightarrow[k]{a.s.}\Qbor{\Y}.$
\end{pro}
 
\subsubsection{Computations for the random north model}
\label{seq:RNM}
We here consider the case  where the $n(s)$ are i.i.d. uniform on $[0,2\pi]$. Let $\Delta^{\RNT}$ be the traveller stage (when going from $s$ to $t$), and let $x^{\RNT}$ be the orthogonal projection of  $\Delta^{\RNT}$ on the line $(s,t)$ and $y^{\RNT}$ that on the orthogonal of $(s,t)$. A simple computation 
gives $`E\l(y_1^{\RNT}\r)=`E\l(y_1^{\RNY}\r)=0$,
\ben
`E\l(x_1^{\RNT}\r)&=& `E\l(x_1^{\DT}\r)\int_{0}^{\theta/2}\frac{\cos(\nu)}{\theta/2}\,d\nu=`E\l(x_1^{\DT}\r)\frac{\sin(\theta/2)}{\theta/2},\\
`E\l(x_1^{\RNY}\r)&=& `E(l_1^{\DY})\int_{0}^{{\theta}}\int_{0}^{{\theta}}\frac{\cos(\nu-\gamma)}{\theta^2}\,d\nu\,d\gamma=`E(l_1^{\DY})\frac{2-2\cos(\theta)}{\theta^2}.
\een
The limiting quotients $\bQ^{\RNT}$ and $\bQ^{\RNY}$ are $\dis\bQ^{\RNT}=\frac{`E(l_1^{\DT})}{`E(x_1^{\RNT})}$, and $\dis\bQ^{\RNY}=\frac{`E(l_1^{\DY})}{`E(x_1^{\RNY})}$.

\subsection{From local to global or why differential equations come into play}
\label{ssec:DE}

In Proposition~\ref{pro:uc} and in~\eref{eq:uc3} a control of the difference between $X_c^{\DT}(k)$ and $Y_c^{\DT}(k)$  and their mean $k`E(x_1^{\DT})/\sqrt{c}$ and 0 is given. If the intensity is $nf$,  the ``local intensity'' around $s$ is $nf(s)$. Roughly, the first stages of the traveller starting from $s$ are close in distribution under $`P_{nf}$ and under $`P_{nf(s)}$. Hence, at the first order, $X_{nf}^{\DT}(k)$ should be close to $k`E(x_1^{\DT})/\sqrt{nf(s)}$: the speed of the traveller depends  on the position, and therefore a differential equation appears. \par
To approximate $\Pos^{\DT}$ by the solution at a differential equation, we will split the traveller's trajectory into some windows corresponding to some sequence of consecutive stages; the windows have to be small enough to keep the approximation of the local intensity $nf$ by $nf(s)$ to be relevant, and large enough to let the approximation $X_{nf}^{\DT}$ by  $k`E(x_1^{\DT})/\sqrt{nf(s)}$ to be relevant too, that is  large enough to let large number type compensations occur.\par
In this section, we present the tools related to the approximation of the traveller position by differential equations. Their proofs are postponed at the Appendix of the paper.\par
 We start with some deterministic considerations. Consider $D$ an open subset of $`R^d$ and let $G:D\to `R^d$ a Lipschitz function. Denote by $\Eq(G,z)$ the following ordinary differential equation
\begin{equation}
\label{eq:Eq}
\Eq(G,z):=\left\{
\begin{array}{rcl}
y(0)&=&z \in D,\\
y'(x)&=&G(y(x)),~~~~~ \textrm{ for }x\geq 0.
\end{array}
\right.
\end{equation}
This class of equations contains the family of equations $\ODE$ defined in~\eref{eq:ODE}. By Cauchy-Lipschitz Theorem, $\Eq(G,z)$ admits a unique solution $y_{sol(G,z)}$, or more simply $y_{\sol}$ when no confusion on $G$ and $z$ is possible.
This solution is defined on a maximal interval $[0,\lambda(G,z))$, where $\lambda(G,z)$ is the hitting time of the border of $D$ by $y_{\sol}$. \par
Before giving a convergence criterion for random trajectories, here is a deterministic criterion very close to the so-called explicit Euler scheme convergence theorem. 
\begin{lem}
\label{lem:determ-version} Let $G:D\to `R^d$ be a Lipschitz function, $z\in D$ fixed, and $\lambda\in[0,\lambda(G,z))$ given.\par
Let  $(a_n)$ and $(c_n)$ two sequences of positive real numbers going to $0$, and $(y_n)$ a sequence of continuous functions from $[0,\lambda+a_n]$ onto $D$ satisfying the following conditions:\\
$a)$ $y_n(0)=z$ for all $n\geq 1$,\\
$b)$ for all $n\geq 1$, $y_n$ is linear between the points $(ja_n, j\in\cro{0,\floor{\lambda/a_n}})$, and the slope of $y_n$ on these windows of size $a_n$ is well approximated by $G\circ y_n$ in the following sense: 
\[\dis\sup_{j=0,\dots,\floor{\lambda/a_n}} \l|\frac{y_n\l((j+1)a_n\r)-y_n\l(ja_n\r)}{a_n}-G(y_n(ja_n))\r| \leq c_n\]
where we denote with an absolute value a norm in $`R^d$. Under these hypothesis, there exists $C_\lambda>0$, such that for $n$ large enough
\[\sup_{x\in[0,\lambda]}\l|y_n(x)-y_{sol(G,z)}(x)\r|\leq C_\lambda\max\{a_n,c_n\}.\]
Moreover, the constants  $C_\lambda$ can be chosen in such a way that the function $\lambda\mapsto C_\lambda$ is bounded on all compact subsets of $[0,\lambda(G,z))$ and does not depend on the initial condition $z\in D$.
\end{lem}
Note that we need $y_n$ to be defined on a slightly larger interval than $[0,\lambda]$ because of border effects.\medskip

We now extend this lemma to the convergence of a sequence of stochastic processes $(\bZ_n)$. 

\begin{cor}
\label{cor:corn} Let $G:D\to `R^n$ be a Lipschitz function, $z\in D$ fixed, and $\lambda\in[0,\lambda(G,z))$ given. \par
Let $(a_n)$, $(b_n)$, $(c_n)$, $(c_n')$ and $(d_n)$ be five sequences of positive real numbers going to 0; let $(\bZ_n)$ be a sequence of continuous stochastic processes from $[0,\lambda+a_n]$ onto $D$ satisfying the following conditions:\\
$a)$ a.s. $\bZ_n(0)=z$, \\
$b)$ the slope of $\bZ_n$ on windows of size $a_n$ is well approximated by $G\circ \bZ_n$, with a large probability:
\begin{equation}
\label{eq:approx2}
\sup_{j\in \cro{0,\lambda /a_n}} `P\l(\l|\frac{\bZ_n\l((j+1)a_n\r)-\bZ_n\l(ja_n\r)}{a_n}-G\l(\bZ_n(ja_n)\r)\r|\geq c_n\r)\leq a_nd_n;
\end{equation}
$c)$ inside the windows, the fluctuations of $\bZ_n$ are small:
\begin{equation}
\label{eq:mort}
\sup_{j\in\cro{0,\lambda/ a_n}}`P\l( \sup_{x\in [ja_n,(j+1)a_n]} \l|\bZ_n\l(x\r)-\bZ_n\l({ja_n}\r)-(x-ja_n)G(\bZ_n(ja_n))\r|\geq c_n'\r)\leq a_n b_n.
\end{equation}
If these three conditions are satisfied, then 
\[`P\l(\sup_{x\in[0,\lambda]} \l|y_{sol(G,z)}(x)-\bZ_n(x)\r|\leq C_\lambda \max\{a_n,c_n,c_n'\}  \r)\geq 1-(\lambda+1)(d_n+b_n)\]
for a function $\lambda\mapsto C_\lambda$ independent of $z$, bounded on every compact subsets of $[0,\lambda(G,z))$.
\end{cor}
Notice that condition $(c)$ contains $(b)$ if $c_n'\leq c_na_n$ which will be the case in the applications we have, but the present presentation allows one to better understand the underlying phenomenon.

\subsection{On the largest stage and the maximum number of points in a ball}
\label{ssec:fb}

The following quantity 
\[\NAVMAX[\theta](S):=\sup_{s \in \cD[a], \eta\in[0,2\pi]} \inf \l\{h:~\#\l(S\cap \l(s+e^{i\eta}\cam(\theta)(h)\r)\r)\geq 1\r\}\]
is a bound on the largest stages length
for all starting points ($s \in \cD[a]$) and targets ($t\in \cD$) and all navigations considered in this paper.
\begin{lem}\label{lem:petitspas}
For any $\theta \in (0,2\pi]$, any $\C>0$, if $n$ is large enough,
\[`P_{nf}\l(\NAVMAX[\theta]\geq n^{\C-1/2}\r)\leq \exp\l(-n^{\C/2}\r).\]
\end{lem}
\proof Consider a tiling of the plane with squares of size $n^{-1/2+\C/2}$. For $n$ large enough, each element of the family $(s+e^{i\eta}\cam(\theta)(n^{-1/2+\C}),s\in\mathbb{C},\eta\in[0,2\pi])$ contains in its interior at least a square of the tiling. It suffices then to show that any square $\Box$ of the tiling intersecting $\cD[a]$ intersects also $\bS$. But, $`P_{nf}(\#S\cap \Box =0) =\exp(-\int_\Box nf(z)dz)\leq \exp(-m_fn^\C)$ since $nf\geq nm_f$, and the area of $\Box$ is $n^{-1+\C}$. Since $O(n^{1-\C})$ squares intersect the bounded domain $\cD[a]$, by the union bound the probability that there exists a square containing no elements of $\bS$ is $O(n^{1-\C}\exp(-m_fn^\C))$. ~$\Box$.\medskip

Now we turn our attention to
 \[\MAXBALL[r](S):=\max \{\# ( S\cap B(x,r))~|~ x\in \cD[a]\},\]
the maximum number of elements of $S$ in a ball with radius $r$ and having its centre in $\cD[a]$.
\begin{lem}\label{lem:max-boule}
 For any $\B>0$, $`e>0$, if $n$ is large enough
\[`P_{nf}\l(\MAXBALL\l[n^{-1/2+\B}\r] \geq n^{2\B+`e}\r)\leq \exp\l(- n^{2\B+`e}\r).\]
\end{lem}
Note that for $s\in\cD$ the mean number of elements in $B(s,n^{-1/2+\B})\cap \bS(nf)$ is $O(n^{2\B})$.\\
\proof Consider a tiling of the plane with squares of size $a_n=n^{\B-1/2}$; denote by $\Squ_a$ the subset of those having a distance to $\cD[a]$ smaller than $a_n$. Any disk $B(x,a_n)$ with $x\in\cD[a]$ intersects a bounded number $d$ of such squares, with $d$ independent of $\B$. Hence for any positive sequence $(b_n)$,
\[\l\{\MAXBALL[a_n](\bS(nf))\geq b_n\r\}\subset \l\{\sup_{\Box \in \Squ_a}  \#\l(\bS(nf)\cap \Box\r) \geq b_n/d\r\}.\]
First $\#\l(\bS(nf)\cap \Box\r) $ is $\Poi(\int_{\Box}nf(u) du)$ distributed, and $\int_{\Box}nf(u) du\leq \lambda':=\int_{\Box}nM_fdu=na_n^2 M_f$. Secondly, for $X\sim \Poi(\lambda)$, and $\lambda'>\lambda$, $`P(X \geq r)= \sup_{t>0} `P(e^{tX}\geq e^{tr})\leq \sup_{t>0}e^{-tr+\lambda'(e^t-1)}.$
 By a simple (classical) optimisation on $t$, one gets  
\begin{equation}
\label{eq:mb}
`P_{nf}\l( \#\l(\bS\cap \Box\r) >b_n/d\r)\leq \exp\l(-\frac{b_n}d \ln\l(\frac{b_n}{dna_n^2 M_f}\r)+\frac{b_n}{d}-na_n^2M_f\r).\end{equation} 
To end, take $a_n=n^{-1/2+\B}$, $b_n=n^{2\B+`e}$ and apply the union bound to the $O(1/a_n^2)$ squares of $\Squ_a$.\cq

For $\C>0, \B>0,`e>0$, consider the following events 
\[\Omega(n,\C):=\l\{\NAVMAX[\theta]\leq n^{\C-1/2}\r\}, \Omega_{n,`e}:=\l\{\MAXBALL\l[n^{\B-1/2}\r]\leq n^{(2+`e) \B}\r\} \cap \Omega(n,\C).\]

\subsection{About the constants in the paper}
\label{ssec:cc}
The aim of this section is to discuss the role of the constants in this paper, and maybe to help the reader to follow more easily the computations.\\
-- we will use $n^{\C-1/2}$ as a bound on $\NAVMAX$; this bound is valid with probability $1-\exp(-n^{\C/2})$, for $n$ large, \\
-- the behaviour of the traveller close to its target is treated in Section~\ref{sec:term}. ``Far to its target'' (or to some particular points in the trajectory)  means $|t-s|\geq n^{\B-1/2}$. When the traveller enters in the final ball $B(t,n^{\B-1/2})$, he will make at most $\MAXBALL[n^{\B-1/2}]$ additional stages to reach $t$, that is at most $n^{3\B}$ stages with probability $1-\exp(-n^{3\B})$. \\
-- In Corollary~\ref{cor:corn}, when the approximation by the solution of a differential equation is needed, a window with size $a_n$ arises. We will take 
\[a_n=n^{\Nu-1/2}.\] 
Taking into account the space normalisation, this corresponds to consider $n^{\Nu}$ stages of the traveller.\par
Some assumptions are made on the relative values of $\Nu, \C,\B$ (for example $2\Nu+3\C<1/2$) in the statements of this paper. In any case, there is a choice of $(\Nu, \C, \B)$ which fulfils all the requirements of all intermediate results, that is 
\begin{equation}
\label{eq:cond1}
0<\Nu<1/4,~~ 2\Nu+3\C<1/2,~~ \C+\Nu<\B,~~ 0<\B<1/4,~~ \C\in(0,1/2)
\end{equation}
as can easily be checked. Another quantity appearing later on is 
$\Nu'\in(0,1/2)$
related to the bound $n^\alpha$, $\alpha<1/8$ appearing in the paper in most of the important theorems.
This $\alpha$ appears to be $\max\{\Nu(\Nu'-1/2), \Nu(\Nu'+1/2)-1/2,\Nu-1/2\}$ where $(\Nu,\C,\B,\Nu')$.
To reach an $\alpha$ close to $(1/8)^-$, take $\Nu=(1/4)^-$, $\C=0^+$, $\B=(1/4)^-$, $\Nu'=0^+$.

\subsection{Navigations in homogeneous / non homogeneous PPP}
\label{ssec:Lc}

\begin{figure}[ht]
\psfrag{D}{$\Delta$}
\psfrag{s}{$s$}\psfrag{d1}{$d_1$}\psfrag{d2}{$d_2$}
\psfrag{t}{$t$}
\psfrag{f}{$\xi(1)$}
\psfrag{f1}{$\xi(2)$}
\centerline{\includegraphics[height=2.6cm]{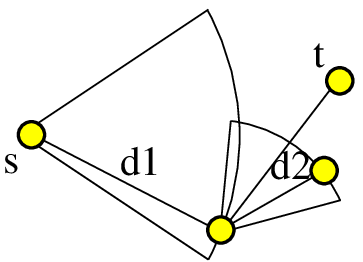}}
\captionn{\label{fig:n-i}Example of configuration of points in the straight Yao navigation leading to an intersection of the decision domains. For \DT\ the decision sectors are simple for $\theta\leq \pi$, and for $\DY$, the condition is $\theta\leq\pi/2$. For $\ST$ and $\SY$, for any $\theta$, the decision sectors can intersect.}
\end{figure}

Here is defined the notion of simple stages, notion related to non-intersecting decision domains. 
\begin{defi} Let $(d_1,\dots,d_k)$ be the $k$ first stages of a traveller going from $s$ to $t$. We say that these stages are $\Nav(s,t)$-simple if the corresponding decision domains do not intersect; this amounts to saying that if the traveller uses as set of possible stops $\{s+\sum_{j=1}^i d_i,~i\leq k\}$, then the decision sectors are non intersecting when going from $s$ to $t$ (see Fig.~\ref{fig:n-i}). A Borelian subset $\Theta$ of $\mathbb{C}^k$ is said to be $\Nav(s,t)$-simple if for any sequence $(d_1,\dots,d_k)\in \Theta$, the stages $(d_1,\dots,d_k)$ are $\Nav(s,t)$-simple.
\end{defi}

Next proposition is really important. It provides around a point $s$, a bound of the deviations of the traveller stages under the non homogeneous $`P_{nf}$  using  the much simpler measure $`P_{nf(s)}$. We were unable to find a sufficient coupling argument; this is a comparison of distributions.\par 
  Let us set the following vectorial notation:
\be  
\Delta^{\DT}[s,j]&:=&\l(\Delta^{\DT}(s,1),\dots,\Delta^{\DT}(s,j)\r);
\ee
in the sequel, $\Delta^{\X}[s,t,j]$ and other similar notation will be used for the other navigation processes $\Nav$ if the target is needed to be specified.
\begin{pro}\label{pro:fundamental2}Let $\X \in \{\Y,\T\}$ and $\theta\leq \pi/3$, or $\X=\SY$ and $\theta < \pi/2$, or $\X=\ST$ and $\theta \leq \pi/2$, or
 $\X= \DT$ and $\theta\leq \pi$, or $\X=\DY$ and $\theta\leq \pi/2$.
 Assume that $(\gamma_n)$ is a sequence of integers such that $\gamma_n\sim n^\Nu$ for some $\Nu\in(0,1/4)$. Let $\C\in(0,1/2)$ such that $2\Nu+3\C<1/2$, and $(s,t)\in \cD'[a]$. For any $\Theta$ Borelian subset of $\mathbb{C}^{\gamma_n}$ $\X(s,t)$-simple,  for $n$ large enough 
\[`P_{nf}\l(\Delta^{\X}[s,t,\gamma_n]\in \Theta\r)\leq  \exp(- n^{\C/2}) +2\ `P_{nf(s)}\l(\Delta^{\X}[s,t,\gamma_n]\in \Theta\r),\]
where ``$t$'' has to be omitted if $\X\in\{\DT,\DY\}$, since the target does not exist in this case. 
\end{pro}
This proposition discusses a very local property since after $n^{w}$ steps the traveller is with great probability at distance $n^{w+C-1/2}<<1$ of its starting point, and then still in ${\cal D}$ for $n$ large.\par
This proposition permits to see that an event with small probability under $`P_{nf(s)}$ (as to observe a large value for $\sum_{i=1}^{\gamma_n} \Delta^{\DT}(s,i)-`E(\sum_{i=1}^{\gamma_n} \Delta^{\DT}(s,i))$ as says Proposition~\ref{pro:uc}) is also small under $`P_{nf}$.
\begin{rem}
The simplicity of the stages is unlikely to arise when $s$ is close to $t$ since the decision domains loose the property to have ``a constant'' direction. For this reason, this proposition will be used only when $s$ is far from $t$. Notice also that the factor 2 in the RHS is taken for convenience, any number greater than 1 also does the job. 
\end{rem}

\prooff{Proposition~\ref{pro:fundamental2}} 
By Lemma~\ref{lem:petitspas}, for $n$ large, $`P_{nf}(\complement\Omega(n,\C))\leq \exp(- n^{\C/2})$. Write 
\[`P_{nf}(\Delta^{\X}[s,t,\gamma_n]\in \Theta)\leq `P_{nf}(\complement\Omega(n,\C))+`P_{nf}(\Delta^{\X}[s,t,\gamma_n]\in \Theta, \Omega(n,\C)).\]
We then just have to bound the second term in the RHS. For this, we will show that  
\be
`P_{nf}(\Delta^{\X}[s,t,\gamma_n]\in \Theta, \Omega(n,\C))&\leq& 2\,`P_{nf(s)}(\Delta^{\X}[s,t,\gamma_n]\in \Theta, \Omega(n,\C))\leq 2\,`P_{nf(s)}(\Delta^{\X}[s,t,\gamma_n]\in \Theta),
\ee
only the first inequality deserving to be proved. We will show more; it is easy to see that under $`P_{nf}$ and under $`P_{nf(s)}$, the distribution of $\Delta^{\X}[s,\gamma_n]$ owns a density with respect to the Lebesgue measure on $\mathbb{C}^{\gamma_n}$, that we denote by $\bD_1$ and $\bD_2$, respectively\footnote{more exactly, the representation of the probability measure of $\Delta^{\X}[s,\gamma_n]$ using densities holds when small stages are considered, small enough such that all decision sectors considered are included in $\cD$, and small enough such that the target is not reached after the $\gamma_n$ first steps}. We will then compare these densities, and bound their ratio on a set of interest. \par
Since the distribution of the stages depend on decision sectors that have a shape depending on $\X$, it is useful here to introduce some notation. Denote by $d^\X(s,t)$ the decision domain (a Camembert section or a triangular domain) corresponding to $\Delta(s,t)^\X$, and by $\partial{d}^\X(s,t)$ the corresponding principal boundaries, namely, depending on $\X$, the side of the triangle which is not on the decision sector, and the arc of circle in the Camembert section. We write with an absolute value the Lebesgue measures of these sets. Let finally ${\sf Seg}^\X(s,t)$ be the length of the segment of the plane starting from $s$, supported by the bisecting line of the decision sector, and whose second end is on $\partial{d}^\X(s,t)$. \par
Under $`P_{nf}$ (resp. $`P_{nf(s)}$), conditionally on $\partial{d}^\X(s,t)$, the density of the law of $\Delta^\X(s,t)$ is proportional to $nf(s+.)$ (resp. is uniform) on $\partial{d}^\X(s,t)$. To characterise the law of $\Delta^\X(s,t)$ it remains to express the distribution of $|{\sf Seg}^\X(s,t)|$; this latter is characterised by the following:
\[`P_{nf}(|{\sf Seg}^X(s,t)|> x)=\exp\l(-\int_{d^\X(s,t,x)} nf(u)du\r),\]
where $d^\X(s,t,x)$ is the decision domain such that $|{\sf Seg}^\X(s,t)|=x$. The density of $|{\sf Seg}^\X(s,t)|$ is then 
\[\l(\int_{\partial d^\X(s,t,x)}nf(u)du\r)\exp\l(-\int_{d^\X(s,t,x)} nf(u)du\r);\]
The density of the same variable under $`P_{nf(s)}$ is obtained by replacing $f(u)$ by the constant value $f(s)$ in this formula.
Finally, under $`P_{nf}$ the density of $\Delta^\X(s,t)$ is 
\begin{equation}\label{ni.1}
g_{nf,s,t}(\delta)
=f(s+\delta)\exp\l(-\int_{d^\X(s,t,x(\delta))} nf(u)du\r)
\end{equation}
for any $\delta$ in the decision sector of the traveller going from $s$ to $t$ using $\X$, where $x(\delta)$ is the length of $|{\sf Seg}^\X(s,t)|$ corresponding to the stage $\delta$ (again, this holds for ``small stages''); under $`P_{nf(s)}$ it is 
\begin{equation}\label{ni.2}
g_{nf(s),s,t}(\delta)=f(s)\exp\l(- \int_{d^\X(s,t,x(\delta))}nf(s) du\r)
\end{equation}
where this last integral value is $nf(s)|d^\X(s,t,x)|=nf(s)|d^\X(.,.,x)|$. Now, let us treat several stages $\sv:=(\delta_1,\dots,\delta_{\gamma_n})\in \Theta$. We will use two ingredients. First the simplicity of $\Theta$ which guaranties the non intersection of the decision domains; this let us use formulas as \eref{ni.1} and \eref{ni.2}, successively. Under the non intersecting condition by successively conditioning on the first stages of the traveller it appears that $\bD_1$ and $\bD_2$ have a multiplicative form:
\ben\label{eq:D2}
{\bf D}_1(\sv)&=&\prod_{i=1}^{\gamma_n} g_{nf(s),s_{i-1},t}(\delta_i), ~~ {\bf D}_2(\sv)=\prod_{i=1}^{\gamma_n} g_{nf,s_{i-1},t}(\delta_i)
\een
where $s_0=s$, $s_j=s_0+\delta_1+\dots+\delta_{j}$. The sequence $(d^\X(s_i,t,x(\delta_i)),i=1,\dots,\gamma_n)$ appearing is the same in both formulas, and then $\bD_1$ and $\bD_2$ are very similar. Let us see why ${\bf D}_2(\sv) \leq 2{\bf D}_1(\sv)$ holds for any element $\sv$ in $\Omega(n,\C)\cap \Theta$ which will be enough to conclude.\par
For this, consider the first factor $\prod_{i=1}^{\gamma_n} f(s_{i-1}+\delta_i)$ appearing in $\bD_2$. If $\sv\in\Omega(n,\C)$, for $m\leq \gamma_n$, since $f$ is Lipschitzian,
\[\prod_{i=1}^{\gamma_n} f(s_{i-1}+\delta_i)\leq \prod_{i=1}^{\gamma_n}f(s)+\alpha_f\,\gamma_nn^{\C-1/2}\leq f(s)^{\gamma_n}(1+c_1 \gamma_n n^{\C-1/2})^{\gamma_n},\] for constants $c_1>0$ and $\alpha_f$, since $f$ is bounded. Let us bound the second factor appearing in  $\bD_2$ : 
\ben\label{eq:eatrz}
\exp\l(-n\int_{d^\X(s_{m-1},t,\delta_m)}f(u)du\r)&=&\exp\l(-n\int_{d^\X(s_{m-1},t,\delta_m)}f(u)-f(s)du\r)\\
&&\times \exp\l(-n\int_{d^\X(s_{m-1},t,\delta_m)}f(s)du\r)
\een
the second term of this  product is simply $\exp\l(-nf(s)|d^\X(s_{m-1},t,\delta_m)|\r)$; let us bound the first one. Using that the area of $d^\X(s_{m-1},t,\delta_m)$ is bounded by $c_2 n^{2\C-1}$ (since $\delta_m\leq n^{\C-1/2}$), and that $f(u)-f(s)$ is greater than $-\alpha_f m n^{\C-1/2}$ (Lipschitz, and since we are in $\Omega(n,\C)$, $|t-s|\leq m\,n^{\C-1/2}$), we get that the LHS of~\eref{eq:eatrz} is bounded by 
\be
\exp\l(-n\int_{d^\X(s_{m-1},t,\delta_m)}f(u)-f(s)du\r) &\leq&\exp\l(c_2 n\times m n^{\C-1/2} n^{2\C-1}\r)=\exp(c_2\,m\,n^{3\C-1/2})
\ee
for some constant $c_2>0$. Finally, putting together the $\gamma_n$ terms involved,
\begin{equation}
{\bf D}_2(\sv) \leq {\bf D}_1(v)\times\l(1+c_1 \gamma_nn^{\C-1/2}\r)^{\gamma_n}\times \prod_{m=1}^{\gamma_n} \exp\l(c_2 m n^{3\C-1/2}\r).
\end{equation}
Recall that $\gamma_n\sim n^{\Nu}\to +\infty$. The first term $(1+c_1 \gamma_nn^{\C-1/2})^{\gamma_n}$ goes to 1 if $\C<1/2$.  The second term goes to 1 if $\sum_{m=1}^{\gamma_n}c_2 m n^{3\C-1/2}$ goes to 0, that is if $2\Nu+3\C-1/2<0$. Finally, taking $\C\in(0,1/2)$ and $\Nu>0$ satisfying $2\Nu+3\C-1/2<0$, we see that $\l(1+c_1 \gamma_nn^{\C-1/2}\r)^{\gamma_n}\times \prod_{m=1}^{\gamma_n} \exp\l(c_2 m n^{3\C-1/2}\r)\to 1$, and thus is less than $2$ for $n$ large enough.~\cq

\subsection{Limiting behaviour for a traveller using \DT\ or \DY}
\label{ssec:lbt}
In this section, we deal with \DT, but \DY\ can be studied similarly.
The aim is to show that when $n\to +\infty$, after rescaling, the function $\Pos_{s}^{\DT}$  under $`P_{nf}$ converges in distribution to $\SOL_{s}^{\Cbis{\T},0}$ the solution of $\ODE(\Cbis{\T},0,s)$ as defined in \eref{eq:ODE}. 
\begin{theo}\label{theo:clef}
 Let $\X= \DT$ and $\theta\leq \pi$ or $\X=\DY$ and $\theta\leq \pi/2$. 
 For any $s\in\cD[a]$, any $\alpha\in(0,1/8)$, for $\lambda \in[0,\lambda(F_{\Cbis{\X},0},s))$, there exists a constant $d>0$ such that for $n$ large enough
\[`P_{nf}\l(\sup_{x\in[0,\lambda]} \l|\Pos^{\DT}_s(x\sqrt{n})-\SOL_s^{\Cbis{\T},0}(x)\r|\leq n^{-\alpha}  \r)\geq 1-\exp(-n^d).\]
\end{theo}
Notice that $\ODE(\Cbis{\T},0,s)$ coincides with $\Eq\l(F_{\Cbis{\T},0},s\r)$ as introduced in~\eref{eq:Eq} (see also the considerations about $\lambda(G,z)$, the hitting time of $\complement \cD$ by the solution of $\Eq(G,z)$). 
Hence, this theorem is a consequence of the following lemma 
\begin{lem}\label{lem:le2}
Let $\lambda\in[0,\lambda(F_{\Cbis{\T},0},s))$, $\gamma_n=\floor{n^\Nu}$ and $a_n=\gamma_n/\sqrt{n}\sim n^{\Nu-1/2}$ be the size of the ``windows''.
 Let $\C,\Nu,\Nu'$ be some positive constants, such that $2\Nu+3\C-1/2<0$ and $\Nu'<1/2$. Then for any $s\in \cD[a]$, the sequence $(x\mapsto\Pos^{\DT}_s(x\sqrt{n}))$ satisfies the hypothesis of Corollary~\ref{cor:corn}, for $G=F_{\Cbis{\T},0}$, $c_n=n^{\Nu(\Nu'-1/2)}/\sqrt{m_f}$, $b_n=\exp(-n^{\min(\Nu\Nu'/2,\C/3)})$, $c_n'=n^{\Nu(\Nu'+1/2)-1/2}/\sqrt{m_f}$, $d_n=b_n$.
Thus 
\[`P_{nf}\l(\sup_{x\in[0,\lambda]} \l|\SOL_s^{\Cbis{\T},0}(x)-\Pos^{\DT}_s(x\sqrt{n})\r|\leq C_\lambda \max\{a_n,c'_n,c_n\}  \r)\geq 1-2(\lambda+1)b_n,\]
for $\lambda\to C_\lambda$ bounded on compact sets.
\end{lem} 
Here, the minimum value for $\max\{a_n,c_n,c'_n\}$ is $n^{-1/8+}$ as explained in Section~\ref{ssec:cc} ($\Nu=(1/4)^-, \C=0^+,\Nu'=0^+$).\par
In the homogeneous PPP $`P_{nc}$ (for some $c>0$), the right order of the variance of $\Pos^{\DT}_s(x\sqrt{n})$ is $1/\sqrt{n}$, since it is a sum of $\sqrt{n}\lambda$ random variables with variance of order $1/n$ by \eref{eq:Poisson1}. Then standard deviations have order $n^{1/4}$. Here the constant $1/8$ arising in the results is not so good, but gives exponential bounds needed here, and are valid also for non homogeneous PPP. \medskip

\prooff{Lemma~\ref{lem:le2}}For short, we write  $\bZ_s^{(n)}(x)$ instead of $\Pos^{\DT}_s(x\sqrt{n})$. 
We will use Proposition~\ref{pro:uc}, Formulas~\eref{eq:uc3} and~\eref{eq:uc4} and will establish some bounds valid for the first  $\gamma_n$ stages of a traveller starting from a generic point $s_0 \in \mathbb{C}$ (that is $\bZ_s^{(n)}(0)=s_0$). Recall that $\Cbis{\T}=`E(x_1^{\DT})$ and consider the following Borelian subset of $\mathbb{C}^{\gamma_n}$,
\begin{equation}
\label{eq:Theta1}
\Theta_n^{(1)}:=\l\{ (\lambda_1,\dots,\lambda_{\gamma_n}) \in \mathbb{C}^{\gamma_n} ,~\sup_{l\in\cro{1,\gamma_n}}\l|\sum_{j=1}^{l} \l(\lambda_j -\frac{\Cbis{\T}}{\sqrt{nf(s_0)}}\r)\r|\leq y_n \r\}
\end{equation}
for $(y_n)$ a sequence that will be fixed later on; notice that $\frac{\Cbis{\T}}{\sqrt{f(s_0)}}=F_{\Cbis{\T},0}(s_0)=F_{\Cbis{\T},0}(\bZ_s^{(n)}(0))$. In term of events, 
\begin{eqnarray*}
\left\{\Delta^\DT[s_0,\gamma_n]\in \Theta_n^{(1)}\right\}&=&\left\{\sup_{l\in\cro{1,\gamma_n}}\l|\sum_{j=1}^{l} \l(\Delta^\DT(s_0,j) -\frac{\Cbis{\T}}{\sqrt{nf(s_0)}}\r)\r|\leq y_n\right\}.
\end{eqnarray*}
Using $\bZ_{s_0}^{(n)}(j/\sqrt{n})=s_0+\sum_{l=1}^{j} \Delta^\DT(s_0,l),$
\ben\label{eq:bic2}
\l\{\Delta^\DT[s_0,\gamma_n]\in\Theta_n^{(1)}\r\} &=&
\l\{\sup_{l\in\cro{1,\gamma_n}}\l| \bZ_{s_0}^{(n)}(l/\sqrt{n})-s_0 -l\frac{\Cbis{\T}}{\sqrt{nf(s_0)}}\r|\leq y_n\r\}.
\een
Therefore, the event $\Theta_n^{(1)}$ contains (for $l=\gamma_n$)
\be
\l\{\l| \frac{\bZ_{s_0}^{(n)}(a_n)-s_0}{a_n} -F_{\Cbis{\T},0}(\bZ_{s_0}^{(n)}(0))\r|\leq y_n/a_n\r\},
\ee
and
\be
\l\{\sup_{x\in[0,a_n]}\l| \bZ_{s_0}^{(n)}(x)-\bZ_{s_0}^{(n)}(0)-(x-0) F_{\Cbis{\T},0}(\bZ_{s_0}^{(n)}(0))\r|\leq y_n\r\}.
\ee
We now take  $y_n=  \frac{x_{\gamma_n}\sqrt{\gamma_n}}{\sqrt{n m_f}}$ for  $x_n\sim n^{\Nu'}$. Since $x_{n}=O(\sqrt{n})$, it fulfils the requirement of Proposition~\ref{pro:uc}. 
Moreover, since $2\Nu+3\C-1/2<0$, the comparison provided by Proposition~\ref{pro:fundamental2} is valid and thus we work for a moment under $`P_{nf(s_0)}$. We get, using also the rescaling~\eref{eq:CE},
\be
`P_{nf(s_0)}\l(\Delta^\DT[s_0,\gamma_n]\in \complement\Theta_n^{(1)}\r)&=&`P_{nf(s_0)}\l(\sup_{l\in\cro{1,\gamma_n}}\l|\sum_{j=1}^{l} \l(\Delta^\DT(s_0,j) -\frac{\Cbis{\T}}{\sqrt{nf(s_0)}}\r)\r|\geq \frac{x_{\gamma_n}\sqrt{\gamma_n}}{\sqrt{n m_f}}\r)\\
&\leq&`P_{1}\l(\sup_{l\in\cro{1,\gamma_n}}\l|\sum_{j=1}^{l}\frac{\Delta^\DT(s_0,j) -\Cbis{\T}}{\sqrt{\gamma_n}}\r|\geq x_{\gamma_n}\r)\\
&\leq& 2\gamma_n\exp(-\gamma' x_{\gamma_n})\leq \exp(-\gamma''n^{\Nu\Nu'})
\ee
for  constants $\gamma'>0,\gamma''>0$, and for $n$ large enough (where it has been used that $\frac{\sqrt{f(s_0)}}{\sqrt{ m_f}}\geq 1$). Take $c_n=\frac{x_{\gamma_n}}{\sqrt{\gamma_nm_f}}\sim n^{\Nu\Nu'-\Nu/2}/\sqrt{m_f}$ (this goes to 0) and $c_n'=\frac{x_{\gamma_n}\sqrt{\gamma_n}}{\sqrt{n m_f}}\sim n^{\Nu\Nu'+\Nu/2-1/2}/\sqrt{m_f}$ (this goes to 0). Then we have established that Formulas~\eref{eq:approx2} and~\eref{eq:mort} in Corollary~\ref{cor:corn} hold true under $`P_{nf}$, with the left most signs ``$\sup$'' deleted, for $j=1$, $G=F_{\Cbis{\T},0}$, $\bZ_n=\bZ^{(n)}$, $b_n=d_n\geq (2 \exp(-\gamma'' n^{\Nu\Nu'})+\exp(- n^{\C/2}))/a_n$. For example we can take $b_n=d_n=\exp(-n^{\min(\Nu\Nu'/2,\C/3)})$.\par
Since these bounds are valid for any starting points $s_0$, and any $j$,  Formulas~\eref{eq:approx2} and~\eref{eq:mort} in Corollary~\ref{cor:corn} hold true in this case with the supremum sign re-established, by Markovianity of the sequence $(\bZ^{n}(j/\sqrt{n}),j\geq 0)$. The assumptions of Corollary~\ref{cor:corn} are satisfied.  The conclusion of this corollary entails those of the present lemma. ~\cq

\subsection{Local representation of navigations using  directed navigations}
\label{ssec:local-with-directed} 
The aim is to represent locally around a point $s$, the first stages of a navigation (cross or straight) under the homogeneous PPP $`P_{nf(s)}$ with the first stages of directed navigation.

\subsubsection*{Local representation of  straight navigation using directed navigation}

\noindent We are here comparing the stages of \ST~and \DT. Same results for $\SY$ holds true also. 
\begin{lem}\label{lem:ST-NI}
Let $S\in \Omega(n,\C)$,  $(s,t)\in \cD'[a]$ such that $|t-s|\geq n^{\B-1/2}$. The $\ST$-decision domains of $\gamma_n\sim n^{\Nu}$ first stages are $\ST(s,t)$-simple, if $\C+\Nu<\B$, for $n$ large enough.
\end{lem}
\proof If $|t-s|\geq n^{\B-1/2}$, since $S\in\Omega(n,\C)$, we have $|t-s_i|\geq n^{\B-1/2}-\gamma_nn^{\C-1/2}\geq n^{\B-1/2}/2$ and $|s-s_i|\leq \gamma_nn^{\C-1/2}$ for $i\leq \gamma_n$ (for $n$ large enough). Therefore, for $i\leq \gamma_n$,  $|\arg(t-s_i)-\arg(t-s)|$ is bounded above by $O(n^{\C+\Nu-\B})$. The angle between the bisecting lines of the decision domains are going to 0 uniformly; then the decision domains are non intersecting.~$\Box$

\begin{lem}\label{lem:c2}
Let $(s,t) \in \cD'[a]$ such that $|t-s|\geq n^{\B-1/2}$ and 
$\gamma_n\sim n^{\Nu}$. Assume that $\C+\Nu<\B$. For any $`e>0$, for any $x>0$, under $`P_{nx}$,  the vectors $\Delta^{\ST}[s,t,\gamma_n]\1_{[0,n^{\C-1/2}]}(\max|\Delta^{\ST}[s,t,\gamma_n]|)$ and \[\l(e^{i\arg(t-s_{j-1})}\Delta^{\DT}(s,j),j=1,\dots,\gamma_n\r)\1_{[0,n^{\C-1/2}]}\l(\max|\Delta^{\DT}[s,\gamma_n]|\r),\] where $s_j=\sum_{i=0}^{j-1}\Delta^{\DT}(s,i)+s$, have same distribution when $n$ is large enough. 
\end{lem}
 Note that under $`P_{nx}$ the stages $\Delta^{\DT}[s,\gamma_n]$ are i.i.d. whereas the stages $\Delta^{\ST}[s,\gamma_n]$ are not. \medskip

\noindent \proof This is a simple consequence that under the hypothesis, the decision domains under \ST~ or \DT~ are simple; then the distribution of the stages in both cases are given by the same computations (based on areas of triangles). ~$\Box$

\subsubsection*{Local  representation of  cross navigation using directed navigation}
%\label{sssec:CN-DN}
We treat here the case $\T$ but $\Y$ can be treated similarly.
For any two points $(s,t)\in \cD$,  let  $K_{s,t}:=\l\lfloor\frac{\arg(t-s)}{\theta}+1/2\r\rfloor$ be an integer $\kappa$ such that $t\in\ov{\kappa}{s}$.
\begin{lem}\label{lem:NIT} Let $\theta\leq \pi/3$. Let $S\in \Omega(n,\C)$, and $(s,t)\in \cD'[a]$ such that $|t-s|\geq n^{\B-1/2}$. The $\T$-decision domains for the $\gamma_n\sim n^{\Nu}$ first stages are $\T(s,t)$-simple, if $\C+\Nu<\B$, for large $n$.
\end{lem}
\proof Two cases have to be considered.~\\
-- When $t$ is far from $\cross(s)$, $K_{s_i,t}$ stays constant for small values of $i$. Hence the bisecting line of the decision domains of the traveller are parallel, and then the domains do not intersect.\\
-- When $t$ is close to $\cross(s)$ but $|s-t|\geq n^{\B-1/2}$, then during the first $n^\Nu$ stages, $t$ remains close  to $\cross(s_i)$, but $|t-s_i|\geq |t-s|-n^{\Nu+\C-1/2}\geq n^{-1/2}(n^\B-n^{\C+\Nu})$). Therefore $K_{s_i,t}$ can take two values  $K_{s,t}$ and $K_{s,t}\pm 1 \mod p_\theta$ depending on the position of $t$ with respect to $\cross(s)$. Therefore the bisecting lines of the decision domains have two possible directions. The angle between these possible directions being $\theta\leq\pi/3$, the decision domains are non intersecting. ~$\Box$ \medskip

We then give a representation of the increments of the cross navigation, immediate since as for Lemma~\ref{lem:c2}, it just relies on the fact that some triangles have the same area (and on Lemma~\ref{lem:NIT}).
\begin{lem}\label{lem:c3} 
Let $(s,t)\in \cD'[a]$ such that $|t-s|\geq n^{\B-1/2}$. Let $\gamma_n\sim n^{\Nu}$ and $\C+\Nu<\B$. 
For any $x>0$, under $`P_{nx}$, the variables
\[\l(e^{i\theta\kappa_j}\Delta^{\DT}(s,j),j=1,\dots,\gamma_n\r)\1_{[0,n^{\C-1/2}]}\l(\max|\Delta^{\DT}[s,\gamma_n]|\r)\]
and $\Delta^{\T}[s,t,\gamma_n]\1_{[0,n^{\C-1/2}]}\l(\max|\Delta^{\T}[s,t,\gamma_n]|\r)$ where $\kappa_j=K_{s+\sum_{l=1}^{j-1}\Delta^{\DT}(s,l),t}$ have same law.
\end{lem} 
Let us examine the consequence of this lemma. As in the proof of Lemma~\ref{lem:NIT} two cases occur.

\subsubsection*{When $t$ is far from $\cross(s)$}
Assume that we are in $\Omega(n,\C)$. And assume that $d(t,\cross(s))\geq  n^{\Nu+\C-1/2}$, in words, the distance between the point $t$ and the set $\cross(s)$  is greater than $n^{\Nu+\C-1/2}$. In this case, $K_{s_i,t}$ is equal to $K_{s,t}$ for $i\leq n^\Nu$.
Therefore the property of the preceding lemma rewrites
\begin{equation}
\label{eq:rep-T-1}
\Delta^{\T}[s,t,\gamma_n]\1_{[0,n^{\C-1/2}]}\l(\max|\Delta^{\T}[s,t,\gamma_n]|\r) =e^{i\theta K_{s,t}}\Delta^{\DT}[s,\gamma_n]\1_{[0,n^{\C-1/2}]}\l(\max|\Delta^{\DT}[s,\gamma_n]|\r),
\end{equation}
and then, under these conditions, \T\ coincides with the \DT\ with direction $\theta K_{s,t}$.

\subsubsection*{When $t$ is close to $\cross(s)$}
%\label{sssec:close-cross}
Assume now that $d(t,\cross(s))\leq  n^{\Nu+\C-1/2}$ but $|t-s|\geq n^{\B-1/2}$.
There exists a unique $k$ such that $d(t,\HL_k(s))\leq  n^{\Nu+\C-1/2}$ for $n$ large enough (since $\B>\Nu+\C$). Consider the line $\bD$ parallel to $\HL_k(s)$ passing via $t$ ($\bD$ is included in $\cross(t)$ if $p_\theta$ is even). We then have $d(\bD,s)\leq  n^{\Nu+\C-1/2}$.\par
\begin{figure}
\psfrag{D}{$\bD$}
\psfrag{s}{$s$}
\psfrag{t}{$t$}
\psfrag{f}{$\xi(1)$}
\psfrag{f1}{$\xi(2)$}
\centerline{\includegraphics[angle=00,height=5cm]{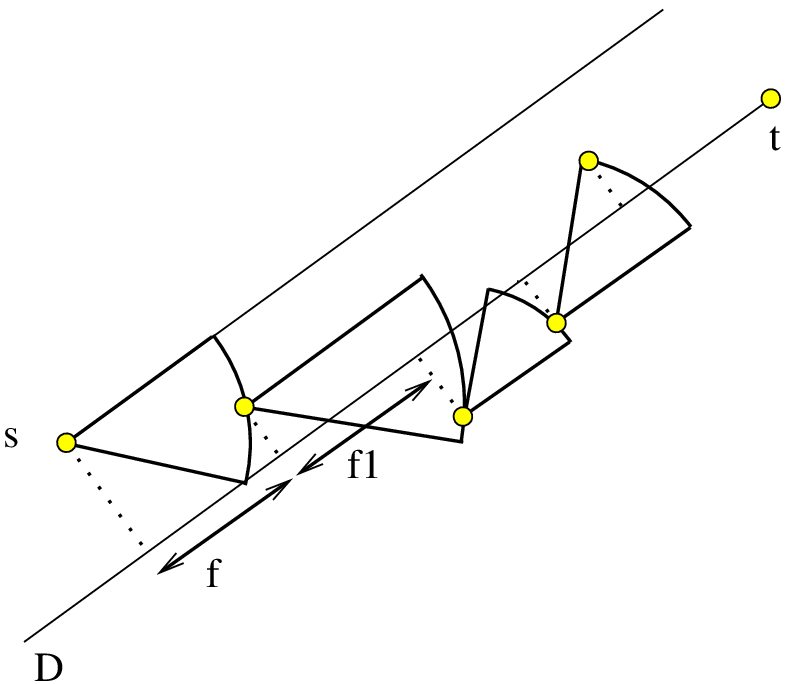}}
\captionn{\label{fig:proj}Progression of the traveller near the line $\bD$. Note that the decision sectors have always a side parallel to $\bD$. }
\end{figure}

Two things are needed to be stated:\\
-- for any $i\leq n^\Nu$, $d(\bD,s_i)\leq  n^{\Nu+\C-1/2}$. In words, if the traveller is close to $\bD$ at some time, it stays close to it afterward. The reason is simple and comes from the second point: the decision domains of the traveller has a border parallel to $\bD$ (see Fig.~\ref{fig:proj}).\\
-- the decision domains of the traveller for these $n^\Nu$ steps have a border parallel to $\bD$, and the other border, of course presents an angle $\theta$ with $\bD$. Therefore, the orthogonal projection of the stages on $\bD$ of all of these stages, have the same distribution (see Fig.~\ref{fig:proj}).

\section{Proofs of the theorems}
\label{sec:PT}

The proofs of our theorems are decomposed in several parts. \\
(a) First, we prove that for any $(s,t)$ fixed, for any navigation $\X$, the function $\Pos_{s,t}^{\X}$ admits a limit (specified in the different theorems of the paper). 
The different costs associated with the path starting by its length -- which indeed appears as a cost, and which can not be handled without knowing the position of the traveller -- is treated afterward.\\
(b) The result for ``one trajectory'' is then extended to all trajectories in several steps: first, the results from (a) arrive with some probability bounds that allows one to handle in once a polynomial number of trajectories (for starting and ending points in a grid). Then the paths between other points are treated by comparison with these trajectories (this is Section~\ref{ssec:GB}).

\subsection{Result for one trajectory}
\label{ssec:RT}

\subsubsection{Straight navigation}
\label{sssec:Surg-straight}

 For any function $g$ and any set $L$, the hitting time of $L$ by $g$ is
$\tau[g](L):=\inf\{x,~g(x)\in L\}$. For example, $\tau\l[\Pos_{s,t}^{\infty,\Cbis{\X}}\r](\{t\})=\Time_{s,t}^{\Cbis{\X}}$. \par
We introduce a uniform big O notation: let $(g_n)$ be a sequence of functions $g_n:A\to `R$, and $(c_n)$ a sequence of real numbers. Notation $g_n=O_{A}(c_n)$ means that $\sup_{y\in A}g_n(y)=O( c_n)$.
\begin{theo}\label{theo:clef2}
Let $\X=\SY$ and $\theta < \pi/2$, or $\X=\ST$ and $\theta \leq \pi/2$. 
For any $\alpha\in(0,1/8)$, there exists $c=c(\alpha)$ such that
\begin{equation}
`P_{nf}\l(\sup_{x\leq \Time_{s,t}^{\Cbis{\X}}} \l|\Pos_{s,t}^{\X}(x\sqrt{n})-\Pos_{s,t}^{\infty,\Cbis{\X}}(x)\r|\geq n^{-\alpha}\r)=O_{\cD'[a]}\l(\exp(-n^{c})\r).
\end{equation}
Moreover there exists $d>0$ such that for $n$ large enough  
\[`P_{nf}\l(\l|\Nb^\X(s,t)/\sqrt{n}-\Time_{s,t}^{\Cbis{\X}}\r|\geq dn^{-\alpha}\r)=O_{\cD'[a]}\l(\exp(-n^{c})\r).\]
\end{theo}
Notice that the second assertion of this theorem is not a direct consequence of the first one since when $\|g_n-g\|_{\infty}\to 0$ it may happen that $\tau[g_n](A)\not\to\tau[g](A)$ for some set $A$. Notice also  for $\X\in\{\ST,\SY\}$, $\Time_{s,t}^{\Cbis{\X}}$ is implicitly known by 
\begin{equation}
\Time_{s,t}^{\Cbis{\X}}:=\inf\l\{~x,\int_{0}^x \frac{e^{i\arg(t-s)}\Cbis{\X}}{\sqrt{f(s+e^{i\arg(t-s)u})}}du=t-s\r\}.
\end{equation}
If $f=c$ is constant this simplifies and we get
\begin{equation}
\Time_{s,t}^{\Cbis{\X}}={|t-s|\sqrt{c}}\,/\,{\Cbis{X}}.
\end{equation}
\prooff{ Theorem~\ref{theo:clef2}} We give the proof in the case $\X=\ST$ the other case $\X=\SY$ is similar. For short, write $\bar{\bZ}_s^{(n)}(x)$ instead of $\Pos_{s,t}^{\ST}(x\sqrt{n})$. For $\B\in(0,1/4)$, and $`e>0$ consider the set $\Omega_{n,`e}$ defined in Section~\ref{ssec:fb}.
By Lemmas~\ref{lem:petitspas} and~\ref{lem:max-boule},
$`P_{nf}(\complement \Omega_{n,`e})\leq \exp(- n^{c_1})$ for some $c_1>0$, for $n$ large enough. Let us assume that $S\in\Omega_{n,`e}$. Consider $T^{n}_{s,t}:=\inf\l\{x,~|\bar{\bZ}^{(n)}_s(x)-t|<n^{\B-1/2}\r\}$ the hitting time of $B(t,n^{\B-1/2})$ by $\bar{\bZ}^{(n)}_s$. By Proposition~\ref{pro:boule}, $T^{n}_{s,t}<+\infty$ almost surely, and for any $x\geq T^{n}_{s,t}$, 
\begin{equation}
\label{eq:nb12}
|\bar{\bZ}_{x}^{(n)}(x)-t|\leq n^{\B-1/2}.
\end{equation}
For $S\in\Omega_{n,`e}$, the total number of stops inside $B(t,n^{\B-1/2})$ is at most $n^{(2+`e)\B}$, which corresponds on the process $\bar{\bZ}^{(n)}$ to a negligible time interval $n^{(2+`e)\B-1/2}=o(1)$ if $(2+`e)\B<1/2$. The space fluctuations of these last stages are at most of order $n^{\B-1/2}$, and then, they are negligible when:
\[\B-1/2<-\alpha.\]
The control of the position of the traveller along the rest of the trajectory will be done by Corollary~\ref{cor:corn}, using together the elements of the proof of Theorem~\ref{theo:clef2}, and Lemma~\ref{lem:c2}.\par 
Assume that $|s-t|>n^{\B-1/2}$,  $\Nu>0$, $\gamma_n=n^\Nu$ and  $y_n=\frac{x_{\gamma_n}\sqrt{\gamma_n}}{\sqrt{n m_f}}$ for  $x_n\sim n^{\Nu'}$, $\Nu'\in(0,1/2)$ as in the proof of Theorem~\ref{theo:clef}. 
By Proposition~\ref{pro:fundamental2}, we know that for a \ST(s,t)-simple set $\Theta$,
\begin{equation}
\label{eq:borne7}
`P_{nf}\l(\Delta^{\ST}[s,t,\gamma_n]\in \Theta\r)\leq  \exp(- n^{\C/2}) +2`P_{nf(s)}\l(\Delta^{\ST}[s,t,\gamma_n]\in \Theta\r).
\end{equation}
We then work under $`P_{nf(s)}$ from now on.
Set, for $(s,t)$ fixed,
\[\app{\Psi}{\mathbb{C}^{\gamma_n}}{\mathbb{C}^{\gamma_n}}{\lambda:=(\lambda_1,\dots,\lambda_{\gamma_n})}{(e^{i\arg(t-s_{k-1}(\lambda))}\lambda_k,k=1,\dots,\gamma_n)\1_{[0,n^{\C-1/2}]}(\max_{l\in\cro{1,\gamma_n}}|\lambda_l|)}\]
where $s_{k-1}(\lambda)=s+\sum_{i=1}^{k-1}\lambda_i$. Hence, Lemma~\ref{lem:c2} says that under $`P_{nf(s)}$, 
\begin{equation}
\label{eq:repu}
\Psi(\Delta^{\DT}[s,\gamma_n])\1_{[0,n^{\C-1/2}]}(\max \Delta^{\DT}[s,\gamma_n])\sur{=}{d}\Delta^{\ST}[s,t,\gamma_n]\1_{[0,n^{\C-1/2}]}(\max|\Delta^{\ST}[s,t,\gamma_n]|).
\end{equation}
For $\Theta_n^{(1)}$ given in \eref{eq:Theta1}, define
\[\Theta_n^{(2)}:=\Theta_n^{(1)}\cap \l\{ (\lambda_1,\dots,\lambda_{\gamma_n}) \in \mathbb{C}^{\gamma_n} ,~\max_{l\in\cro{1,\gamma_n}} |\lambda_l|<n^{\C-1/2}\r\},\]
By Lemmas~\ref{lem:le2} and~\ref{lem:petitspas}, for some $\gamma''$ and for $n$ large enough,
\begin{equation}
\label{eq:law}
`P_{nf(s)}\l(\Delta^{\DT}[s,\gamma_n]\in \Theta_n^{(2)}\r)\geq 1-\exp(-\gamma''n^{\Nu\Nu'})-\exp(-n^{\C/2}).
\end{equation}
If $\Delta^{\DT}[s,\gamma_n]\in\Theta_n^{(2)}$ then it is also in $\Theta_n^{(1)}$ and then most of the equalities or set inclusions of the proof of Lemma~\ref{lem:le2} can be recycled here, starting from $\left\{\Delta^{\DT}[s_0,\gamma_n]\in \Theta_n^{(2)}\right\}$
\begin{eqnarray*}
&=&\left\{\max_{l\in\cro{1,\gamma_n}}\l|\sum_{j=1}^{l} \l(\Delta^{\DT}(s_0,j) -\frac{\Cbis{\T}}{\sqrt{nf(s_0)}}\r)\r|\leq y_n,\max \Delta^{\DT}[s_0,\gamma_n]<n^{\C-1/2} \right\}.
\end{eqnarray*}
For $\Delta^{\DT}[s_0,\gamma_n]\in \Theta_n^{(2)}$,  for $\bZ^{(n)}$ defined as in the proof of Lemma~\ref{lem:le2},  since $|e^{i\arg(t-s)}|=1$,
\ben\label{eq:rot1bis}
\l| e^{i\arg(t-s)}\frac{\bZ_{s_0}^{(n)}(a_n)-s_0}{a_n} -e^{i\arg(t-s)}F_{\Cbis{\T},0}(\bZ_{s_0}^{(n)}(0))\r|&\leq& y_n/a_n,\\\label{eq:rot2bis}
\sup_{x\in[0,a_n]}\l| e^{i\arg(t-s)}(\bZ_{s_0}^{(n)}(x)-\bZ_{s_0}^{(n)}(0))-(x-0)e^{i\arg(t-s)} F_{\Cbis{\T},0}(\bZ_{s_0}^{(n)}(0))\r|&\leq& y_n.
\een
Consider now the increments $\widehat{\Delta}[s,t,\gamma_n]:=\Psi(\Delta^{\DT}[s,\gamma_n]).$ 
Using \eref{eq:law} and \eref{eq:repu}, for any Borelian $\Theta$,
\[|`P_{nf(s)}(\widehat{\Delta}[s,t,\gamma_n]\in\Theta)-`P_{nf(s)}(\Delta^{\ST}[s,t,\gamma_n]\in\Theta)|\leq \exp(-\gamma''n^{\Nu\Nu'})+\exp(-n^{\C/2}).\]
Up to this exponentially small probability, we may work with $\widehat{\Delta}[s,t,\gamma_n]$ instead of $\Delta^{\ST}[s,t,\gamma_n]$.  For $\Delta^{\DT}[s,\gamma_n]\in\Theta_n^{(2)}$ and $k\leq \gamma_n$,  let us bound
\[d:=\l|\sum_{j=1}^k(\widehat{\Delta}(s,t,j)-e^{i\arg(t-s)}\Delta^{\DT}(j))\r|.\]
Writing $S_j=\sum_{m=1}^{j}\Delta^{\DT}(j)$,
\ben\nonumber
d&=&\l|\sum_{j=1}^ke^{i\arg(t-S_{j-1})}\Delta^{\DT}(j)-e^{i\arg(t-s)}\Delta^{\DT}(j)\r|\\
\label{eq:comp1} &\leq&\max\{|e^{i\arg(t-S_{j-1})}-e^{i\arg(t-s)}|,j\in\cro{0,k}\}\times  \sum_{j=1}^k|\Delta^{\DT}(j)|.
\een
Since $\Delta^{\DT}[s,\gamma_n]\in\Theta_n^{(2)}$, the second  term in the RHS of \eref{eq:comp1} is bounded by $k\Cbis{\T}/\sqrt{nf(s)}+y_n$ which is smaller than $O_{\cD'[a]}(a_n+y_n)=O_{\cD'[a]}(n^{\Nu-1/2})$ (since $y_n=o(a_n)$). Now, to control the maximum, we compare $\arg(t-S_j)$ with $\arg(t-s)$. For $\Delta^{\DT}[s,\gamma_n]\in \Theta_n^{(2)}$, for $j\in\cro{0,\gamma_n}$, 
\begin{equation}
\label{eq:y1}
|S_j-s|\leq \frac{\Cbis{\T}}{m_f}\gamma_nn^{-1/2}+y_n=O_{\cD'[a]}(n^{\Nu-1/2}+y_n)=O_{\cD'[a]}(n^{\Nu-1/2})
\end{equation}
(since $y_n=o(a_n)$) and then \begin{equation}
\label{eq:y2}
|t-S_j|\geq n^{\B-1/2}-|s-S_j|\geq n^{\B-1/2}/2
\end{equation} 
for $n$ large enough since $\Nu<\B$ (uniformly in $s\in{\cD}(a)$). Using that $|e^{ia}-1|=O(|a|)$ we get
 \[\max_{j\in\cro{0,\gamma_n}}\l|e^{i(\arg(t-S_{j-1}(\lambda))-\arg(t-s))}-1\r|=O_{\cD'[a]}(y_n/n^{\B-1/2}).\]
(The tangent of the angle $\widehat{S_{j-1},t,s}$ is $O(y_n/|t-S_j|)$). Then, 
\[d=O_{\cD'[a]}(n^{\Nu\Nu'+\Nu/2-\B+\Nu-1/2})=O_{\cD'[a]}(n^{\Nu\Nu'+3\Nu/2-\B-1/2}). \]
Set  $\widehat{\bZ}_{s}^{(n)}(j/\sqrt{n})=s+\sum_{l=1}^{j} \widehat{\Delta}(s,t,l)$. Again, for any Borelian set $\Theta$, 
\begin{equation}
\label{eq:compf}
|`P_{nf(s)}(\widehat{\bZ}_{s}^{(n)}\in\Theta)-`P_{nf(s)}(\bar{\bZ}_{s}^{(n)}\in\Theta)|\leq  \exp(-\gamma''n^{\Nu\Nu'})+\exp(-n^{\C/2}).
\end{equation}
We have  
\begin{equation}
\label{eq:zzd}
\sup_{x\in[0,a_n]}|\widehat{\bZ}_{s}^{(n)}(x)-e^{i\arg(t-s)}\bZ^{(n)}(x)|\leq d=O_{\cD'[a]}(n^{\Nu\Nu'+3\Nu/2-\B-1/2}).
\end{equation}
Hence, using \eref{eq:rot1bis} and \eref{eq:rot2bis}, if  $\Delta^{\DT}[s,\gamma_n]\in \Theta_n^{(2)}$, then
\[\l\{\l|\frac{\widehat{\bZ}_{s_0}^{(n)}(a_n)-s_0}{a_n} -e^{i\arg(t-s)}F_{\Cbis{\T},0}(\widehat{\bZ}_{s_0}^{(n)}(0))\r|\leq y_n/a_n+d/a_n\r\},\]
and
\[\l\{\sup_{x\in[0,a_n]}\l| (\widehat{\bZ}_{s_0}^{(n)}(x)-\widehat{\bZ}_{s_0}^{(n)}(0))-(x-0)e^{i\arg(t-s)} F_{\Cbis{\T},0}(\widehat{\bZ}_{s_0}^{(n)}(0))\r|\leq y_n+d\r\}.\]
By \eref{eq:law}, this occurs with a probability exponentially close to 1 under $`P_{nf(s)}$, and then this is also true for $\bar{\bZ}^{(n)}$ under $`P_{nf(s)}$  by~\eref{eq:compf}, and then for $\bar{\bZ}^{(n)}$ under $`P_{nf}$ by~\eref{eq:borne7}. \par
We are now in situation to use Corollary~\ref{cor:corn} on the process  $\bar{\bZ}^{(n)}$ under $`P_{nf}$. The corresponding value of $\max\{a_n,c_n,c_n'\}$ is $M_n:=\max\{n^{\Nu-1/2}, (y_n+O_{\cD'[a]}(n^{\Nu\Nu'+3\Nu/2-\B-1/2}))/a_n\}=\max\{n^{\Nu-1/2},
O_{\cD'[a]}(n^{\Nu\Nu'-\Nu/2}),O_{\cD'[a]}(n^{\Nu\Nu'+\Nu/2-B})\}$ and the probability $b_n$ and $d_n$ are smaller than $e^{-n^c}$ for some $c>0$, for $n$ large enough. We want $M_n$ to be as small as possible. For this we choose $\Nu'$ close to 0, and since $\B<1/4$, the maximum is obtained by taking $\Nu=\B^-$. In this case, $\max\{a_n,c_n,c_n'\}=n^{-\B/2^-}$. Now, the conclusion of the theorem holds if $-\B/2^-<-\alpha, \B-1/2<-\alpha$. Hence, $\alpha$ must be chosen in $(0,1/8)$ for the existence of $\B, `e', \Nu$ satisfying all the requirements of the present proof and then at any time less or equal to $\Time_{s,t}^{\Cbis{\X}}$ as long as $\Pos_{s,t}^{\X}(x\sqrt{n})$ is outside $B(t,n^{\B-1/2})$. By what is said above, at time $\Time_{s,t}^{\Cbis{\X}}$, the traveller is in $B(t,n^{-\alpha})$ with probability $1-\exp(-n^c)$; therefore, by Proposition~\eref{pro:boule} after this time, the traveller will come closer to $\{t\}$ at each stage. This implies the first assertion of the theorem. \par
 We now pass to the proof of the second assertion: it remains to show that the number of steps of the traveller to reach the target once in $B(t,n^{-\alpha})$ is negligible before $\sqrt{n}$. From what is said below \eref{eq:nb12}, we only need to control the number of stages needed to enter in $B(t,n^{\B-1/2})$  from time $\Time_{s,t}^{\Cbis{\X}}$, where we know that the traveller is in $B(t,n^{-\alpha})$.
The argument below Equation \eref{eq:nb12} will allow us to see that at most  $dn^{-\alpha+1/2}$ steps will be needed for a constant $d$, with probability exponentially close to 1.  For this, we observe that  Lemma~\ref{lem:c2} can still be used, as well as the sets $\Theta_n^{(2)}$ define above. The sum of the length of the increments after $c_n$ steps (once in $B(t,n^{-\alpha})$ with $c_n\geq \gamma_n$ is at least $dc_nn^{-1/2}\pm \frac{c_n}{\gamma_n}y_n$ for a constant $d$ (that is $dc_nn^{-1/2}$ at the first order) with probability exponentially close to 1. Hence, by Proposition~\ref{pro:boule} (4) the number of steps needed to traverse a distance at most $n^{-\alpha}$ is at most $d n^{-\alpha+1/2}$ for a constant $d$.  ~$\Box$

\subsubsection{Cross navigation}
\label{sssec:Surg-cross}
\begin{theo}\label{theo:clef3}
Let $\X \in \{\Y,\T\}$ and $\theta\leq \pi/3$. For any $\alpha\in(0,1/8)$, there exists $c=c(\alpha)$, 
\[`P_{nf}\l(\sup_{x\in[0,\tau_{s,t}^\X]} \l|\Pos_{s,t}^{\X}(x\sqrt{n}) - \Pos_{s,t}^{\infty,\Cbis{\X},\Cbor{\X}}(x)\r|\geq n^{-\alpha}\r)=O_{\cD'[a]}\l(\exp(-n^{c})\r),\]
where $\tau_{s,t}^\X=\Time_{s,I(s,t)}^{\Cbis{\X}}+\Time_{I(s,t),t}^{\Cbor{\X}}$. Moreover, there exists  $d>0$ such that for $n$ large enough
\[`P_{nf}\l(\l|\Nb^\X_{s,t}/\sqrt{n}-\tau_{s,t}^\X\r|\geq dn^{-\alpha}\r)=O_{\cD'[a]}\l(\exp(-n^{c})\r).\]
\end{theo}
Again, for $\X\in\{\T,\Y\}$, $\tau_{s,t}^\X$ is implicitly known since

\be \Time_{s,I(s,t)}^{\Cbis{\X}}&:=&\inf\l\{~x,\int_{0}^x \frac{e^{i\arg(I(s,t)-s)}\Cbis{\X}\,du}{\sqrt{f(s+e^{i\arg(I(s,t)-s)u})}}=I(s,t)-s\r\}\\
\Time_{I(s,t),t}^{\Cbor{\X}}&:=&\inf\l\{~x,\int_{0}^x \frac{e^{i\arg(t-I(s,t))}\Cbor{\X}\,du}{\sqrt{f(I(s,t)+e^{i\arg(t-I(s,t))u})}}=t-I(s,t)\r\}.
\ee
If $f=c$ is constant, this gives
\[\tau_{s,t}^\X=\sqrt{c}\l(\frac{|I(s,t)-s|}{\Cbis{X}}+\frac{|t-I(s,t)|}{\Cbor{X}}\r).\]

\prooff{Theorem~\ref{theo:clef3}} We here consider the case $\X=\T$, the case $\X=\Y$ being similar.\par 
Recall  the contents of Lemma~\ref{lem:c3}, and of the last two paragraph of Sections~\ref{ssec:local-with-directed}. The proof starts as for the proof of Theorem~\ref{theo:clef2} using also Proposition~\ref{pro:boule}, and again we consider only the case where $|s-t|>n^{\B-1/2}$ and the portion of the traveller trajectory which is outside the ball $B(t, n^{\B-1/2})$ as in the proof of Theorem~\ref{theo:clef2}. Consider, the set $\Omega_{n,`e}$ (defined in~Section~\ref{ssec:fb}).
When $t$ is far from $\cross(t)$, that is if $d(t,\cross(s_0))\geq n^{\Nu+\C-1/2}$, by~\eref{eq:rep-T-1} the \T\ coincides exactly with \DT\ with direction $\theta K_{s,t}$ for $n^\Nu$ stages, provided that $S\in \Omega_{n,`e}$ (since only $n^\Nu$ stages of size at most $\NAVMAX[\theta]$ are concerned). Therefore, all the properties and bounds obtained for this case, particularly in the proof of Theorem~\ref{theo:clef} holds true here also. Hence, the traveller will stay close to ${\bf B}_{K_{s,t}}$ the bisecting line of $\ov{K_{s,t}}{s}$, its fluctuation around this line being larger than $n^{-\alpha}$ for $\alpha\in(0,1/8)$ with probability exponentially small. Therefore, the trajectory of the traveller from its starting point till a neighbourhood of $I_{s,t}$ will be the solution of $\ODE(\Cbis{\T},\arg(I_{s,t}-s),s_0)$.\par

Assume now that the traveller satisfies $d(t,\cross(s_0))\leq  n^{\Nu+\C-1/2}$ (this can occur at the beginning of its trip or after some sequences of $n^\Nu$ consecutive stages). Recall the considerations of the last paragraph of Section~\ref{ssec:local-with-directed} and also observe Fig.~\ref{fig:proj}. If for some $j$, $d(t,\cross(s_j))\leq  n^{\Nu+\C-1/2}$, then this will remain true till the traveller enters in the final ball $B(t, n^{\B-1/2})$.  Let us describe more precisely, what happens when $d(t,\cross(s_j))$ is small. Denote by $\bD$ the parallel to the branch of $\cross(s_i)$ being close to $t$, passing via $t$.
In order to control the position of the traveller, knowing that it is close to $\bD$, an orthogonal projection on $\bD$ is used. The progression of this projection on $\bD$ measures the progress of the traveller toward $t$. We will not enter into the details, everything works with respect to the orthogonal projection on $\bD$ as in the case of directed navigation, since the projection $\xi$ owns also some exponential moments. 
Therefore, once close to $\Delta$, the movement of the traveller will asymptotically be ruled by $\ODE(\Cbor{\T},\arg(t-I_{s,t}),s_0)$.
Now, taking into account that $\Nu+\C-1/2\leq \B-1/2$, then necessarily the traveller will enter in the ball $B(t, n^{\B-1/2})$ where it will make a negligible number of steps, with negligible fluctuations (see discussion below Equation~\eref{eq:nb12}).~$\Box$

\subsection{Analysis of the traveller costs}
\label{ssec:aco}
In Section~\ref{ssec:QI}, we introduced the cost $\Cost_H^{\Nav}(s,t):=\sum_{j=1}^{\Nb(s,t)} H(\Delta_j^{\Nav})$ related to the traveller journey from $s$ to $t$, associated with an elementary cost function $H$, and $\Nav$. 
If $H$ is the modulus function $H:x\mapsto |x|$, then  $\Cost_H^{\Nav}(s,t)=|\Path^{\Nav}(s,t)|$, if $H:x\mapsto 1$ then $\Cost_H^{\Nav}(s,t)=\Nb(s,t)$, already discussed in Theorems~\ref{theo:clef2} and~\ref{theo:clef3}. In the sequel be only consider the case $H_g:x\mapsto |x|^g$ for some $g\geq 0$. Other functions could certainly be studied following the same steps, but the present case covers the applications we have (discussions around the cases $g\in[2,4]$ appear in \cite{LWW01}).\par
Let us discuss a bit at the intuitive level. Under $`P_{nf}$, locally around position $s_0$, a stage $\Delta^{\DT}$ is close in distribution to $\frac{\Delta_1^{\DT}}{\sqrt{nf(s_0)}}$ a rescaled stage under $`P_1$. For a regular function $H$, $H(\Delta^{\DT})$ is close to  $H(\Delta_1^{\DT}/\sqrt{nf(s_0)})$. Two main points have to be noticed so far. The other point is that at the first order, under $`P_{nf}$, $H(\Delta^{\DT})$ depends on the behaviour of $H$ near 0. Functions $H$ ``that are regular near 0'' are needed to get simple asymptotic behaviours. This justifies the choice of the class of functions $H_g$. The contribution to the cost of $H(\Delta^{\DT})$ depends on the position of the traveller. Again a differential equation appears: it is important to consider the pair 
\[\l((\Pos^{\Nav}(s,t,i),{\sf cost}^{\Nav}(s,t,i)),i=1,\dots,k\r),\]
where ${\sf cost}^{\Nav}(s,i)=\sum_{j=1}^i H(\Delta_j^{\Nav})$ since the cost can not be studied independently of the position.\par

Another remark concerns the case $H:x\mapsto 1$, in which case ${\sf cost}^{\Nav}(s,j)=j$. In this case, defining
\[\CC^{(n)}(j/\sqrt{n}):={\sf cost(s,j)}/\sqrt{n}=j/\sqrt{n},~~~~~\textrm{ for }j\geq 0,\]
and by linear interpolation between the points $(j/\sqrt{n},j\geq 0)$, on any interval $[0,\lambda]$, $\l(\CC^{(n)},n\geq 0\r)$ converges uniformly to $c_{\sol}:=y\mapsto y$. Therefore, the pair $(\CC^{(n)},\Pos_{s,t}^{\Nav}(.\sqrt{n}))$ converges to $(c_{\sol},\Pos_{s,t}^{\infty,c})$ under the same conditions under those $\Pos_{s,t}^{\Nav}(.\sqrt{n})$ converges to $\Pos_{s,t}^{\infty,c}$. Even if at the first glance this convergence seems to entail that of $\Nb^{\Nav}(s,t)$ this is not immediately the case as observed in the proofs of Theorem~\ref{theo:clef2} and~\ref{theo:clef3}, since the convergence of functions does not imply the convergence of hitting times. Here the same phenomena arise for other cost functions. \par
Let us now come back to the case where $H=H_g$, for some $g>0$.  By the scaling argument, we have for any $c>0$, 
$H_g(\Delta_c^\X)\sur{=}{(d)}{|\Delta_1^\X|^g}/{c^{g/2}},$
and then, further, under $`P_c$,
\[\sum_{j=1}^{a_n}H_g(\Delta_c^\X(j))\sur{=}{(d)}\sum_{j=1}^{a_n}{|\Delta_1^\X(j)|^g}/{c^{g/2}}.\]
Since the r.v. $|\Delta_1^\X(j)|$ have exponential moments (see Section~\ref{sssec:DLN}) so do $|\Delta_1^\X(j)|^g$ for any $g>0$. Therefore by the law of large numbers for any $a_n\to+\infty$,
\[\sum_{j=1}^{a_n}{|\Delta_1^\X(j)|^g}/{a_n}\as `E(|\Delta_1^\X(j)|^g)\]
and Petrov's Lemma allowing to control the deviation around the mean can be used (see Section \ref{sssec:DLN}). Hence, at the first order, $\sum_{j=1}^{a_n}H_g(|\Delta_c^\X(j)|)$ is close to $\frac{a_n}{c^{g/2}} `E(|\Delta_1^\X(j)|^g)$, and this for any $c$, including the case where $c=c_n$ depends on $n$.
Under $`P_{nf}$, if the traveller is at position $s$, then $|\Delta_{nf}^\X(1)|^g$ is close in distribution to $|\Delta_1^\X|^g/(nf(s))^{g/2}$, since $\Nb(s,t)$ has order $\sqrt{n}$,  set
\[\bC^{(n),\X}(j/\sqrt{n}):={\sf cost^\X(s,j)}n^{g/2-1/2}.\]
In order to extend Theorem~\ref{theo:clef} to the pair (Position,Cost), we need to introduce a system of ODE. We already saw that at least up to some decompositions by parts of the trajectories, the limiting position of the traveller was the solution of $\ODE(\lambda,\nu,s)$ for some parameters $(\lambda,\nu)$ depending on the details of the studied navigation. For any $s_0\in \cD$, $c_0 \in \mathbb{R}$, $q\in\mathbb{R}$ consider the following system 
\bq \label{eq:m}
\ODE^{(2)}(\lambda,\nu,q,s_0,c_0):=
\left\{\begin{array}{ccl} \rho(0)&=&s_0,~~ {\mathcal{C}}(0)=c_0,\\
    \dis\frac{\partial \rho(x)}{\partial x}&=&F_{\lambda,\nu}(\rho(x))=\dis\frac{\lambda e^{i\nu}}{f(\rho(x))^{1/2}},\\
\dis\frac{\partial {\mathcal{C}}(x)}{\partial x}&=& \dis\frac{q}{f(\rho(x))^{g/2}}.
\end{array}\right.
\eq
The existence and uniqueness of a solution $(\rho,\mathcal{C})$ to this system is guarantied by Cauchy-Lipschitz Theorem. The function $\rho$ is $\SOL_{s_0}^{\lambda,\nu}$ since the conditions on $\rho$ coincide with $\ODE(\lambda,\nu,s_0)$.
The function $\mathcal{C}$ has clearly a simple integral representation using $\rho$ and $f$:
\be
{\mathcal{C}}(y)&=& c_0+q  \int_{0}^y \frac{dx}{f(\rho(x))^{g/2}}.
\ee
Since 
\be
\rho(y)&=&s_0+\int_{0}^y F_{\lambda,\nu}(\rho(x)) \,dx
=s_0+ \lambda e^{i\nu} \int_{0}^y \frac{dx}{f(\rho(x))^{1/2}}
\ee
in the case $g=1$ this immediately leads to
\bq
\mathcal{C}(y)=c_0+q\lambda^{-1} (\rho(y)-s_0)e^{-i\nu}.
\eq
If $g\neq 1$, $\rho$ and $\mathcal{C}$ are related by a linear formula only if $f$ is constant since in this case  $\rho$ and $\mathcal{C}$ are linear functions.

\paragraph{Limits in the directed case}

We explain in this case only the appearance of the limiting differential equation. 
\begin{lem}  Let $\X=\DT$ and $\theta\leq \pi$, or $\X=\DY$ and $\theta\leq \pi/2$. The pair of processes $(\Pos_{s,t}^{\X}(x\sqrt{n}),\bC^{(n),\X}(x))_{x\in[0,\lambda]}$ satisfies the assumption of Corollary~\ref{cor:corn}, for $b_n=n^{\Nu\Nu'-\Nu/2}/\sqrt{m_f}$, $c_n=\exp(-n^{\min(\Nu\Nu'/2,c)})$,$c_n'=c_n$, $d_n=b_n$. Therefore it converges to $(\Pos_{s,t}^{\infty,\Cbis{\X}}(x),{\sf cost}^{\infty,\X}(x))_{x\in[0,\lambda]}$  solution  $\ODE^{(2)}(\Cbis{\X},0,\Cbis{\X,g},s_0,0)$  for $\lambda \in[0,\lambda(F_{\Cbis{\X},0},s))$, where $\Cbis{\X,g}=`E(|\Delta_1^\X|^g)$.
\end{lem}

\proof The proof uses the ideas of the proof of Theorem~\ref{theo:clef} (we will use below $\Theta_n^{(1)}$ and $\bZ^{(n)}$, defined in its proof). Consider again the set $\Theta_n^{(1)}$ as defined in~\eref{eq:Theta1}, and introduce the following Borelian subset of $\mathbb{C}^{\gamma_n}$:
\begin{equation}
\label{eq:Theta4}
\Theta_n^{(4)}:=\l\{ (\lambda_1,\dots,\lambda_{\gamma_n}) \in \mathbb{C}^{\gamma_n},~\sup_{l\in\cro{1,\gamma_n}}\l|\sum_{j=1}^{l} \l( n^{g/2-1/2}|\lambda_j|^g -\frac{\Cbis{T,g}}{\sqrt{n}(f(s_0))^{g/2}}\r)\r|\leq y_n \r\}
\end{equation}
for $(y_n)$ a sequence. Using $\dis\bC^{(n)}(a_n)=c_0+ n^{g/2-1/2}\sum_{l=1}^{\gamma_n}|\Delta^\DT(s_0,l)|^g$ we get 
\[\{\Delta^\DT[s_0,\gamma_n]\in \Theta_n^{(4)}\}=\l\{\sup_{l\in\cro{1,\gamma_n}}\l| \bC_{s_0}^{(n)}(l/\sqrt{n})-c_0 -\frac{l}{\sqrt{n}}\frac{\Cbis{T,g}}{f(s_0)^{g/2}}\r|\leq y_n\r\},\]
and everything works as in the proof of Theorem~\ref{theo:clef}, in particular,
using also the rescaling~\eref{eq:CE}, 
\be
`P_{nf(s_0)}\l(\Delta^\DT[s_0,\gamma_n]\in \complement\Theta_n^{(4)}\r)
&=&`P_{nf(s_0)}\l(\sup_{l\in\cro{1,\gamma_n}}\l|\sum_{j=1}^{l}  n^{\frac{g-1}{2}}|\Delta^\DT(s_0,j)|^g -\frac{\Cbis{T,g}}{\sqrt{n}f(s_0)^{\frac{g}2}}\r|\geq  y_n\r)\\
&\leq&`P_{1}\l(\sup_{l\in\cro{1,\gamma_n}}\l|\sum_{j=1}^{l}\frac{|\Delta^\DT(s_0,j)|^g -\Cbis{T,g}}{\sqrt{\gamma_n}}\r|\geq y_n  f(s_0)^{\frac{g}2}\frac{\sqrt{n}}{\sqrt{\gamma_n}}\r).
\ee
Take $y_n=x_{\gamma_n}\frac{\sqrt{n}}{\sqrt{\gamma_n}m_f^{g/2}}$, for $x_n\sim n^{\Nu'}$ and then, the proof ends as the one of Theorem~\ref{theo:clef}. Notice that the deviations of $\bC^{(n)}$ are of the same order as that of $\bZ^{(n)}$ in the proof of Theorem~\ref{theo:clef}. Now, using that $\Delta^\DT \in \Theta_n^{(1)}\cap \Theta_n^{(4)}$ with a probability exponentially close to 1 and therefore the conclusion of the theorem holds true for the pair $(\bZ^{(n)}, \bC^{(n)})$. ~\cq

\paragraph{Limits in the straight case}

Here $\X\in\{\ST,\SY\}$. Consider $(\SOL_{s}^{\Cbis{\X},\arg{t-s}},C^{\X,g})$ the solution of $\ODE^{(2)}(\Cbis{\X},\arg(t-s),\Cbis{\X,g},s,0)$, and consider $\Time_{s,t}^{\Cbis{\X}}$.
The limiting cost will be $\bC^{\X,g}_{s,t}:=C^{\X,g}\l(\Time_{s,t}^{\Cbis{\X}}\r).$  Notice that if $g=1$, then $C^{\X,1}(x)=`E(|\Delta_1^\X|)|\Pos_{s,t}^{\Cbis{\X}}-s|$.
\begin{theo}
Let $\X=\SY$  and $\theta < \pi/2$ or $\X=\ST$ and $\theta \leq \pi/2$.
 For any $\alpha\in(0,1/8)$, any $\beta>0$, for any $\lambda>0$, for $n$ large enough
\[`P_{nf}\l(\sup_{(s,t)\in \cD'[a]}  
\l|\frac{\Cost_{H_g}^{\Nav}(s,t)}{n^{1/2-g/2}}-\bC^{\X,g}_{s,t}\r|\geq n^{-\alpha}\r)\leq n^{-\beta}.\] 
\end{theo}
The proof follows the step of that of Theorem~\ref{theo:clef2}: first a proof for $(s,t)$ fixed is obtained, then the proof is extended to a sub-grid of $\cD[a]^2$ and then to all pairs using the arguments of Section \ref{ssec:GB}. \par
The proof for $(s,t)$ fixed is similar to that of Theorem~\ref{theo:clef2}: the contribution to the cost of the stages of the traveller outside the final ball $B(t,n^{-\alpha})$ is provided by the solution of a differential equation. Then, when the traveller enter in the final balls $B(t,n^{-\alpha})$ and then in $B(t,n^{\B-1/2})$, we use again that these final contributions are negligible and affect the total cost up to a negligible amount (smaller than $n^{-\alpha}$ with a huge probability).  %\par

\paragraph{Limits in the cross case}
Here $\X\in\{\T,\Y\}$. This time again, one must use the decomposition of the limiting path at the point $I(s,t)$.
From $s$ to $I(s,t)$, denote by $(\SOL_{s}^{\Cbis{\X},\arg{I(s,t)-s}},C_1^{\X,g})$ the solution of $\ODE^{(2)}(\Cbis{\X},\arg(I(s,t)-s),\Cbis{X,g},s_0,0)$, on the time interval $[0,\Time_{s,I(s,t)}^{\Cbis{\X}}]$; between $I(s,t)$ and $t$, let $(\SOL_{I(s,t)}^{\Cbor{\X},\arg{t-I(s,t)}},C_2^{\X,g})$ be the solution of $\ODE^{(2)}(\Cbor{\X},\arg(t-I(s,t)),$ $\Cbis{X,g},I(s,t),0)$  on the time interval  $[0,\Time_{I(s,t),y}^{\Cbor{\X}}]$. 
The limiting cost will be in this case 
\[\bC^{\X,g}_{s,t}:=C_1^{\X,g}\l(\Time_{s,I(s,t)}^{\Cbis{\X}}\r)+C_2^{\X,g}\l(\Time_{I(s,t),t}^{\Cbor{\X}}\r).\]
 
\begin{theo} For any $\theta\leq \pi/3$,for any $\X\in\{\T,\Y\}$, any $\alpha\in(0,1/8)$, any $\beta>0$, for any $\lambda>0$, for $n$ large enough
\[`P_{nf}\l(\sup_{(s,t)\in \cD'[a]}  
\l|\frac{\Cost_{H_g}^{\Nav}(s,t)}{n^{1/2-g/2}}-\bC^{\X,g}_{s,t}\r|\geq n^{-\alpha}\r)\leq n^{-\beta}.\] 
\end{theo}

\subsection{Globalisation of the bounds}
\label{ssec:GB}

We have obtained some bounds for the position and for some cost functions of a traveller going from $s$ to $t$. We here desire to prove some uniform bounds since in the main theorems a supremum on $(s,t)$ lies inside the considered probabilities. The number of possible pairs $(s,t)$ being infinite, the union bound is not sufficient here. We adopt a two points strategy to get the uniformity needed. First we get the uniformity for pairs $(s,t)$ where $s$ and $t$ belong to a sub-grid of $\cD[a]$
\[\Grid_n(c_0,a):=n^{-c_0}\mathbb{Z}^2\cap \cD[a].\] 
The cardinality of $\Grid_n(c_0,a)$ being $O(n^{2c_0})$, by the union bound, any theorem of the form:
\[\textrm{for all }(s,t)\in \cD'[a], ~~`P_{nf}\l(h(\Path_{s,t}^{\Nav})\in A\r)=O( u_n)\] 
where $u_n$ does not depend on $s,t$, and where $h$ is any function of the paths  has the following corollary
\[`P_{nf}\l(\sup_{(s,t)\in \Grid_n(c_0,a)} h(\Path_{s,t}^{\Nav})\in A\r)=O( n^{2c_0}u_n).\]
In other words, if the probability concerning one path is exponentially small (this is the case for  most of our theorems concerning one trajectory), then it 
is still exponentially small when considering altogether all paths starting and ending 
in $\Grid_n(c_0,a)$.\par

The second point of our strategy is the following. Take any $(s,t) \in \cD'[a]$ say in the 
complementary of $\Grid_n(c_0,a)$.
We will show that 
with a probability close to 1, the trajectory from $s$ to $t$ can be split with a huge probability in at most 13 parts  
(uniformly on $(s,t)$), such that on each part the 
path of the traveller coincides with a part of a path of a traveller starting 
and ending on the grid. This will be sufficient to conclude, since the theorems we have control 
the behaviour of the path of a traveller all along, and then on the aforementioned parts. \par
In order to do so, we introduce $\Squ^*_{n}(c_0)$ the set of squares of the plane with side 
length $n^{-c_0}$, having their vertices in $n^{-c_0}\mathbb{Z}^2$ and at least one of them in
 $\Grid_{n}(c_0,a)$. \par
Additionally, consider four tilings of the plane with squares with length
\[a_n:=n^{\B-1/2},\]
the three last ones being obtained from the first one by the translation of $a_n/2,$ $i a_n/2$ 
and  $a_n(1+i)/2$, respectively. By $\Squ_a^1(\B), \Squ_a^2(\B), \Squ_a^3(\B), \Squ_a^4(\B)$ we denote the subsets 
of the squares of each of these tilings having a distance to $\cD[a]$ smaller than $a_n$ and by 
$\Squ_a(\B)$ their union. If  $n$ is large enough, the union of the squares of $\Squ_a(\B)$ contains 
$\cD[a]$ and are included in the interior of $\cD$ (for example in $\cD[a/2]$), and observe also 
that any disk $B(x,a_n/4)$ with $x\in\cD[a]$ is totally included in a square of $\Squ_a(\B)$. \par

With any $(s,t)\in\cD'[a]$ we will associate two points $(s_g,t_g)$ belonging to $\Grid_{n}(c_0,a)$ 
as follows. First  $t_g$ is the upper left corner of the square of $\Squ^*_{n}(c_0)$ containing $t$ 
(therefore $t_g$ depends on $n$ and $c_0$); if $t$ is in $\Squ^\star(c_0)$ then take $t_g=t$. The point $s_g$ is given by the following lemma.

\begin{lem}
\label{lem:rmin} Let $s\in \cD$, and consider one of the navigation processes presented in this paper using as set of stopping places $S$. For any $\rho < r_{\min}(S):=\min\{|x-y|, x,y \in S, x\neq y\}$ and for any pair  $(s', t)$ of $S\times\cD$ there exists a point $s_g$ in $\rho\mathbb{Z}^2\cap \cD$ such that $\Nav(s_g, t)=\Nav(s,t)=s'$, in other words the path from $s$ and from $s_g$ to $t$ coincides from the second stage.
\end{lem}
\proof First observe that around $s'$ there is a disk of radius $r_{\min}(S)$ containing no other vertex of $S$.
For a straight navigation process let us define $R=\dovc(s',s'-t)(r_{\min}(S)/2)$ and for a cross navigation process, let us define $R=\ovc{s}{-k}(r_{\min}(S)/2)$ where $k$ is such that $t \in \sect[s',k]$.
\begin{figure}[ht]
 \psfrag{s}{$'$}
 \psfrag{R}{$R$}
 \psfrag{t}{$t$}
 \psfrag{rmin}{$r_{\min}(S)$}
 \psfrag{demirmin}{$r_{\min}(S)/2$}
 \psfrag{demitheta}{$\theta/2$}
 \psfrag{r}{$r$}
\centerline{\includegraphics[width=8cm]{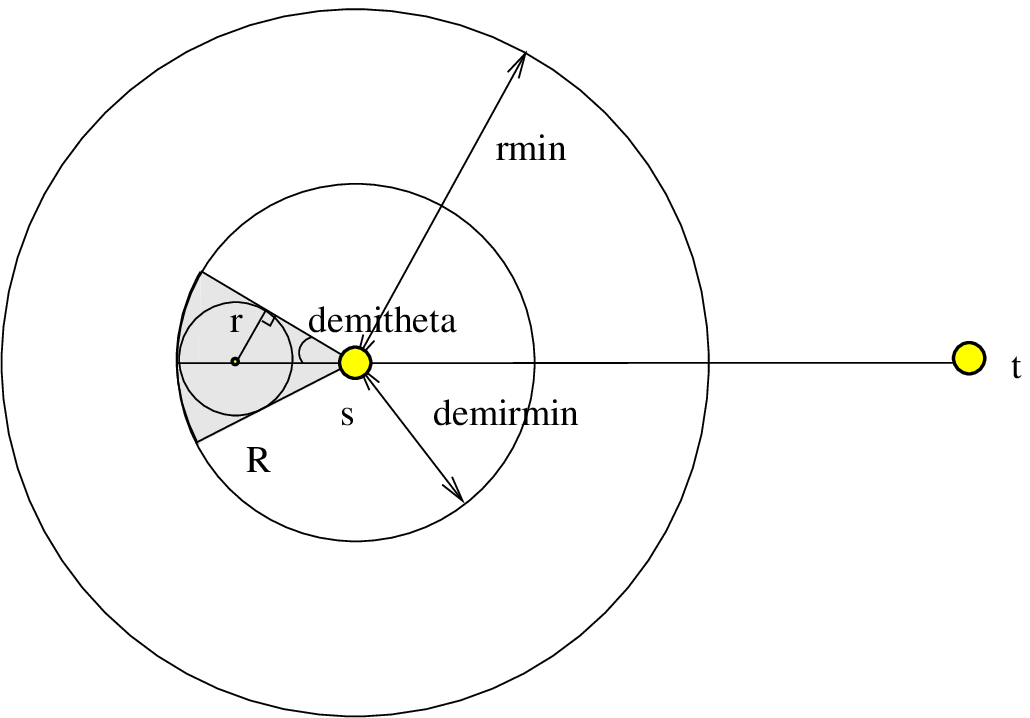}}
\captionn{\label{fig:rmin} Illustration of the proof of Lemma~\ref{lem:rmin}}
\end{figure}
 By construction any point $s$ of $R$ is such that $\Nav(s, t)=s'$. $R$ contains a disk of radius $r=\frac{r_{\min}(S)}{2}\frac{\sin (\theta/2)}{1+\sin (\theta/2)} > \frac{\theta}{4\pi}$ (see Fig.~\ref{fig:rmin}). This disk intersects at least one point of  $\rho\mathbb{Z}^2\cap \cD$.~$\Box$
\begin{lem} For any $c>0$, 
$`P_{nf}(r_{min}\leq n^{-c})=O(n^{-2c+2}).$
\end{lem}
\proof With a $`P_{nf}$ probability exponentially close to 1, $|\#\bS|$ is smaller than $c_1n$ for some $c_1>0$ (see the proof of Lemma~\ref{lem:max-boule}). By the union bound
\[`P_{nf}(r_{min}\leq \lambda ~|~|\#\bS|=m)\leq m^2 `P_{nf}(|V_1-V_2|\leq \lambda)\]
for $V_1,V_2$ independent with density $f$ on $\cD$. Conditioning on $V_1$, it is easily seen that this probability is smaller than $c_2 \lambda^2$ (for some $c_2$ independent on $\lambda$). 
Finally, we get $`P_{nf}(r_{min}\leq \lambda)\leq c_1^2n^2c_2\lambda^2+`P_{nf}(|\#\bS|\geq c_1n)$, which leads easily to the stated result.~$\Box$ \medskip

A consequence of the two previous lemmas is that if $c_0$ is large enough, with a probability $O(n^{-2c_0+2})$ any path from $s$ to $t$ coincides with a path from a point $s_g$ to $t$ for a point $s_g$ in $\Grid_n(c_0,a)$ up to the starting position (and then, up to the first step, smaller than $\NAVMAX$, which will then be uniformly negligible at the scale we are working in). Hence, to approximate a path from $s$ to $t$ with a path from $s_g$ to $t_g$, most of the difficulties come from the target.\medskip

Given a source $s$, a \emph{section} of the $\Path(s,t)=(s_0,\dots, s_{\Nb(s,t)})$ is a sub-path $(s_{j_1},\dots, s_{j_2})$ for which $\Nav(s, t, j)=\Nav(s_{j-1},t_g)$ for each $1\leq j_1 < j  < j_2$, that is a part of $\Path(s,t)$ matching a path toward $t_g$. Given a source $s$ and a target $t$, a $\Squ_a(\B)$-\emph{black box} is a sub-path $(s_{j_1},\dots, s_{j_2})$ of $\Path(s,t)$ such that $s_{j_1}$ and $s_{j_2}$  are both in the same square of $\Squ_a(\B)$. 
For $k\in\mathbb{N}$, we will say that a path admits a $(k,\B,c_0)$ decomposition if it can be decomposed into at most $k$ sections and $k$ $\Squ_a(\B)$ black boxes. 
\begin{pro}
\label{pro:globalization}
Let $\B>0$ and $\rho>0$ be fixed. For any $\X\in\{\Y,\T\}$ and $\theta\leq \pi/3$, or $\X=\ST$ and $\theta\leq \pi/2$, or $\X=\SY$ and $\theta <  \pi/2$, there exists $c_0>0$ such that, 
$$`P_{nf}\l(\forall (s,t)\in\cD'[a],  \Path^{\X}(s,t) \text{admits a  $(13,\B,c_0)$ decomposition}\r) \geq  1-o(n^{-\rho}).$$
\end{pro}
Before proving this proposition, let us examine how it entails, together with the already proved results concerning one trajectory, the theorems of this paper.

\subsection{Consequences}

From the above discussion, the results concerning one trajectory as Theorems~\ref{theo:clef2} and~\ref{theo:clef3} can be extended as follows: a supremum on $(s,t)$ in $\Grid(c_0,a)$ can be added inside the probabilities present in these theorems. It remains to see what is the price to take the supremum on all $(s,t)$. We saw above that with a polynomially close to 1 probability, for any $(s,t)$, $\Path(s,t)$ can be decomposed in at most $13$ sections, and $13$ black boxes. On each section, $\Path(s,t)$ coincides with a section of $\Path(s,t_g)$. Therefore the fluctuations of the position functions (or the cost functions) on each of these sections  are smaller than $n^{-\alpha}$ with a probability exponentially close to 1, since this is the case for the trajectories between points of $\Grid_n(c_0,a)$. To see this, assume for example that 
\[\sup_{x\in[0,\lambda]} \l|\Pos^{\X}_{s,t}(x\sqrt{n})-\Pos^{\infty,\Cbis{\X}}_{s,t}(x)\r| \leq `e.\] Therefore, for any increasing sequence $x_1,\dots,x_{2k}$ (where $x_{2i-1}\sqrt{n}$ and $x_{2i}\sqrt{n}$ have to be understood as the beginning and ending time of the sections), one has
\[\sup_{0\leq x_1\leq \dots \leq x_{2k}\leq \lambda}\l|\sum_{i=1}^k\l|\Pos^{\X}_{s,t}(x_{2i}\sqrt{n})-\Pos^{\X}_{s,t}(x_{2i-1}\sqrt{n})\r|-\sum_{i=1}^k\l|\Pos^{\infty,\Cbis{\X}}_{s,t}(x_{2i})-\Pos^{\infty,\Cbis{\X}}_{s,t}(x_{2i-1})\r|\r| \leq k`e;\] Hence, a global control of the path ensures a good control of the sections, provided that the number of sections is small. This is the case here, since we have 13 sections at most with a large probability.  Now, the contribution of the black boxes (the stages between times $x_{2i}\sqrt{n}$ and $x_{2i+1}\sqrt{n}$) have to be controlled. For this, notice first that the space fluctuations for any black box $(s_{j_1},\dots,s_{j_2})$ is small: $\max\{|s_j-s_{j_1}|, j\in\cro{j_1,j_2}\}\leq c a_n$ for a constant $c$ depending only on $\theta$; indeed, since $s_{j_1}$ and $s_{j_2}$ are in the same square $\Box$ of $\Squ_a(\B)$, and since at each stage the traveller gets closer to the target, and since the angle between its local direction and the direction to the target is  smaller than $\theta$, it must stay in a domain with area at most $O(a_n^2)$ (which is included in the ball $B(\circ,c{a_n})$ where $\circ$ is the centre of  $\Box$).\begin{figure}[ht]
 \psfrag{ct}{$C_t$}
 \psfrag{C}{$C$}
 \psfrag{theta}{$\theta$}
 \psfrag{t}{$t$}
\centerline{\includegraphics[width=5cm]{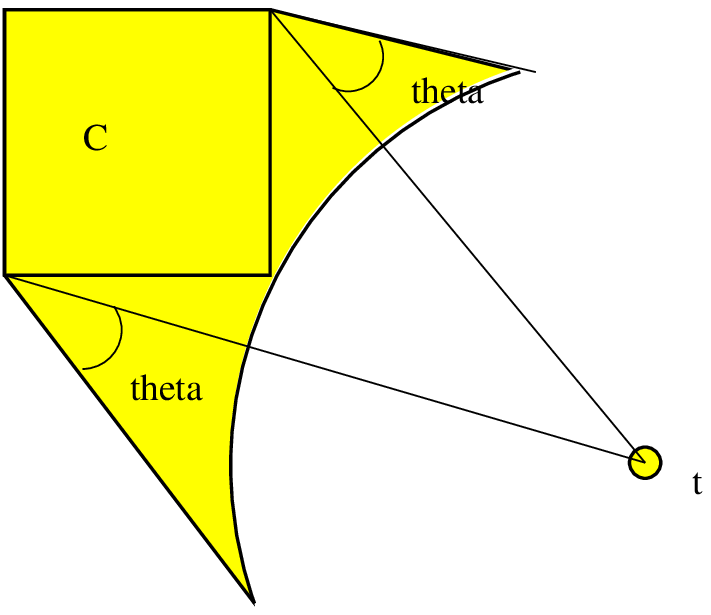}}
\vskip -2cm ~

\captionn{\label{fig:blackbox} Area around $\Box$ containing the range of the black-box $(s_{j_1},\dots,s_{j_2})$}
\end{figure}
Moreover, let $l(\X)$ be the maximum length (number of stages) of all black boxes for the algorithm  $\X\in\{\T,\Y,\ST,\SY\}$:
\[l(\X):=\max\{j_2-j_1, (s_{j_1},\dots,s_{j_2}) \textrm{ black-box}\}\]
is then bounded by $\MAXBALL(ca_n,a)$, and then 
\[`P(l(\X)\geq a_n^2n^{1+`e})=`P(l(\X)\geq n^{2\B+`e})\leq \exp(-n^d)\]
for some $d>0$, if $n$ is large enough.
Since $2\B+`e$ can be chosen smaller than $1/2$ in all proofs, at most $k$ black box, of negligible size $n^{2\B+`e}$ are concerned. The contribution of these black boxes to the time normalisation of the processes is negligible, as well as the space normalisation: at most $a_n^2n^{1+`e}$ contributions of size $n^\C$ with size smaller than $\NAVMAX$ or $H(\NAVMAX/\sqrt{n})$ are (uniformly for all $(s,t)$) negligible.

\subsection{Proof of Proposition~\ref{pro:globalization}}

In order to bound the number of black boxes on a trajectory, we need to understand when a 
navigation decision differs when navigating to $t$ or to $t_g$. We will see that this is 
related to some geometrical features of some straight lines, called below ``rays'', issued from the points of $S$.\par
\begin{figure}[ht]
 \psfrag{s}{$s$}
 \psfrag{s1}{$s_1$}
 \psfrag{s2}{$s_2$}
 \psfrag{s3}{$s_3$}
 \psfrag{s4}{$s_4$}
 \psfrag{s5}{$s_5$}
 \psfrag{s6}{$s_6$}
 \psfrag{r1}{$r_1$}
 \psfrag{r2}{$r_2$}
 \psfrag{r3}{$r_3$}
 \psfrag{r4}{$r_4$}
 \psfrag{r5}{$r_5$}
 \psfrag{r6}{$r_6$}
 \psfrag{r7}{$r_7$}
 \psfrag{r8}{$r_8$}
 \psfrag{r10}{$r_{10}$}
 \psfrag{theta}{$\theta$}
\centerline{\includegraphics[width=14cm]{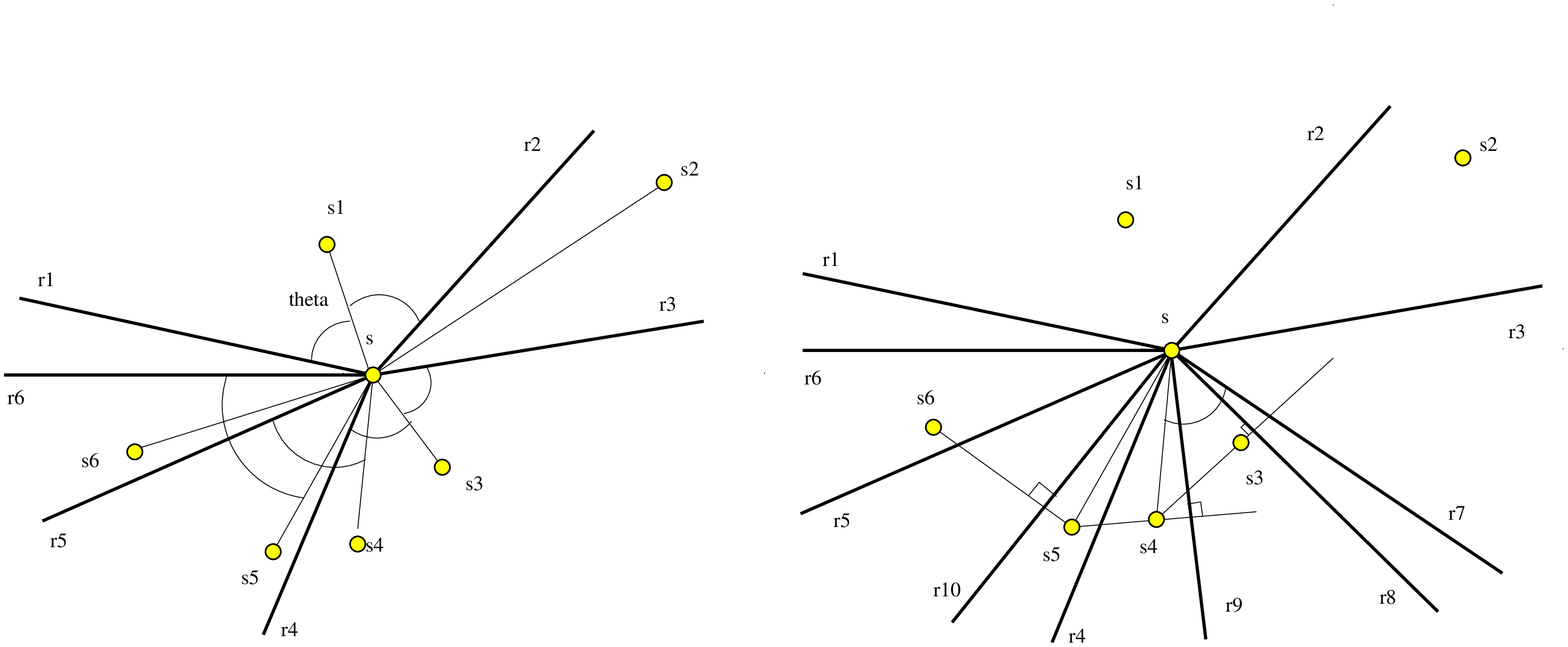}}
\captionn{\label{fig:star} On the left an example of a navigation-star in \SYN. For a vertex $t$ for enough between $r_i$ and $r_{i+1}$, $\SY(s,t)=s_i$ (assuming that $r_7=r_1$). On the right, an example of a navigation-star in \STN. In this case the situation is a bit more complex. For instance, $\ST(s,t)=s_3$ if $t$ is between $r_3$ and $r_{7}$ or between $r_8$ and $r_4$. If $t$ is between $r_7$ and $r_8$ (and far enough) then $\ST(s,t)=s_4$.}
\end{figure}

 A \emph{star} is a 
collection of half-lines (we also use the term \emph{rays}) issued from the same point called the \emph{centre}. 
A \emph{navigation-star} of $s \in S$ is a cyclic ordered list $(r_1, r_2, \dots, r_k)$ of half-lines starting at $s$ such that 
for any two points $t,t'$ of $\cD$ between two consecutive half-lines then if $\Nav(s,t)\neq t$ and $ \Nav(s,t')\neq t'$ then
$\Nav(s,t)=\Nav(s,t')$; in other words, if $t$ and $t'$ are far enough from $s$ (at distance at least $2\NAVMAX$ suffices) 
then the first stop starting from $s$ is the same whether if the target point is $t$ or $t'$ (see Fig.\ref{fig:star}).

In order to build navigation stars we now define different types of rays; they will be used to control the decompositions of straight or cross navigation paths in the sequel:
\begin{itemize}
\item[--] Given a point $s$ of $S$ and $k \in\cro{1,p_{\theta}}$, the $k$th \emph{type-1 ray} around $s$ is simply $r_1(s,k)=\HL_k(s)$, 
\item[--] Let $s$ be a point of $S$. We construct two rays associated with each points $s'$ of $S$. Denote by $r_2^+(s,s')$ (resp. $r_2^-(s,s')$)  the half-line obtained by a rotation of $\theta/2$ (resp. $-\theta/2$) with centre $s$ of the half-line $[s,s')$. These rays are called \emph{type-2 rays}.
\item[--] Let $s$ be a point of $S$. For any pair $(s',s'')$ of $S^2$, we denote by $r_3(s,s')$  the half-line $[s,t)$ such that $[s,t)$ is orthogonal to $(s',s'')$. This rays are called \emph{type-3 rays}.
\end{itemize}

We now examine the navigation stars of the different algorithms we have. It may help to notice that if $NS$ is a navigation star around $s$ then any set of rays containing $NS$ is also a navigation star.
\begin{itemize}

\item {\bf Cross navigations}. For cross navigations, the set of type-1 rays centred at $s$, in other word, $\cross(s)$, is a navigation star of $s$.

\item {\bf \SYN}. 
Starting from $s$ for a target $t$ far enough, $\SY(s,t)$ depends only on $\arg(t-s)$. When $t$ is turning around $s$, $\SY(s,t)$ changes when the nearest point of $s$ in $\dov{(s,t)}\cap S$ is changing. This potentially happens when the angle between $(s,t)$ and $(s,s')$ is $\pm \theta/2$ for a point $s'$ in $S$ close enough to $s$, in any case for $s'$ such that $|s-s'|\leq\NAVMAX$. Hence the set of type-2 rays centred at $s$ is a navigation star of $s$ for \SYN.

\item {\bf \STN}. Additionally to the rays presented in the previous point, another case has to be treated. Notice that $\ST(s,t)$ is the point of $\dov(s,t)\cap S$ that has the closest to $s$ orthogonal projection on the bisecting line of $\dov(s,t)\cap S$. When $t$ is turning around $s$, $\ST(s,t)$ is potentially changing when the order of the orthogonal projections of the elements of $\dov{(s,t)}\cap S$ on the bisecting line of $\dov(s,t)$ is changing. This may happens when the line that passes via two elements $s'$ and $s''$ of $S$, not too far from $s$ (such that $|s-s'|\leq\NAVMAX$ and $|s-s''|\leq\NAVMAX$), is orthogonal to the line $(s,t)$. Hence the set of type-2 rays and type-3 rays centred at $s$ is a navigation star of $s$ for \STN.

\end{itemize}
A navigation star containing only the rays defined in the previous paragraph is called \emph{standard}. A \emph{constellation} (resp. \emph{standard constellation}) of $S$ is the union of the navigation stars (resp. standard navigation stars) of every $s\in S$.
 \\
Given  $C_t \in \Squ^*_{n}(c_0)$ and $C \in \Squ_a(\B)$, denote by $\Numb_i(C,C_t)$ the number of type-$i$ rays that intersects $C$ and that are defined with points $s$, $s'$ and $s''$ in $S\cap C$ (depending on the case, $s'$ and $s''$ can be non necessary). Define also, 
\ben 
\label{eq:pf1}\Numb^{\T}(C,C_t)&=&\Numb^{\Y}(C,C_t)=\Numb_1(C,C_t),\\\label{eq:pf2}\Numb^{\SY}(C,C_t)&=&\Numb_2(C,C_t),\\
\label{eq:pf3}\Numb^{\ST}(C,C_t)&=&\Numb_2(C,C_t)+\Numb_3(C,C_t)
\een
 as one can guess in view of the above discussion. We also set $\Numb(C,C_t)=\sum_{i=1}^3 \Numb_i(C,C_t)$ a bound on all the $\Numb^{\X}(C,C_t)$ for $\X\in\{\ST,\SY,\T,\Y\}$.
\begin{lem} 
\label{lem:proba-cross}
Let $\B>0$ and $c_0>0$ be fixed. For any algorithm $\X\in\{\ST,\SY,\T,\Y\}$,
\[\sup\limits_{\substack{(C,C_t)\in \Squ_a(\B)\times\Squ^*_{n}(c_0)\\ d(C,C_t)\geq a_n}} `P_{nf}\l(\Numb^{\X}(C,C_t)\geq 1\r) = O\l(n^{4\B+2 -c_0}\r).\]
\end{lem}
\proof All along the proof,  $C_t \in \Squ^*_{n}(c_0)$ and $C \in \Squ_a(\B)$ are such that $d(C_t,C)\geq a_n$. We first control $\Numb_1(C,C_t)$. Let $St(C_t)$ the set of points $s \in \cD[a]$ such that $\cross(s)$ crosses $C_t$. This set forms  also a cross.
Clearly for any $i\in\{1,2,3\}$ 
\be 
`P_{nf}(\Numb_i(C,C_t)\geq 1) &\leq& `P_{nM_f}(\Numb_i(C,C_t)\geq 1).
\ee
We then obtain the bound under $`P_{nM_f}$ instead than $`P_{nf}$. Now, 
\be`P_{nM_f}(\Numb_1(C,C_t)\geq 1) &\leq&`P_{nM_f}\l(\# \l(\bS \cap C \cap St(C_t)\r) \geq 1\r)= O(n |C \cap St(C_t)|)
\ee
where $|A|$ denote the area of $A$ (since, for an integer random variable $X$, $`P(X\geq 1)\leq `E(X)$). Observing that $St(C_t)$ is composed by $p_\theta$ stripes having width $O(n^{-c_0})$ and that $C$ has  diameter $\sqrt{2}n^{\B-1/2})$, we have    $ |C \cap St(C_t)|=O(n^{\B-1/2-c_0})$, and then
\bq`P_{nf}(\Numb_1(C,C_t)\geq 1)= O(n^{1/2+\B-c_0}).
\eq
Let us pass to the control of $\Numb_2$.  For $p$ the centre of $C_t$, we have 
$ `P_{nf}(\Numb_2(C,C_t)\geq 1)$
\[ \leq`P_{nM_f}\l(\exists X_1,X_2 \in \bS \cap C, \l||\arg(X_2-X_1)-\arg(p-X_1)| - \frac\theta2\r| \leq  \arctan \l(\frac{n^{-c_0}}{\sqrt{2}|p-X_1|}\r)\r).\]
Since $|X_1 - p| > n^{\B-1/2}$, $\arctan \l(\frac{n^{-c_0}}{\sqrt{2}|p-X_1|}\r) \leq  \arctan  \l(n^{-c_0-\B+1/2}/\sqrt{2}\r)$. Moreover using~\eref{eq:mb},  $`P_{nM_f}(|\bS \cap C| > n^{2\B+d}) \leq \exp(-n^{2\B+d})$ for any $d>0$, provided that $n$ is large enough. Recall now, that under $`P_{nM_f}$, knowing that $\#(\bS\cap C)=k$, the set $\bS\cap C$ is distributed as a sample of $k$ i.i.d. uniform random variables in $C$. Therefore, for any function $g$ and measurable set $A$, $`P_{nM_f}(\exists X_1,X_2 \in \bS, X_1\neq X_2,  g(X_1,X_2)\in A)$
\be
&\leq &`P_{nM_f}(\#(\bS\cap C >n^{2\B+d}))\\
&&+\sum_{k\leq n^{2b+d}}`P_{nM_f}(\exists X_1,X_2 \in \bS\cap C, g(X_1,X_2)\in A| \#\bS=k) `P_{nM_f}(\#\bS=k)\\
& \leq &`P_{nM_f}(\#(\bS\cap C >n^{2\B+d})) +(n^{2\B+d})^2 `P(g(U_1,U_2)\in A)
\ee
for two i.i.d. uniform random variables $U_1$ and $U_2$ in $C$. 

Therefore, $`P_{nM_f}(\Numb_2(C,C_t)\geq 1)$ is bounded by:
\[ O((n^{2\B+d})^2)`P\l( \l||\arg(U_2-U_1)-\arg(p-U_1)| - \frac\theta2\r| \leq  \arctan \l(\frac{n^{-c_0}}{\sqrt{2}|p-U_1|}\r)\r) + \exp(-n^{2\B+d}).\]
A simple picture shows that $U_1$ given, $U_2$ must lie in a triangle included in $C$, with diameter smaller than $\sqrt{2}a_n$ and  with an angle smaller than $O(\frac{n^{-c_0-\B+1/2}}{\sqrt{2}})$. Since the density of $U_2$ is $1/a_n^2$, for any $U_1$ given in $C$, the probability is $O(n^{-c_0-\B+1/2}a_n/a_n^2)=O(n^{-c_0-2\B+1})$ (with a constant in the $O$ sign uniform for $U_1$ in $C$). Then
\be
 `P_{nf}(\Numb_2(C,C_t)\geq 1)&=&  O(n^{4\B+2d}) O(n^{-c_0-2\B+1}) +  \exp(-n^{2\B+d})
                              =  O(n^{2\B+2d-c_0+1}).
\ee
Taking $d$ small enough, we get $`P_{nf}(\Numb_2(C,C_t)\geq 1) = O(n^{2\B-c_0+2}).$ \par
Let us pass to the control of $\Numb_3$. The quantity $`P_{nf}(\Numb_3(C,C_t)\geq 1)$ is bounded by 
\[ `P_{nM_f}\l( \exists X_1,X_2,X_3 \in \bS \cap C,  \l|\arg(X_3-X_2)-\arg(p-X_1) - \frac{\pi}2\r| \leq  \arctan \l( \frac{n^{-c_0-\B+1/2}}{\sqrt{2}} \r)\r).\]
Again, this is 
\[
O((n^{2\B+d})^3)`P\l( \l|\arg(U_3-U_2)-\arg(p-U_1) - \frac{\pi}2\r| \leq  \arctan \l( n^{-c_0-\B+1/2} \r)\r)\]
for some $U_1,U_2$ and $U_3$ i.i.d. uniform in $C$. For $(U_1,U_2)$ given, $U_3$ must be in a subset of $C$ with Lebesgue measure $O(n^{-c_0-\B+1/2}a_n)$. Taking into account the density of $U_3$, we get
\[`P_{nf}(\Numb_3(C,C_t)\geq 1)=O(n^{4\B-c_0+2})\]
for $d$ chosen small enough.\par
We conclude this proof using~\eref{eq:pf1},~\eref{eq:pf2},~\eref{eq:pf3}, and the union bound.~$\Box$
\medskip

For any $C_t \in \Squ^*_{c_0}$ let$\BB({C_t})= \{ C \in \Squ_a(\B),~\Numb(C,C_t)\geq 1\}$
be the number of squares in $\Squ_a(\B)$ from which is issued at least one ray crossing $C_t$. Notation $\BB$ is chosen
to make apparent that the elements of $\BB(C_t)$ are considered as black boxes related to $C_t$ later on.
\begin{lem}
\label{lem:boxes}
Let $\B>0$ be fixed and $\rho>0$ be fixed. There exists a $c_0>0$ such that,
\[`P_{nf}\l(\exists C_t\in \Squ^*_{c_0}, \#\BB(C_t) \geq 12 \r) = o(n^{-\rho}).\]
\end{lem}
\proof Notice that under $`P_{nM_f}$, since $\bS \cap C_1$ and $\bS\cap C_2$ are independent variables for $C_1\cap C_2=\varnothing$, for any fixed $j\in \{1,2,3,4\}$, any $C_t$ and $i$ fixed, 
the family of variables $(\Numb_i(C,C_t), C\in \Squ_a(\B)^j)$ are i.i.d.. Hence, for each $j\in\{1,2,3,4\}$
the variables $\l(\1_{\sum_{i=1}^3 \Numb_i(C,C_t)\geq 1},C\in\Squ_a(\B)^j\r)$ are i.i.d. Bernoulli distributed. Let  
\[\BB_j({C_t})= \{ C \in \Squ_a(\B)^j,~\Numb(C,C_t)\geq 1\}.\]
We have, by the union bound,
\be
 `P\l(\exists C_t \in \Squ^*_{c_0}, \#\BB_j(C_t)\geq 3\r) &\leq& \sum_{C_t\in\Squ^*_{c_0}} \binom{\Squ_a(\B)^j}{3} `P_{nM_f}(\Numb(C,C_t)\geq 1)^3\\
&=&O(n^{3-6\B+2c_0})\sup_{C\in \Squ_a(\B),C_t\in \Squ^*_{c_0}} `P_{nM_f}(\Numb(C,C_t)\geq 1)^3.
\ee
From Lemma~\ref{lem:proba-cross} this last term is $O(n^{3(4\B+2 -c_0)})$. Then, this is $o(n^{-\rho})$ for $c_0$ chosen large enough. Now, to conclude, write 
\[`P_{nf}\l(\exists C_t\in \Squ^*_{c_0}, \#\BB(C_t) \geq 12 \r)\leq \sum_{j=1}^4  `P_{nM_f}\l(\exists C_t \in \Squ^*_{c_0}, \#\BB_j(C_t)\geq 3\r).~~~\Box ~\]

\prooff{Proposition~\ref{pro:globalization}} First observe that rays are defined by at most 3 
points belonging to a ball of radius at most $\NAVMAX$. From Lemma~\ref{lem:petitspas}  we know 
that $\NAVMAX \leq n^{-1/2+\B}/4$ with high probability. Hence the probability that 12 rays cross a 
square $C_t$ is bounded by
\[`P_{nM_f}\l(\exists C_t\in \Squ^*_{c_0}, \#\BB(C_t) \geq 12 \r)+`P_{nf}(\NAVMAX \geq n^{-1/2+\B}/4),\]
which is $O(n^{-\rho})$ as wished, for $\B>0$ and $c_0$ large enough. This means that for any target $t$, for 
any point $s$ outside these 12 squares plus the square of $\Squ_a(\B)$ that contains $B(t,n^{-1/2+\B})$: 
 $\Nav(s,t)= \Nav(s,t_g)$. For each $\Path(s,t)$ and each $C$ of these (at most) 13 squares we consider the sub-path  $P_C$ (of $\Path(s,t)$) between the first stage that is inside $C$ and the last stage that is inside $C$. Each of these portions of the trajectory forms a black box. 

For every point $s'$ of the path $\Path(s,t)$ outside each of these $13$ black boxes, none of the rays of the 
navigation star of $s'$ crosses the square $C_t$, hence $\Nav(s',t)=\Nav(s',t_g)$. Hence path is partitioned 
into at most 13 black boxes and at most 13 sections. $\Box$

\section{Appendix}

\prooff{Lemma~\ref{lem:determ-version}} By hypothesis, for a certain triangular array $(c_n'(j), j=1,\dots,\floor{\lambda/a_n})$ satisfying  $|c'_n(j)|\leq c_n$ for any $j\leq \floor{\lambda/a_n}$, the following representation of $y_n$ holds:
\[y_n((j+1)a_n)=y_n\l(ja_n\r)+a_nG(y_n(ja_n))+a_nc_n'(j).\]
Hence $y_n$  appears to be a perturbed version of the explicit Euler scheme, used to approximate a solution of $\Eq(G,z)$, which is defined by $z_n(0)=z$,
\[z_n((j+1)a_n)=z_n\l(ja_n\r)+a_nG(z_n(ja_n)),~~~~ j\leq \floor{\lambda/a_n}.\]
We review here the standard argument leading to the conclusion, namely the convergence of $(z_n)$ to $y_{\sol}$ with a bound on the speed of convergence, and derive the same result for $(y_n)$. First, by Cauchy-Lipschitz Lemma, $\Eq(G,z)$ has a unique solution, denoted $y_{\sol}$ on $[0,\lambda]$ (for $\lambda<\lambda(G,z)$ the maximum domain on which one can extend this solution). 
Now define for any $x\in [0,\lambda]$
\begin{equation}
R_{a_n}(x)=a_n^{-1}(y_{\sol}(x+a_n)-y_{\sol}(x))-G(y_{\sol}(x)).
\end{equation}
Assume that $G$ is $\alpha$-Lipschitz and bounded by $\beta$ (it is bounded since it is the maximum of a Lipschitz function on the bounded domain $ D$); then $y_{\sol}$ is clearly $\beta$-Lipschitz. We have
\ben\label{eq:ineq-R}
\sup_{x \in [0,\lambda-a_n]}\l|R_{a_n}(x)\r|&=&\sup_{x \in [0,\lambda-a_n]}\l|\int_{x}^{x+a_n} \frac{G(y_{\sol}(u))du}{a_n}-G(y_{\sol}(x))\r|\\
&\leq & \sup_{x \in [0,\lambda-a_n]}\alpha\int_{x}^{x+a_n} \frac{\l|y_{\sol}(u)-y_{\sol}(x)\r| du}{a_n}
\leq  \frac{\alpha\beta}{2}a_n. 
\een
Let $\zeta_n$ be a perturbation of the explicit Euler scheme:
\[\zeta_n((j+1)a_n)=\zeta_n\l(ja_n\r)+a_nG(\zeta_n(ja_n))+`e_n(j),~~~~ j\leq \floor{\lambda/a_n}.\]
We have
\be
\zeta_n((j+1)a_n)-z_n((j+1)a_n)&=&\zeta_n(ja_n)-z_n(ja_n)+a_n\l(G(\zeta_n(ja_n))-G(z_n(a_n)\r)+`e_n(j)\\
|\zeta_n((j+1)a_n)-z_n((j+1)a_n)|&\leq&(1+\alpha a_n)|\zeta_n(ja_n)-z_n(ja_n)|+|`e_n(j)|
\ee
and therefore, for any $j$
\be
|\zeta_n(ja_n)-z_n(ja_n)|&\leq&(1+\alpha a_n)^j|\zeta_n(0)-z_n(0)|+\sum_{m=0}^{j-1}(1+\alpha a_n)^{j-1-m}|`e_n(m)|
\ee
which leads (using $(1+\alpha a_n)^k \leq e^{ka_n}$ for any $k\geq 0$) that for $j$ such that $j\leq \floor{\lambda/a_n}+1$, 
\be
|\zeta_n(ja_n)-z_n(ja_n)|&\leq &e^{ja_n}\l(|\zeta_n(0)-z_n(0)|+\sum_{m=0}^{j-1}|`e_n(m)|\r)\\
&\leq&e^{a_n(1+\floor{\lambda/a_n})}\l(|\zeta_n(0)-z_n(0)|+\sum_{m=0}^{\floor{\lambda/a_n}}|`e_n(m)|\r).
\ee
Now to end, observe that
\[y_{\sol}((j+1)a_n)=y_{\sol}(ja_n)+a_n G(y_{\sol}(ja_n))+a_nR_{a_n}(ja_n).\]
This is a perturbed version of the Euler explicit scheme, and then by~\eref{eq:ineq-R}, for $n$ large enough, 
\begin{equation}
\label{eq:contri1}
\sup_{j,~j\leq \floor{\lambda/a_n}+1}\l|y_{\sol}(ja_n)-z_n(ja_n)\r|\leq e^{a_n(1+\floor{\lambda/a_n})}\sum_{m=0}^{\floor{\lambda/a_n}} a_n^2\alpha\beta/2\leq \lambda e^\lambda\alpha\beta a_n,
\end{equation}
since $e^{a_n}<2$ for $n$ large enough.\par
The same work, using the comparison between $z_n$ and $y_n$ leads, for $n$ large enough, to
\begin{equation}
\label{eq:contri2}
\sup_{j,~ \floor{\lambda/a_n}+1}|y_n(ja_n)-z_n(ja_n)|\leq e^{a_n(1+\floor{\lambda/a_n})}\sum_{m=0}^{\floor{\lambda/a_n}} a_nc'_n(j)\leq 2e^{\lambda}\lambda\max_{j\leq \lambda/a_n}|c'_n(j)|\leq  2\lambda e^{\lambda} c_n.
\end{equation}
Finally summing~\eref{eq:contri1} and~\eref{eq:contri2} we get  $\sup_{x\in[0,\lambda]}|y_{\sol}(x)-y_n(x)|\leq C_\lambda\max\{a_n, c_n\},$
for $C_\lambda=\lambda e^{\lambda}(2+\alpha\beta)+3\beta$ (this is obtained at first for the points $x\in\{ja_n,~ j\leq \floor{\lambda/a_n}+1\}$, then for all $x$ in the interval using that on an interval of size $a_n$ the fluctuation of $y_{\sol}$ are bounded by $a_n \beta$,  and those of $y_n$ are bounded by $\beta a_n + c_n$. Then, the map $\lambda\mapsto C_\lambda$ has the properties stated in the lemma. 
\cq

\prooff{Corollary~\ref{cor:corn}}  Let $(\bZ_n)$ be a sequence of processes satisfying the hypothesis. Define, for each $n$, a process $\bY_n$ coinciding with $\bZ_n$ at the points $(ja_n, j=0,\dots,\floor{\lambda/a_n}+1)$, and which is linear between them. Introduce the event
 \ben\label{eq:Yn}
\Omega_n&:=&\bigcap_{j=0}^{\floor{\lambda/a_n}}\l\{\l|\frac{\bY_n\l((j+1)a_n\r)-\bY_n\l(ja_n\r)}{a_n}-G(\bY_n(ja_n))\r|\leq c_n\r\}
\een 
according to Lemma~\ref{lem:determ-version}, $\sup_{x\in [0,\lambda]}|\bY_n(x)-y_{\sol}(x)| \leq D_\lambda\max\{a_n,c_n\}$ on $\Omega_n$, for a function $\lambda\mapsto D_\lambda$ bounded on compact sets; moreover by the union bound, 
\begin{equation}
`P(\Omega_n)\geq 1-(\floor{\lambda/a_n}+1)a_n d_n.
\end{equation}
We then also get immediately $`P\l(\sup_{x\in [0,\lambda]}|\bY_n(x)-y_{\sol}(x)| \leq `e\r)\sous{\to}n 1$, for any $`e>0$, that is the convergence in probability for the uniform norm of $\bY_n$ to $y_{\sol}$. Define now 
\[\Omega'_n:=\bigcap_{j=0}^{\floor{\lambda/a_n}}\l\{\sup_{x\in [ja_n,(j+1)a_n]} \l|\bZ_n\l(x\r)-\bZ_n\l({ja_n}\r)-(x-ja_n)G(\bZ_n(ja_n))\r|\leq c_n'\r\}.\]
Again, by the union bound $`P(\Omega_n')\geq 1-(\floor{\lambda/a_n}+1)a_n b_n$. Assume now that we are on $\Omega_n\cap \Omega_n'$. Since $\bY_n$ is linear between in $[ja_n,(j+1){a_n}]$ , and since $\bZ_n$ and $\bY_n$ coincide at the points $(ja_n,j\geq 0)$, then $\bZ_n$ satisfies also~\eref{eq:Yn}. Hence, for any $t\in [ja_n,(j+1){a_n}]$, 
\[\bY_n(x)=\bZ_n(ja_n)+(x-ja_n)G(\bZ_n(a_n)).\]
And thus, on $\Omega_n'$, we have
$|\bZ_n(x)-\bY_n(x)|\leq c_n',~~~~ \textrm{ for any }x\in[0,\lambda]$
and then, $|\bZ_n(x)-y_{\sol}(x)|\leq c_n'+D_\lambda\max\{a_n,c_n\}$ on $\Omega_n\cap \Omega_n'$. Therefore
\be
`P\l(\sup_{x\in[0,\lambda]} |\bZ_n(x)-y_{\sol}(x)|\leq c_n'+D_\lambda\max\{a_n,c_n\}\r)&\geq& 1-(\floor{\lambda/a_n}+1)a_n (b_n+d_n)\\
&\geq&  1-(\lambda+a_n) (b_n+d_n).~~~\textrm{\cq}
\ee 

\small

\bibliographystyle{hplain}
\bibliography{biblio}

\end{document}